\documentclass{article}

\usepackage{amsmath}
\usepackage{amssymb}
\usepackage{amsthm}
\usepackage[all]{xy}

\newtheorem{thm}{Theorem}[section]
\newtheorem{prop}[thm]{Proposition}
\newtheorem{lem}[thm]{Lemma}
\newtheorem{cor}[thm]{Corollary}
\newtheorem{ithm}{Theorem}

\theoremstyle{definition}
\newtheorem{dfn}[thm]{Definition}

\theoremstyle{remark}
\newtheorem{rem}{Remark}
\newtheorem*{acknowledgments}{Acknowledgments}

\newcommand{\C}{\mathbb{C}}
\newcommand{\R}{\mathbb{R}}
\newcommand{\T}{\mathbb{T}}
\newcommand{\Z}{\mathbb{Z}}

\newcommand{\g}{\mathfrak{g}}
\newcommand{\ha}{\mathfrak{h}}
\newcommand{\su}{\mathfrak{su}}
\renewcommand{\AA}{\mathfrak{A}}

\newcommand{\A}{\mathcal{A}}
\newcommand{\F}{\mathcal{F}}
\newcommand{\G}{\mathcal{G}}
\newcommand{\sS}{\mathcal{S}}
\newcommand{\U}{\mathcal{U}}
\newcommand{\V}{\mathcal{V}}
\newcommand{\bF}{\bar{\F}}

\newcommand{\Ad}{\mathrm{Ad}}
\newcommand{\id}{\mathrm{id}}
\newcommand{\Hom}{\mathrm{Hom}}
\newcommand{\Lie}{\mathrm{Lie}}
\newcommand{\Tr}{\mathrm{Tr}}

\newcommand{\Cech}{\v{C}ech {}}
\newcommand{\im}{\sqrt{\! - \! 1}}
\renewcommand{\d}{\partial}

\newcommand{\hGamma}{\h{\Gamma}}
\renewcommand{\l}{\ell}
\renewcommand{\epsilon}{\varepsilon}
\def\u#1{ \underline{#1} }
\def\h#1{ \widehat{#1} }
\def\til#1{ \tilde{#1} }

\title{Reduction of strongly equivariant bundle gerbes with 
connection and curving}

\author{Kiyonori Gomi
\thanks{The author's research is supported by Research Fellowship of the Japan Society for the Promotion of Science for Young Scientists.}}

\date{}

\begin{document}

\maketitle

\begin{abstract}
From a certain strongly equivariant bundle gerbe with connection and curving over a smooth manifold on which a Lie group acts, we construct under some conditions a bundle gerbe with connection and curving over the quotient space. In general, the construction requires a choice, and we can consequently obtain distinct stable isomorphism classes of bundle gerbes with connection and curving over the quotient space. A bundle gerbe naturally arising in Chern-Simons theory provides an example of the reduction.
\end{abstract}

\tableofcontents


\section{Introduction}
\label{sec:introduction}

By the theory of Kostant \cite{Ko} and Weil \cite{We}, the degree two cohomology classes of a smooth manifold are geometrically realized by Hermitian line bundles or by circle bundles through their characteristic classes. Similarly, the degree three cohomology classes admit various geometric realizations. The notion of \textit{gerbes}, invented by Giraud \cite{Gi}, is one of such realizations. Brylinski investigated theory of gerbes in detail, and provided a number of applications to geometry and topology \cite{Bry1}. In particular, connective structures and curvings, which are notions of connections on gerbes, are due to his work.

\medskip

In \cite{Go2}, the author studied, under certain conditions, the relationship between \textit{equivariant gerbes} \cite{Bry2,Bry1} with connective structure and curving over a smooth manifold $M$, on which a Lie group $G$ acts, and gerbes with connective structure and curving over the quotient space $M/G$. For convenience, we include here the result in a concise fashion.

\begin{ithm}[\cite{Go2}] \label{ithm:coh_reduction_gerbe}
Suppose that the action of $G$ on $M$ is free and locally trivial, and that we can make the quotient space $M/G$ into a smooth manifold in such a way that the projection map $q : M \to M/G$ is smooth.

(a) There exists a gerbe with connective structure and curving over $M/G$ whose pull-back under $q : M \to M/G$ is isomorphic to a $G$-equivariant gerbe with connective structure and curving over $M$ given, if and only if an obstruction class vanishes.

(b) For a $G$-equivariant gerbe with connective structure and curving over $M$ whose obstruction class vanishes, the isomorphism classes of such gerbes with connective structure and curving over $M/G$ as in (a) are, in general, not unique.
\end{ithm}

We notice that an obstruction similar to Theorem \ref{ithm:coh_reduction_gerbe} (a) exists for a $G$-equivariant principal $\T$-bundle with connection, where $\T = \{ u \in \C |\ |u| = 1 \}$ is the unit circle. In fact, the obstruction is the \textit{moment} \cite{B-V} associated with a $G$-equivariant $\T$-bundle with connection. However, there is no counterpart of Theorem \ref{ithm:coh_reduction_gerbe} (b) in considering equivariant $\T$-bundles: there is exactly one isomorphism class of a $\T$-bundle over $M/G$ whose pull-back is equivariantly isomorphic to a given $G$-equivariant $\T$-bundle with connection whose moment vanishes.

The proof of Theorem \ref{ithm:coh_reduction_gerbe} is rather abstract: we use certain cohomology groups developed on the basis of Brylinski's work \cite{Bry2}. So it is left as a problem to construct concretely a gerbe with connective structure and curving over $M/G$ from a $G$-equivariant gerbe with connective structure and curving over $M$.

\medskip

The purpose of the present paper is to answer this problem by using so-called \textit{strongly equivariant bundle gerbes} \cite{Ma-S,Me}. The notion of \textit{bundle gerbes} was introduced by Murray \cite{Mu}, as the other geometric realization of a degree three cohomology class of a smooth manifold. In fact, the stable isomorphism classes of bundle gerbes over $M$ are classified by $H^3(M, \Z)$. As in the case of gerbes, bundle gerbes admit differential geometric structures called connection and curving. In \cite{Ma-S}, Mathai and Stevenson defined an equivariant version of bundle gerbes. As they remarked, their definition is too strong in a sense. Later, in the preprint of \cite{Me}, Meinrenken proposed a general formulation of equivariant bundle gerbes. The equivariant bundle gerbes of Mathai and Stevenson are recovered in a special case, and are called \textit{strongly equivariant bundle gerbes}.

\medskip

Now we explain the results in this paper: a construction of a bundle gerbe with connection and curving over $M/G$ from a strongly $G$-equivariant bundle gerbe with connection and curving over $M$. In the following, we call such a construction a \textit{reduction} (or \textit{quotient}) of a strongly equivariant bundle gerbe with connection and curving.

\medskip

First of all, we recall the definition of bundle gerbes, following \cite{Me}. A bundle gerbe $\G = (Y, P, s)$ over $M$ consists of a surjective submersion $\pi : Y \to M$, a $\T$-bundle $P \to Y^{[2]}$ and a section $s : Y^{[3]} \to \delta P$ such that $\delta s = 1$. (See the figure below.) Here we denote by $Y^{[p]}$ the $p$-fold fiber product of $Y$, and define $\delta P$ as $\delta P = \pi_1^*P \otimes \pi_2^*P^{\otimes -1} \otimes \pi_3^*P$, using the projection $\pi_i : Y^{[3]} \to Y^{[2]}$ omitting the $i$th factor. The $\T$-bundle $\delta \delta P$ and the section $\delta s$ are defined in the same way. Note that $\delta \delta P$ is canonically isomorphic to the trivial $\T$-bundle over $Y^{[4]}$.  A connection on $\G = (Y, P, s)$ is defined to be a connection 1-form $\nabla \in \im A^1(P)$ on the circle bundle $P$ such that $s^*(\delta \nabla) = 0$. A curving for a connection $\nabla$ is defined to be a 2-form $f \in \im A^2(Y)$ such that $\delta f = F(\nabla)$, where $F(\nabla)$ stands for the curvature of $\nabla$.

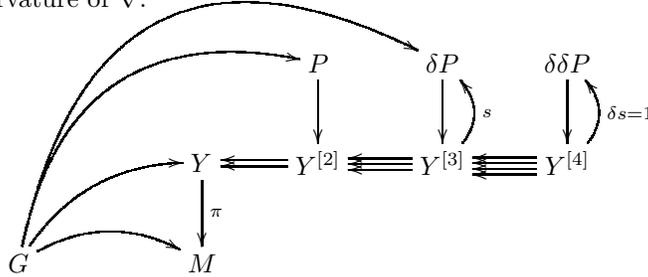
$$
\xymatrix{
&  &
&
P \ar[d] &
\delta P \ar[d] &
\delta \delta P \ar[d] &
\\
&  &
Y^{   } \ar[d]^{\pi} &
Y^{[2]} \ar[l]<-0.3ex> \ar@<0.3ex>[l] & 
Y^{[3]} \ar@<0.1ex>[l] \ar@<-0.4ex>[l] \ar@<0.6ex>[l] 
        \ar@/_1pc/[u]_{s}&
Y^{[4]} \ar[l] \ar@<-0.5ex>[l] \ar@<0.5ex>[l] \ar@<1ex>[l]
        \ar@/_1pc/[u]_{ \delta s = 1} \\
G \ar@/^5pc/[uurrrr] \ar@/^3pc/[uurrr] \ar@/^1pc/[urr]  \ar@/^1pc/[rr] &  &
M &
&
&
}
$$

Let a Lie group $G$ act on $M$ by left. A bundle gerbe $\G = (Y, P, s)$ is said to be strongly $G$-equivariant, when the action of $G$ on $M$ lifts to $Y$, the induced action on $Y^{[2]}$ also lifts to $P$ by bundle isomorphisms, and the section $s$ is $G$-invariant. A connection $\nabla$ and a curving $f$ are said to be $G$-invariant, when they are $G$-invariant differential forms.

\medskip

Next we describe the obstruction class to the reduction. Let $\G = (Y, P, s)$ be a strongly $G$-equivariant bundle gerbe equipped with a $G$-invariant connection $\nabla$ and a $G$-invariant curving $f$. We define a map $\til{\lambda} : Y^{[2]} \to \Hom(\g, \im\R)$ by $\langle X | \til{\lambda}(y_1, y_2) \rangle = \nabla(p; X^*)$, where $p \in P$ is a point lying on the fiber of $(y_1, y_2) \in Y^{[2]}$, $X^* \in T_pP$ is the tangent vector generated by the infinitesimal action of $X \in \g = \Lie G$, and $\langle \ | \ \rangle : \g \otimes \Hom(\g, \im\R) \to \im\R$ is the natural contraction. Note that the value $\til{\lambda}(y_1, y_2)$ is independent of the choice of the point $p$ on the fiber. The map $\til{\lambda}$ is the moment (\cite{B-V}) associated with $(P, \nabla)$. We introduce the following vector spaces:
\begin{align*}
\mathcal{C} &= 
A^1(M, \g^*) \oplus C^{\infty}(G, A^0(M, \g^*)), \\
\mathcal{Z} &= 
\left\{
(E, \zeta) \in \mathcal{C} \Big|
\begin{array}{l}
g^*E - \Ad_gE = d \zeta(g), \\
\Ad_g\zeta(h) - \zeta(gh) + h^*\zeta(g) = 0. \\
\end{array}
\right\}, \\
\mathcal{B} & =
\left\{
(E, \zeta) \in \mathcal{C} \Bigg|
\begin{array}{l}
\mu \in A^0(M, \g^*), \\
E = d \mu, \\
\zeta(g) = g^*\mu - \Ad_g \mu.
\end{array}
\right\}. 
\end{align*}
It is known \cite{Ma-S} that there is a map $\lambda : Y \to \Hom(\g, \im\R)$ such that $\delta \lambda = \til{\lambda}$. When we choose such a map $\lambda$, we define $(E, \zeta) \in \mathcal{Z}$ by setting
\begin{align}
\langle X | \pi^*E \rangle 
&= 
\frac{1}{2\pi\im}
\left( \langle X | d \lambda \rangle + \iota_{X^*} f \right), 
\label{iformula:E} \\
\langle X | \pi^*\zeta(g) \rangle
&=
\frac{1}{2\pi\im}
\langle X | g^*\lambda - \Ad_g \lambda \rangle.
\label{iformula:zeta}
\end{align}
Though $(E, \zeta) \in \mathcal{Z}$ above depends on the choice of $\lambda$, the element $[E, \zeta] \in \mathcal{Z}/\mathcal{B}$ is independent of the choice. We define $\beta_G(\G, \nabla, f) \in \mathcal{Z}/\mathcal{B}$ by $\beta_G(\G, \nabla, f) = [E, \zeta]$. The next theorem implies that $\beta_G(\G, \nabla, f)$ is the obstruction to the reduction of $(\G, \nabla, f)$: 

\begin{ithm} \label{ithm:obstruction}
Let $G$ and $M$ be as in Theorem \ref{ithm:coh_reduction_gerbe}. There exists a bundle gerbe with connection and curving over $M/G$ whose pull-back under $q : M \to M/G$ is stably isomorphic to a given strongly $G$-equivariant bundle gerbe with connection and curving $(\G, \nabla, f)$ over $M$, if and only if $\beta_G(\G, \nabla, f) = 0$ in $\mathcal{Z}/\mathcal{B}$.
\end{ithm}

We remark that the obstruction $\beta_G(\G, \nabla, f) \in \mathcal{Z}/\mathcal{B}$ originates in the study of an equivariant version of the smooth Deligne cohomology which classifies the isomorphism classes of equivariant gerbes with connective structure and curving \cite{Go1}. We also remark that the expression of $\beta_G(\G, \nabla, f)$ owes its simplicity to that of strongly equivariant bundle gerbes. 

\medskip

Now we perform a reduction of a strongly $G$-equivariant bundle gerbe with connection and curving $(\G, \nabla, f)$. We assume that the action of $G$ on $M$ is as in Theorem \ref{ithm:coh_reduction_gerbe}. We also put the same assumption on the action of $G$ on $Y$. 

In this case, a simple construction of a bundle gerbe $\bar{\G} = (\bar{Y}, \bar{P}, \bar{s})$ over $M/G$ is known \cite{Ma-S}: we put $\bar{Y} = Y/G$ and $\bar{P} = P/G$. Then the $G$-invariant section $s : Y^{[3]} \to \delta P$ corresponds to $\bar{s} : \bar{Y}^{[3]} \to \delta \bar{P}$. However, $\nabla$ and $f$ do not descend to give a connection and a curving on $\bar{\G}$ directly, since these $G$-invariant differential forms do not necessarily vanish in the direction of $G$.

Let us suppose that the obstruction class $\beta_G(\G, \nabla, f)$ vanishes. Then we can choose a map $\lambda : Y \to \Hom(\g, \im\R)$ such that $\delta \lambda = \til{\lambda}$ and $(E, \zeta) = 0$, where $(E, \zeta) \in \mathcal{Z}$ is defined by (\ref{iformula:E}) and (\ref{iformula:zeta}). By the assumption, $q : M \to M/G$ gives rise to a principal $G$-bundle, on which $G$ acts by left. Let us choose and fix a connection $\Xi \in A^1(M, \g)$ on this bundle. We define a 1-form $\kappa \in \im A^1(Y)$ by $\kappa = \langle \pi^*\Xi | \lambda \rangle$, and consider the strongly $G$-equivariant bundle gerbe with connection and curving $(\G, \nabla - \delta \kappa, f - d \kappa)$. We notice that $(\G, \nabla - \delta \kappa, f - d \kappa)$ and $(\G, \nabla, f)$ are stably isomorphic. Because $\nabla - \delta \kappa$ and $f - d \kappa$ vanish in the direction of $G$, they descend to give a connection $\bar{\nabla}$ and a curving $\bar{f}$ on $\bar{\G}$, respectively. This construction is our reduction of $(\G, \nabla, f)$:

\begin{ithm} \label{ithm:reduction_sebgcc}
The pull-back of $(\bar{\G}, \bar{\nabla}, \bar{f})$ is stably isomorphic to $(\G, \nabla, f)$ as a $G$-equivariant bundle gerbe with connection and curving over $M$.
\end{ithm}

As is apparent, our construction depends on the choice of a map $\lambda$. Let $\lambda'$ be the other choice, and $(\bar{\G}, \bar{\nabla}', \bar{f}')$ the bundle gerbe with connection and curving obtained by using $\lambda'$. In general, $\bar{\nabla}'$ coincides with $\bar{\nabla}$, while $\bar{f}'$ differs from $\bar{f}$.

\begin{ithm} \label{ithm:difference}
$(\bar{\G}, \bar{\nabla}, \bar{f})$ and $(\bar{\G}, \bar{\nabla}, \bar{f}')$ are stably isomorphic, if and only if $\lambda' - \lambda$ is the moment of a $G$-equivariant $\T$-bundle with flat connection.
\end{ithm}

We also shown that, if $\lambda$ is fixed, then the stable isomorphism class of $(\bar{\G}, \bar{\nabla}, \bar{f})$ is independent of the choice of the connection $\Xi$.

\bigskip

One example of our method of reduction is given by the action of $G = \T$ on $M = SU(2)$ through the diagonal embedding $\T \to SU(2)$. In this case, the reduction yields all the stable isomorphism classes of bundle gerbes with connection and curving over $M/G \cong S^2$ by varying the choice of $\lambda$. The other example concerns Chern-Simons theory. We study a bundle gerbe over the space of connections on the trivial $SU(2)$-bundle over $S^1$, on which the gauge transformation group acts naturally. The reduction yields a bundle gerbe with connection and curving whose 3-curvature is the canonical 3-form on $SU(2)$. This example is also related to the \textit{Chern-Simons line bundle} \cite{F1} and to a \textit{quasi-Hamiltonian space} introduced by Alekseev, Malkin and Meinrenken \cite{A-M-M}.

\bigskip

This paper is organized as follows. Section \ref{sec:EDC}--\ref{sec:EBG} are preliminaries. In Section \ref{sec:EDC}, we briefly review the smooth Deligne cohomology group and its equivariant version \cite{Go1,Go2}. In Section \ref{sec:bundle_gerbes}, we recall bundle gerbes and some basic facts related. In Section \ref{sec:EBG}, we introduce Meinrenken's equivariant bundle gerbes, and explain the relation with strongly equivariant bundle gerbes. 

After these preliminaries, we construct, in Section \ref{sec:characteristic_class}, an equivariant smooth Deligne cohomology class for an equivariant bundle gerbe with connection and curving. This is thought of as a characteristic class, and is a key to the study of $\beta_G(\G, \nabla, f)$. In Section \ref{sec:obstruction}, we show Theorem \ref{ithm:obstruction} by using results in Section \ref{sec:EDC} and Section \ref{sec:characteristic_class}. The reduction of strongly equivariant bundle gerbes is given in Section \ref{sec:reduction}. We prove Theorem \ref{ithm:reduction_sebgcc} and \ref{ithm:difference} in the section. A reduction of ``trivializations'' for equivariant bundle gerbes is also given. Finally, in Section \ref{sec:example}, we provide examples of reductions.

\bigskip

\begin{acknowledgments}
I would like to thank T. Kohno and M. Furuta for useful discussions and valuable suggestions. 
\end{acknowledgments}


\section{Review of the smooth Deligne cohomology}
\label{sec:EDC}

In this section, we recall the smooth Deligne cohomology group \cite{Bry1,De-F,E-V} and its equivariant version \cite{Go1,Go2}. 

\subsection{Smooth Deligne cohomology groups}

As a convention of this paper, a ``smooth manifold'' means a paracompact smooth manifold modeled on a Hausdorff locally convex topological vector space. We also assume the existence of a partition of unity. We always work in the smooth category, so we often drop the word ``smooth.''

\medskip

Let $M$ be a smooth manifold. We denote by $\u{\T}_M$ the sheaf of germs of functions on $M$ which take its values in the unit circle $\T = \{ z \in \C |\ |z| = 1 \}$, and by $\u{A}^q_M$ the sheaf of germs of $\R$-valued differential $q$-forms on $M$.

\begin{dfn}
Let $N$ be a non-negative integer.

(a) We define the \textit{smooth Deligne complex} $\F(N)_M$ by the following complex of sheaves on $M$:
$$
\u{\T}_M \stackrel{\frac{1}{2\pi\im}d\log}{\longrightarrow}
\u{A}^1_M \stackrel{d}{\longrightarrow} 
\u{A}^2_M \stackrel{d}{\longrightarrow} 
\cdots \stackrel{d}{\longrightarrow}
\u{A}^N_M \longrightarrow
0 \longrightarrow \cdots.
$$
where the sheaf $\u{\T}_M$ is located at degree 0 in the complex. 

(b) We define the \textit{smooth Deligne cohomology group} $H^p(M, \F(N)_M)$ by the hypercohomology group of the smooth Deligne complex.
\end{dfn}

We often omit the subscripts of $\u{\T}_M, \u{A}^q_M$ and $\F(N)_M$.

\begin{prop}[\cite{Bry1}] \label{prop:Deligne_coh:manifold}
Let $N$ be a positive integer.

(a) If $0 \le p < N$, then $H^p(M, \F(N))$ is isomorphic to $H^p(M, \T)$, 
the cohomology with coefficients in the constant sheaf $\T$.

(b) If $p = N$, then $H^N(M, \F(N))$ fits into the following exact sequences:
\begin{gather*}
0 \longrightarrow
H^N(M, \T) \longrightarrow
H^N(M, \F(N)) \longrightarrow
A^{N+1}(M)_0 \longrightarrow 0, \\
0 \longrightarrow
A^N(M) / A^N(M)_0 \longrightarrow
H^N(M, \F(N)) \longrightarrow
H^{N+1}(M, \Z) \longrightarrow 0, 
\end{gather*}
where $A^q(M)_0$ is the group of closed integral $q$-forms on $M$. 

(c) If $N < p$, then $H^p(M, \F(N))$ is isomorphic to 
$H^p(M, \u{\T}) \cong H^{p+1}(M, \Z)$.
\end{prop}


\subsection{Simplicial manifolds and simplicial sheaves}

We recall simplicial manifolds \cite{Se} and the notion of a sheaf on a simplicial manifold (a simplicial sheaf, for short) \cite{De}. 

A \textit{simplicial manifold} $X_{\bullet}$ is a sequence of manifolds $\{ X_p \}_{p \ge 0}$ together with face maps $\d_i : X_{p+1} \to X_p, \ (i = 0, \ldots, p+1)$ and degeneracy maps $s_i : X_p \to X_{p+1}, \ (i = 0, \ldots, p)$ obeying the following relations:
\begin{align}
\d_i \circ \d_j &= \d_{j-1} \circ \d_i, \quad (i < j), 
\label{rel1} \\
s_i \circ s_j &= s_{j+1} \circ s_i, \quad (i \le j), 
\label{rel2} \\
\d_i \circ s_j &= 
\begin{cases}
s_{j-1} \circ \d_i, & (i < j), \\
\id, & (i = j, j+1), \\
s_j \circ \d_{i-1}, & (i > j+1).
\end{cases}
\label{rel3}
\end{align}

The simplicial manifold that we will use is the one associated to a Lie group action. Let $G$ be a Lie group, which is not necessarily compact. When $G$ acts on a smooth manifold $M$ (by left), we define a simplicial manifold $G^{\bullet} \times M$ by the sequence $\{ G^p \times M \}_{p \ge 0}$. The face map $\d_i : G^{p+1} \times M \to G^p \times M$ is given by
$$
\d_i(g_1, \ldots, g_{p+1}, x) =
\begin{cases}
(g_2, \ldots, g_{p+1}, x), 
& i = 0 \\
(g_1, \ldots, g_{i-1}, g_i g_{i+1}, g_{i+2}, \ldots, g_{p+1}, x), 
& i = 1, \ldots, p \\
(g_1, \ldots, g_p, g_{p+1} x), 
& i = p + 1,
\end{cases}
$$
and the degeneracy map $s_i : G^p \times M \to G^{p+1} \times M$ is given by
$$
s_i(g_1, \ldots, g_p, x) = (g_1, \ldots, g_i, e, g_{i+1}, \ldots, g_p, x).
$$

We note that the \textit{realization} \cite{Se} of $G^{\bullet} \times M$ is the homotopy quotient \cite{A-B} $(EG \times M)/G$, where $EG$ is the total space of the universal bundle for $G$.

\medskip

Let $X_{\bullet} = \{ X_p \}_{p \ge 0}$ be a simplicial manifold, and $\sS_{\bullet} = \{ \sS_p \}_{p \ge 0}$ a sequence of sheaves such that $\sS_p$ is a sheaf on $X_p$. When there are homomorphisms $\til{\d_i} : \d_i^{-1}\sS_p \to \sS_{p+1}$ and $\til{s}_i : s_i^{-1} \sS_{p+1} \to \sS_p$ obeying the same relations as (\ref{rel1}), (\ref{rel2}) and (\ref{rel3}), we call $\sS_{\bullet}$ a \textit{simplicial sheaf} on $X_{\bullet}$. For example, the sequence $\u{\T}_{X_{\bullet}} = \{ \u{\T}_{X_p} \}_{p \ge 0}$ forms a simplicial sheaf by taking natural homomorphisms $\til{\d}_i$ and $\til{s}_i$.

The notion of a \textit{complex of simplicial sheaves} is defined in a similar fashion.

\begin{dfn}[\cite{Bry3,Bry2}]
We define an open cover of a simplicial manifold $X_{\bullet}$ by a sequence of open covers $\U^{\bullet} = \{ \U^{(p)} \}_{p \ge 0}$ such that:

(a) $\U^{(p)} = \{ U^{(p)}_{\alpha^{(p)}} \}_{\alpha^{(p)} \in \AA^{(p)} }$ is an open cover of $X_p$;

(b) the index set $\AA^{(p)}$ forms a simplicial set $\AA^{\bullet} = \{ \AA^{(p)} \}_{p \ge 0}$;

(c)
we have $\d_i( U^{(p+1)}_{\alpha^{(p+1)}} ) \subset U^{(p)}_{\d_i( \alpha^{(p+1)} )}$ and $s_i( U^{(p)}_{\alpha^{(p)}} ) \subset U^{(p+1)}_{s_i( \alpha^{(p)} )}$.
\end{dfn}

We can find in \cite{Bry2,Me} a construction of an open cover of $G^{\bullet} \times M$.

\medskip

Let $\sS_{\bullet} = \{ \sS_p \}_{p \ge 0}$ be a complex of simplicial sheaves, where $\sS_p$ is a complex of sheaves on $X_p$. We denote by $\sS_p^{[k]}$ the sheaf located at degree $k$ in the complex, and by $\til{d} : \sS_p^{[k]} \to \sS_p^{[k+1]}$ the coboundary operator. Let $\U^{\bullet} = \{ \U^{(p)} \}_{p \ge 0}$ be an open cover  of $X_{\bullet}$. We define a triple complex $(K^{i, j, k}, \d, \check{\delta}, \til{d})$ by
\begin{equation}
K^{i, j, k} =
\prod_{U^{(i)}_{\alpha^{(i)}_0}, \ldots, U^{(i)}_{\alpha^{(i)}_j}}
\Gamma(U^{(i)}_{\alpha^{(i)}_0 \ldots \alpha^{(i)}_j}, \sS_i^{[k]}),
\end{equation}
where we write $U^{(i)}_{\alpha^{(i)}_0 \ldots \alpha^{(i)}_j}$ for the intersection $U^{(i)}_{\alpha^{(i)}_0} \cap \cdots \cap U^{(i)}_{\alpha^{(i)}_j}$. The coboundary operator $\d : K^{i, j, k} \to K^{i+1, j, k}$ is given by $\d = \sum_{l=0}^{i+1}(-1)^l\d^*_l$, the $\check{\delta} : K^{i, j, k} \to K^{i, j+1, k}$ is the \Cech coboundary operator, and $\til{d} : K^{i, j, k} \to K^{i, j, k+1}$ is induced from the complex $(\sS_i^{[k]}, \til{d})$. From the triple complex, we obtain the total complex by putting $C^m(\U^{\bullet}, \sS_{\bullet}) = \oplus_{m = i + j + k}K^{i, j, k}$, where the total coboundary operator is defined by $D = \d + (-1)^i \check{\delta} + (-1)^{i+j}\til{d}$ on the component $K^{i, j, k}$. We denote by $H^m(\U^{\bullet}, \sS_{\bullet})$ the cohomology of this total complex. Now the cohomology $H^m(X_{\bullet}, \sS_{\bullet})$ of the complex of simplicial sheaves $\sS_{\bullet}$ is defined to be
$$
H^m(X_{\bullet}, \sS_{\bullet}) =
\varinjlim H^m(\U^{\bullet}, \sS_{\bullet}),
$$
where the direct limit is taken over the ordered set of open covers of $X_{\bullet}$. 

When each open cover $\U^{(p)}$ is a \textit{good cover} \cite{B-T,Bry1} of $X_p$, we call $\U^{\bullet}$ a \textit{good cover} of $X_{\bullet}$. If $\U^{\bullet}$ is a good cover, then there is a natural isomorphism $H^m(\U^{\bullet}, \sS_{\bullet}) \cong H^m(X_{\bullet}, \sS_{\bullet})$.


\subsection{Equivariant smooth Deligne cohomology groups}

Let $N$ be a non-negative integer.

\begin{dfn}
Let $G$ be a Lie group acting on a smooth manifold $M$. We define a complex of simplicial sheaves $\F(N)_{G^{\bullet} \times M}$ on $G^{\bullet} \times M$ by the sequence of smooth Deligne complex $\{ \F(N)_{G^p \times M} \}_{p \ge 0}$, where the homomorphisms $\til{\d}_i : \d_i^{-1} \F(N)_{G^p \times M} \to \F(N)_{G^{p+1} \times M}$ and $\til{s}_i : s_i^{-1} \F(N)_{G^{p+1} \times M} \to \F(N)_{G^p \times M}$ are the natural ones.
\end{dfn}

We often omit the subscripts of $\F(N)_{G^{\bullet} \times M}$, $\F(N)_{G^p \times M}$, etc.

The cohomology $H^m(G^{\bullet} \times M, \F(N))$ can be identified with the ordinary smooth Deligne cohomology on the quotient space under some assumptions.

\begin{prop}[\cite{Go2}] \label{prop:quotient_Deligne}
Let $G$ be a Lie group acting on a smooth manifold $M$. We assume that the action is free and locally trivial, and that the quotient space $M/G$ is a smooth manifold in such a way that the projection map $q : M \to M/G$ is smooth. For a non-negative integer $N$ the projection map induces an isomorphism of groups
$$
q^* : \ H^N(M/G, \F(N)_{M/G}) \longrightarrow 
H^N(G^{\bullet} \times M, \F(N)_{G^\bullet \times M}).
$$
\end{prop}

This proposition suggests that $H^m(G^{\bullet} \times M, \F(N))$ is rather unsuitable for the classification of equivariant geometric objects, such as equivariant principal $\T$-bundles with connection. So we introduce the other cohomology group in the below.

\bigskip

For a moment, we fix a non-negative integer $i$. We have an obvious fibration $\pi : G^i \times M \to G^i \times pt$, where $pt$ is the manifold consisting of a single point. For a positive integer $p$, we define a subsheaf $F^p\!\u{A}^q_{G^i \times M}$ of $\u{A}^q_{G^i \times M}$ by setting $F^p\!\u{A}^q_{G^i \times M} = \pi^{-1}\u{A}^p_{G^i \times pt} \otimes \u{A}^{q-p}_{G^i \times M}$, where the tensor product is taken over $\pi^{-1}\u{A}^0_{G^i \times pt}$.

For an open subset $U \subset G^i \times M$, the group $F^p\!\u{A}^q_{G^i \times M}(U)$ consists of those $q$-forms $\omega$ on $U$ satisfying $\iota_{V_1} \cdots \iota_{V_{q-p+1}} \omega = 0$ for tangent vectors $V_1, \ldots, V_{q-p+1}$ at $x \in U$ such that $\pi_*V_k = 0$. If $\{ g_j \}$ and $\{ x_k \}$ are systems of local coordinates of $G$ and $M$ respectively, then the $q$-form $\omega$ has a local expression
$$
\omega = 
\sum_{r \ge p} 
\sum_{\substack{j_1, \ldots, j_r \\ k_1, \ldots, k_{q-r}}}
f_{J, K}(g, x) 
dg_{j_1} \wedge \cdots \wedge dg_{j_r} \wedge
dx_{k_1} \wedge \cdots \wedge dx_{k_{q-r}}.
$$

We define the sheaf of germs of \textit{relative} $q$-forms with respect to the fibration $G^i \times M \to G^i \times pt$ by $\u{A}^q_{rel} = \u{A}^q_{G^i \times M} / F^1\!\u{A}^q_{G^i \times M}$. We also define a complex of sheaves $\bF(N)_{G^i \times M}$ on $G^i \times M$ by
$$
\bF(N)_{G^i \times M} : \
\u{\T} \stackrel{\frac{1}{2\pi\im}d\log}{\longrightarrow}
\u{A}^1_{rel} \stackrel{d}{\longrightarrow} 
\u{A}^2_{rel} \stackrel{d}{\longrightarrow} 
\cdots \stackrel{d}{\longrightarrow}
\u{A}^N_{rel} \longrightarrow
0 \longrightarrow \cdots.
$$

\begin{dfn}
Let $G$ be a Lie group acting on a smooth manifold $M$. We define a complex of simplicial sheaves $\bF(N)_{G^{\bullet} \times M}$ on $G^{\bullet} \times M$ by the sequence $\{ \bF(N)_{G^i \times M} \}_{i \ge 0}$.The homomorphisms $\til{\d}_i : \d_i^{-1} \bF(N)_{G^p \times M} \to \bF(N)_{G^{p+1} \times M}$ and $\til{s}_i : s_i^{-1} \bF(N)_{G^{p+1} \times M} \to \bF(N)_{G^p \times M}$ are the natural ones.
\end{dfn}

In \cite{Go1,Go2}, the cohomology groups $H^m(G^{\bullet} \times M, \bF(N))$ are called \textit{equivariant smooth Deligne cohomology groups}.

Note that, if the topology on $G$ is discrete, then we have $\F(N) = \bF(N)$, so that the cohomology $H^m(G^{\bullet} \times M, \bF(N))$ coincides with $H^m(G^{\bullet} \times M, \F(N))$. 

\begin{rem}
In \cite{L-U}, Lupercio and Uribe introduced the \textit{Deligne cohomology group for the orbifold $M/G$} by the cohomology $H^m(G^{\bullet} \times M, \F(N))$, where $G$ is a finite group. Since the topology on $G$ is discrete, we have $H^m(G^{\bullet} \times M, \F(N)) \cong H^m(G^{\bullet} \times M, \bF(N))$.
\end{rem}

\medskip

When the Lie algebra of $G$ is non-trivial, the cohomology $H^m(G^{\bullet} \times M, \F(N))$ differs from $H^m(G^{\bullet} \times M, \bF(N))$ generally. To see a precise relation between these cohomology groups, we introduce a subcomplex of $\F(N)_{G^{\bullet} \times M}$ as follows. 

Let us fix again a non-negative integer $i$ for a moment. By means of the fibration $G^i \times M \to G^i \times pt$, we can form a subcomplex $F^1\!\F(N)_{G^i \times M}$ of the ordinary smooth Deligne complex $\F(N)_{G^i \times M}$:
$$
\xymatrix@C=15pt@R=2pt{
\u{\T} \ar[r] &
\u{A}^1 \ar[r] &
\u{A}^2 \ar[r] &
\cdots \ar[r] &
\u{A}^{N-1} \ar[r] &
\u{A}^N \ar[r] &
0 \ar[r] &
\cdots, \\
\cup &
\cup &
\cup &
\cdots &
\cup &
\cup &
\cup & \\
0 \ar[r] &
F^1\!\u{A}^1 \ar[r] &
F^1\!\u{A}^2 \ar[r] &
\cdots \ar[r] &
F^1\!\u{A}^{N-1} \ar[r] &
F^1\!\u{A}^N \ar[r] &
0 \ar[r] &
\cdots. \\
}
$$
Now we define the subcomplex $F^1\!\F(N)_{G^{\bullet} \times M}$ by $\{ F^1\!\F(N)_{G^i \times M} \}_{i \ge 0}$. 

Notice that, for each $i$, the complex $\bF(N)_{G^i \times M}$ is obtained as the quotient $\F(N)_{G^i \times M}/F^1\!\F(N)_{G^i \times M}$. Hence $\bF(N)_{G^{\bullet} \times M}$ fits into the the following short exact sequence:
\begin{equation}
0 \longrightarrow
F^1\!\F(N)_{G^{\bullet} \times M} \longrightarrow
\F(N)_{G^{\bullet} \times M} \stackrel{\varphi}{\longrightarrow}
\bF(N)_{G^{\bullet} \times M} \longrightarrow 0.
\label{exact_seq:simplicial_complex}
\end{equation}
This induces a long exact sequence of cohomology groups:
$$
\xymatrix@C=7pt@R=2pt{
&
\quad \quad \quad \quad \quad \quad \quad \quad \cdots \ar[r] &
H^{m-1}(G^{\bullet} \! \times \! M, \F(N)) \ar[r]^{\varphi} &
H^{m-1}(G^{\bullet} \! \times \! M, \bF(N)) \\
\ar[r] &
H^m(G^{\bullet} \! \times \! M, F^1\!\F(N)) \ar[r] &
H^m(G^{\bullet} \! \times \! M, \F(N)) \ar[r]^{\varphi} &
H^m(G^{\bullet} \! \times \! M, \bF(N)) \\
\ar[r] &
H^{m+1}(G^{\bullet} \! \times \! M, F^1\!\F(N)) \ar[r] &
H^{m+1}(G^{\bullet} \! \times \! M, \F(N)) \ar[r] &
\cdots. \quad \quad \quad  \quad \quad \quad
}
$$
We denote the Bockstein homomorphism in the exact sequence above by 
$$
\beta : 
H^m(G^{\bullet} \times M, \bF(N)) \longrightarrow
H^{m+1}(G^{\bullet} \times M, F^1\!\F(N)). 
$$
Because $G^i \times M$ is assumed to admit a partition of unity for each $i$, the cohomology $H^m(G^{\bullet} \times M, F^1\!\F(N))$ is computed as the $m$th cohomology of the double complex $(L^{i, j}, \d, \til{d})$ given by
\begin{equation}
L^{i, j} = 
\left\{
\begin{array}{cc}
F^1\!A^j(G^i \times M), & (1 \le j \le N), \\
0, & \text{otherwise}.
\end{array}
\right. \label{double_complex:F1FN}
\end{equation}
As is clear, $H^m(G^{\bullet} \times M, F^1\!\F(N))$ is a vector space over $\R$. Since $F^1\!A^j(M) = 0$, we have $H^m(G^{\bullet} \times M, F^1\!\F(N)) = 0$ for $m = 0, 1$.

\medskip

In the reminder, we express $H^m(G^{\bullet} \times M, F^1\!\F(N))$ in a more accessible style. Let $\g$ be the Lie algebra of $G$, $\g^*$ its dual space, $\langle \ | \ \rangle : \g \otimes \g^* \to \R$ the natural contraction, and $A^q(M, \g^*)$ the vector space of $\g^*$-valued $q$-forms on $M$. By the (co)adjoint, the Lie group $G$ acts on $\g^*$ by left:
$$
\langle X | \Ad_g f \rangle = \langle \Ad_{g^{-1}}X | f \rangle,
$$
where $X \in \g$, $f \in \g^*$ and $g \in G$. This induces a left action of $G$ on $A^q(M, \g^*)$. On the other hand, the left action of $G$ on $M$ induces a right action on $A^q(M, \g^*)$ by the pull-back: $\xi \mapsto g^*\xi$ for $\xi \in A^q(M, \g^*)$. 

\begin{lem}[\cite{Go2}] \label{lem:iso_H2F1F1_H2F1F2}
There are isomorphisms
\begin{align*}
H^2(G^{\bullet} \times M, F^1\!\F(1)) 
& \cong
\left\{ \mu \in A^0(M, \g^*) |\ 
g^*\mu = \Ad_g \mu \right\}, \\
H^2(G^{\bullet} \times M, F^1\!\F(2)) 
& \cong
\{ \mu \in A^0(M, \g^*) |\ 
g^*\mu = \Ad_g \mu, \ d\mu = 0 \}.
\end{align*}
\end{lem}

\begin{proof}
We only give maps inducing the isomorphisms in this lemma, and refer the reader to \cite{Go2} for detail. Because $F^1\!A^1(M) = F^1\!A^2(M) = 0$, we have
\begin{align*}
H^2(G^{\bullet} \times M, F^1\!\F(1)) 
& =
\{ \alpha \in  F^1\!A^1(G \times M) |\ \d \alpha = 0 \}, \\
H^2(G^{\bullet} \times M, F^1\!\F(2)) 
& =
\left\{
\alpha \in F^1\!A^1(G \times M) |\ \d \alpha = 0, \ d \alpha = 0
\right\}.
\end{align*}
For $\alpha \in F^1\!A^1(G \times M)$, we define $\mu \in A^0(M, \g^*)$ by $\langle X | \mu(x) \rangle = \alpha((e, x); X \oplus 0)$, where $x \in M$ and $X \in \g = T_eG$. Then the assignment $\alpha \mapsto f$ gives the isomorphisms.
\end{proof}

\begin{cor}[\cite{Go2}] \label{cor:vanishing_H2F1F2}
If $[\g, \g] = \g$, then $H^2(G^{\bullet} \times M, F^1\!\F(2)) = 0$.
\end{cor}

\begin{proof}
Suppose that $\mu \in A^0(M, \g^*)$ satisfies $g^*\mu = \Ad_g\mu$ and $d \mu = 0$. By these relations, we see that $\langle [X, Y] | \mu \rangle = 0$ for all $X, Y \in \g$. Thus, if $[\g, \g] = \g$, then $\mu = 0$, so that $H^2(G^{\bullet} \times M, F^1\!\F(2)) = 0$.
\end{proof}

A map $\zeta : G \to A^0(M, \g^*)$ induces a map $G \times M \to \g^*$ by $(g, x) \mapsto \zeta(g)(x)$. We say that $\zeta$ is \textit{smooth} when the induced one is. Let $C^{\infty}(G, A^0(M, \g^*))$ be the vector space of smooth maps $\zeta : G \to A^0(M, \g^*)$.

\begin{lem} \label{lem:iso_H3F1F2}
There is an isomorphism $H^3(G^{\bullet} \times M, F^1\!\F(2)) \cong \mathcal{Z}/\mathcal{B}$, where $\mathcal{Z}$ and $\mathcal{B}$ are defined by
\begin{align}
\mathcal{C} & = 
A^1(M, \g^*) \oplus C^{\infty}(G, A^0(M, \g^*)), \\
\mathcal{Z} & = 
\left\{
(E, \zeta) \in \mathcal{C} \Big|
\begin{array}{l}
g^*E - \Ad_gE = d \zeta(g), \\
\Ad_g\zeta(h) - \zeta(gh) + h^*\zeta(g) = 0. \\
\end{array}
\right\}, \\ 
\mathcal{B} & = 
\left\{
(E, \zeta) \in \mathcal{C} \Bigg|
\begin{array}{l}
\mu \in A^0(M, \g^*), \\
E = d \mu, \\
\zeta(g) = g^*\mu - \Ad_g \mu.
\end{array}
\right\}. \label{subspace:coboundary}
\end{align}
\end{lem}

\begin{proof}
Recall that $H^3(G^{\bullet} \times M, F^1\!\F(2))$ is the cohomology of (\ref{double_complex:F1FN}). Let $[\alpha, \beta]$ be a class in the cohomology, where $\alpha \in F^1\!A^2(G \times M)$ and $\beta \in F^1\!A^1(G^2 \times M)$. We define $H\!\beta \in  F^1\!A^1(G \times M)$ by setting 
$$
H\!\beta((g, x); gX \oplus V) = \beta((g, e, x); 0 \oplus X \oplus 0),
$$ 
where we expressed a tangent vector at $g \in G$ as $gX \in T_gG$ using $X \in T_eG = \g$. Note that cocycles $(\alpha, \beta)$ and $(\alpha - d H\!\beta, \beta + \d H\!\beta)$ induce the same cohomology class. Now we define $E \in A^1(M, \g^*)$ and $\zeta \in C^{\infty}(G, A^0(M, \g^*))$ by
\begin{align*}
\langle X | E(x; V) \rangle
& = 
(\alpha - d H\!\beta)((e, x); X \oplus 0, 0 \oplus V), \\
\langle X | \zeta(g)(x) \rangle
& = 
(\beta + \d H\!\beta)((e, g, x); X \oplus 0 \oplus 0).
\end{align*}
By a lengthy calculation, we can see that the cocycle condition for $(\alpha, \beta)$ implies that $(E, \zeta)$ belongs to $\mathcal{Z}$. We next consider the cochain $(\alpha, \beta)$ of the form $(\alpha, \beta) = (- d\gamma, \d \gamma)$ for an element $\gamma \in F^1\!A^1(G \times M)$. In this case, we define $\mu \in A^0(M, \g^*)$ by
$$
\langle X | \mu(x) \rangle = \gamma((e, x); X \oplus 0).
$$
By means of $\mu$ above, we can show that the $(E, \zeta) \in \mathcal{Z}$ defined by $(\alpha, \beta)$ belongs to $\mathcal{B}$. Therefore we obtain a well-defined homomorphism 
$$
\Phi :\  H^3(G^{\bullet} \times M, F^1\!\F(2)) \longrightarrow 
\mathcal{Z}/\mathcal{B}
$$ 
by $\Phi([\alpha, \beta]) = [E, \zeta]$. In order to prove that $\Phi$ is an isomorphism, it suffices to construct the inverse homomorphism. For $(E, \zeta) \in \mathcal{C}$ we define $\alpha \in F^1\!A^2(G \times M)$ and $\beta \in F^1\!A^1(G^2 \times M)$ by
\begin{align*}
\alpha((g, x); gX \oplus V, gX' \oplus V')
& =
\langle X | E(x; V') \rangle - 
\langle X' | E(x; V) \rangle \\
& +
\langle [X, X'] | \zeta(e)(x) \rangle, \\
\beta((g_1, g_2, x); g_1X_1 \oplus g_2 X_2 \oplus V)
& =
\langle X_1 | \zeta(g_2)(x) \rangle.
\end{align*}
If $(E, \zeta)$ belongs to $\mathcal{Z}$, then $(\alpha, \beta)$ is a cocycle. Suppose that $(E, \zeta)$ is expressed as in (\ref{subspace:coboundary}) by a function $\mu \in A^0(M, \g^*)$. Then we define $\gamma \in F^1\!A^1(G \times M)$ by
$$
\gamma((g, x); gX \oplus V) = \langle X | \mu(x) \rangle.
$$
We can verify that the $(\alpha, \beta)$ defined by $(E, \zeta)$ is of the form $(-d \gamma, \d \gamma)$. Therefore we obtain a homomorphism $\Psi : \mathcal{Z}/\mathcal{B} \to H^3(G^{\bullet} \times M, F^1\!\F(2))$ by $\Psi([E, \zeta]) = [\alpha, \beta]$. Note that, if $(E, \zeta) \in \mathcal{Z}$, then we have 
$$
\langle [X, X'] | \zeta(e)(x) \rangle 
=
\langle X | E(x; {X'}^* ) \rangle -
\langle X | d\zeta((e, x); X' \oplus 0) \rangle.
$$
Thus, we can see that $\Psi$ is the inverse of $\Phi$.
\end{proof}

\begin{cor}
If $G$ is compact, then we can take a cocycle of the form $(E, 0) \in \mathcal{Z}$ as a representative of a class in $H^3(G^{\bullet} \times M, F^1\!\F(2)) \cong \mathcal{Z}/\mathcal{B}$.
\end{cor}

\begin{proof}
Suppose that a cocycle $(E, \zeta) \in \mathcal{Z}$ is given. Because $G$ is compact, we can take an invariant measure $dg$ on $G$. We assume that the measure is normalized. By the cocycle condition for $\zeta$, we have $\zeta(h) = \Ad_h(\Ad_{(gh)^{-1}}\zeta(gh)) - h^*(\Ad_{g^{-1}}\zeta(g))$. Thus, if we define $\mu \in A^0(M, \g^*)$ by $\mu = \int_{g \in G}\Ad_{g^{-1}}\zeta(g) dg$, then $\zeta(h) = \Ad_h\mu - h^*\mu$, so that $[E, \zeta] = [E + d\mu, 0]$.
\end{proof}


\section{Bundle gerbes}
\label{sec:bundle_gerbes}

In this section, we define bundle gerbes \cite{Mu,Mu-S}, following Meinrenken's formulation \cite{Me}. We also recall some basic properties of bundle gerbes.



\subsection{A complex with no cohomology}

Let $M$ be a smooth manifold. Suppose that we have the other manifold $Y$ and a surjective submersion $\pi : Y \to M$ admitting local sections. The $p$-fold fiber product of $\pi : Y \to M$ is defined by $Y^{[p]} = \{ (y_1, \ldots, y_p) |\ \pi(y_1) = \cdots = \pi(y_p) \}$, and the projections $\pi_i : Y^{[p]} \to Y^{[p-1]}, \ (i = 1, \ldots, p)$ by omitting the $i$th factor. We define a homomorphism $\delta : A^q(Y^{[p-1]}) \to A^q(Y^{[p]})$ by $\delta = \sum_{i = 1}^p(-1)^{i-1}\pi_i^*$. It is easy to verify that $\delta \delta = 0$. 

\begin{lem}[\cite{Mu}] \label{lem:Murray}
For $q \ge 0$, the following sequence is exact:
$$
0 \to 
A^q(M) \stackrel{\pi^*}{\to}
A^q(Y) \stackrel{\delta}{\to}
A^q(Y^{[2]}) \stackrel{\delta}{\to}
A^q(Y^{[3]}) \stackrel{\delta}{\to} \cdots.
$$
\end{lem}

When a principal $\T$-bundle $P \to Y^{[p]}$ is given, we define a principal $\T$-bundle $\delta P \to Y^{[p+1]}$ by setting $\delta P = \pi_1^*P \otimes \pi_2^*P^{\otimes -1} \otimes \pi_3^*P \otimes \cdots \otimes \pi_p^*P^{\otimes (-1)^{p-1}}$. Notice that $\delta \delta P \to Y^{[p+2]}$ is canonically isomorphic to the trivial $\T$-bundle over $Y^{[p+2]}$. In the sequel, similar notations will be used often: a section $s : Y^{[p]} \to P$ induces $\delta s : Y^{[p+1]} \to \delta P$, and a connection $\nabla \in \im A^1(P)$ on $P$ induces a connection $\delta \nabla$ on $\delta P$. 


\subsection{Bundle gerbes}

\begin{dfn}[\cite{Mu,Mu-S}] \label{dfn:bundle_gerbe}
A \textit{bundle gerbe} $\G = (Y, P, s)$ over $M$ consists of a surjective submersion $\pi : Y \to M$ admitting local sections, a principal $\T$-bundle $P \to Y^{[2]}$, and a section $s : Y^{[3]} \to \delta P$ such that $\delta s = 1$.
\end{dfn}

For a bundle gerbe $\G = (Y, P, s)$, the \textit{inverse} $\G^{\otimes -1}$ to $\G$ is defined by $(Y, P^{\otimes -1}, s^{\otimes -1})$. When $\G' = (Y', P', s')$ is the other bundle gerbe, the \textit{product} $\G \otimes \G'$ is defined by $(Y \times_\pi Y', P \otimes P', s \otimes s')$, where $Y \times_\pi Y' \to M$ is the fiber product of $Y$ and $Y'$.

In \cite{Me}, trivializations for bundle gerbes are called ``pseudo line bundles.'' (The notion of pseudo line bundles appears in \cite{Bry-M1} as trivializations for a certain gerbes.) As principal $\T$-bundles are used in Definition \ref{dfn:bundle_gerbe}, we introduce trivializations for bundle gerbes by the name of ``pseudo $\T$-bundles.''

\begin{dfn}[\cite{Mu,Mu-S}]
A \textit{pseudo $\T$-bundle} $(R, v)$ for a bundle gerbe $\G = (Y, P, s)$ consists of a principal $\T$-bundle $R \to Y$, and a section $v : Y^{[2]} \to \delta R^{\otimes -1} \otimes P$ such that $\delta v = s$. 
\end{dfn}

A bundle gerbe admitting a pseudo $\T$-bundle is said to be \textit{trivial}.

For a bundle gerbe $\G$, there exists a cohomology class $\delta (\G) \in H^2(M, \u{\T}) \cong H^3(M, \Z)$ called the \textit{Dixmier-Douady class} of $\G$. Under the inverse and the product, the class behaves as $\delta(\G^{\otimes -1}) = - \delta(\G)$ and $\delta(\G \otimes \G') = \delta(\G) + \delta(\G')$.

\begin{prop}[\cite{Mu}]
A bundle gerbe $\G$ is trivial if and only if $\delta (\G) = 0$.
\end{prop}

For bundle gerbes, there are two notions of equivalence: \textit{isomorphism} and \textit{stable isomorphism}. The former notion is stronger than the latter.

\begin{dfn}[\cite{Mu,Mu-S}]
Let $\G = (Y, P, s)$ and $\G' = (Y', P', s')$ be bundle gerbes over $M$.

(a) We define an \textit{isomorphism} from $\G$ to $\G'$ by a pair $(\varphi, \til{\varphi})$ consisting of a fiber preserving diffeomorphism $\varphi : Y \to Y'$ and an isomorphism of $\T$-bundles $\til{\varphi} : P \to P$ covering $\varphi^{[2]} : Y^{[2]} \to (Y')^{[2]}$ such that $\delta \til{\varphi} \circ s = s' \circ \varphi^{[3]}$. Here we denote by $\varphi^{[p]} : Y^{[p]} \to (Y')^{[p]}$ the natural map induced from $\varphi$, and by $\delta \til{\varphi}: \delta P \to \delta P'$ the isomorphism of $\T$-bundles induced from $\til{\varphi}$.

(b) We define a \textit{stable isomorphism} from $\G$ to $\G'$ by a pseudo $\T$-bundle $(R, v)$ for $\G^{\otimes -1} \otimes \G'$.
\end{dfn}

By the behavior of the Dixmier-Douady class under the inverse and the product, $\G$ and $\G'$ are stably isomorphic if and only if $\delta (\G) = \delta(\G')$. Using this fact, we can see that the stable isomorphism is indeed an equivalence relation on bundle gerbes. Since isomorphic bundle gerbes have the same Dixmier-Douady class, they are stably isomorphic. However, there are bundle gerbes which are stably isomorphic but not isomorphic \cite{Mu}.

Note that the stable isomorphism classes of bundle gerbes constitute a group by the product and the inverse.

\begin{prop}[\cite{Mu-S}] \label{prop:classification_bg}
The assignment $\G \mapsto \delta(\G)$ induces an isomorphism from the group of stable isomorphism classes of bundle gerbes over $M$ to the cohomology group $H^2(M, \u{\T}) \cong H^3(M, \Z)$.
\end{prop}

Let $\Phi : M' \to M$ be a smooth map. If $\pi : Y \to M$ is a surjective submersion admitting local sections, then so is the natural map $\Phi^*Y \to M'$. Clearly, we have $(\Phi^*Y)^{[p]} = \Phi^*(Y^{[p]})$. For a bundle gerbe $\G = (Y, P, s)$ over $M$, the \textit{pull-back} $\Phi^*\G$ is defined to be $(\Phi^*Y, (\til{\Phi}^{[2]})^*P, (\til{\Phi}^{[3]})^*s)$, where $\til{\Phi}^{[p]} : \Phi^*Y^{[p]} \to Y^{[p]}$ is the natural map covering $\Phi : M' \to M$. The pull-back of a pseudo $\T$-bundle is also defined in a similar way: let $(R, v)$ be a pseudo $\T$-bundle for $\G$. Then $(\til{\Phi}^*R, (\til{\Phi}^{[2]})^*v)$ is a pseudo $\T$-bundle for $\Phi^*\G$, which we denote by $\Phi^*(R, v)$. We note that the Dixmier-Douady class is natural under the pull-back operation.


\subsection{Connection, curving and 3-curvature}

\begin{dfn}[\cite{Mu,Mu-S}]
Let $\G = (Y, P, s)$ be a bundle gerbe over $M$.

(a) A \textit{bundle gerbe connection} (or \textit{connection}, for short) $\nabla$ on $\G$ is a connection $\nabla \in \im A^1(P)$ on the principal $\T$-bundle $P$ such that $s^*(\delta \nabla) = 0$.

(b) A \textit{bundle gerbe curving} (or \textit{curving}, for short) $f$ for $\nabla$ is a 2-form $f \in \im A^2(Y)$ such that $\delta f = F(\nabla)$, where $F(\nabla)$ is the curvature of $\nabla$.

(c) The \textit{3-curvature} $\Omega$ of a bundle gerbe $\G$ with a connection $\nabla$ and a curving $f$ is the unique 3-form $\Omega \in \im A^3(M)$ such that $\pi^* \Omega = df$.
\end{dfn}

As is known \cite{Mu}, there always exists a connection on a bundle gerbe given, and a curving for a connection given. In general, there are various choices of a connection and a curving. For example, if $f$ is a curving for $\nabla$, then  so is $f + \pi^*\sigma$, where $\sigma \in \im A^2(M)$.

The inverse and the product for bundle gerbes with connection and curving are defined in a similar way: $(\G, \nabla, f)^{\otimes -1} = (\G^{\otimes -1}, \nabla^{\otimes -1}, -f)$, $(\G, \nabla, f) \otimes (\G', \nabla', f') = (\G \otimes \G', \nabla \otimes \nabla', f + f')$.

\begin{dfn}[\cite{Mu,Mu-S}]
Let $(\G, \nabla)$ be a bundle gerbe with connection over $M$, and $(R, v)$ a pseudo $\T$-bundle for $\G = (Y, P, s)$. We define a \textit{connection} $\eta$ on $(R, v)$ to be a connection $\eta \in \im A^1(R)$ on the principal $\T$-bundle $R$ such that $v^*(\delta \eta^{\otimes -1} \otimes \nabla) = 0$. 
\end{dfn}

A bundle gerbe with connection and curving $(\G, \nabla, f)$ is said to be \textit{trivial} if it admits a pseudo $\T$-bundle $(R, v)$ with a connection $\eta$ such that $F(\eta) = f$.

For a bundle gerbe with connection and curving $(\G, \nabla, f)$ over $M$, there exists a cohomology class $\delta(\G, \nabla, f) \in H^2(M, \F(2))$. Under the inverse and the product, this class behaves as $\delta((\G, \nabla, f)^{\otimes -1}) = - \delta(\G, \nabla, f)$ and $\delta((\G, \nabla, f) \otimes (\G', \nabla', f')) = \delta(\G, \nabla, f) + \delta(\G', \nabla', f')$.

\begin{prop}[\cite{Mu}] \label{ref:trivial_bgcc}
A bundle gerbe with connection and curving $(\G, \nabla, f)$ is trivial if and only if $\delta(\G, \nabla, f) = 0$.
\end{prop}

\begin{dfn}[\cite{Mu,Mu-S}]
Let $(\G, \nabla, f)$ and $(\G', \nabla', f')$ be bundle gerbes with connection and curving over $M$. 

(a) We define an \textit{isomorphism} from $(\G, \nabla, f)$ to $(\G', \nabla', f')$ by an isomorphism $(\varphi, \til{\varphi}) : \G \to \G'$ such that $\til{\varphi}^*\nabla' = \nabla$ and $\varphi^*f' = f$.

(b) We define a \textit{stable isomorphism} from $(\G, \nabla, f)$ to $(\G', \nabla', f')$ by a pseudo $\T$-bundle $(R, v)$ with a connection $\eta$ for $(\G, \nabla, f)^{\otimes -1} \otimes (\G', \nabla', f')$ such that $F(\eta) = f' - f$.
\end{dfn}

\begin{prop}[\cite{Mu-S}] \label{prop:classification_bgcc}
The assignment $\G \mapsto \delta(\G, \nabla, f)$ induces an isomorphism from the group of stable isomorphism classes of bundle gerbes with connection and curving over $M$ to $H^2(M, \F(2))$.
\end{prop}

In Proposition \ref{prop:Deligne_coh:manifold} (b), we have two surjections $H^2(M, \F(2)) \to A^3(M)_0$ and $H^2(M, \F(2)) \to H^3(M, \Z)$. The image of $\delta(\G, \nabla, f)$ under the first surjection is $\frac{-1}{2\pi\im} \Omega$, where $\Omega$ is the 3-curvature of $(\G, \nabla, f)$. The image under the second surjection is the Dixmier-Douady class $\delta(\G)$.

As is noticed, we have a choice of a curving for a bundle gerbe connection. Let $i : A^2(M)/A^2(M)_0 \to H^2(M, \F(2))$ be the injection in the second exact sequence in Proposition \ref{prop:Deligne_coh:manifold} (b). For a 2-form $\sigma \in \im A^2(M)$, we have $\delta(\G, \nabla, f + \pi^*\sigma) = \delta(\G, \nabla, f) + i(\frac{-1}{2\pi\im}\sigma)$. Hence Proposition \ref{ref:trivial_bgcc} leads to:

\begin{cor} \label{cor:different_curving}
The bundle gerbes with connection and curving $(\G, \nabla, f)$ and $(\G, \nabla, f + \pi^*\sigma)$ are stably isomorphic if and only if there exists a connection $\eta$ on a principal $\T$-bundle $R \to M$ such that $F(\eta) = \sigma$.
\end{cor}

The notion of the pull-back of bundle gerbes with connection and curving is defined in the same manner as that of bundle gerbes. The class $\delta(\G, \nabla, f)$ is also natural under the pull-back operation.


\section{Equivariant bundle gerbes}
\label{sec:EBG}

In this section, we introduce the notion of \textit{equivariant bundle gerbes} formulated by Meinrenken \cite{Me}. As is mentioned in Section \ref{sec:introduction}, the notion includes \textit{strongly equivariant bundle gerbes}, i.e.\@ the equivariant bundle gerbes studied by Mathai and Stevenson \cite{Ma-S}. We also give the definition of strongly equivariant bundle gerbes, and explain the relation with equivariant bundle gerbes.


\subsection{Differential forms in the direction of $M$}

Let $M$ be a smooth manifold on which a Lie group $G$ acts by left. Usually, we do not assume $G$ to be compact. In order to define equivariant bundle gerbe, we use the simplicial manifold $G^{\bullet} \times M = \{ G^p \times M \}_{p \ge 0}$.

Let $Y_{\bullet} = \{ Y_p \}_{p \ge 0}$ be the other simplicial manifold. We denote by $\pi : Y_{\bullet} \to G^{\bullet} \times M$ a simplicial surjective submersion admitting local sections, that is, a sequence $\{ Y_p \to G^p \times M \}_{p \ge 0}$ of surjective submersions which admit local sections for each $p$ and are compatible with the face and degeneracy maps of $Y_{\bullet}$ and $G^{\bullet} \times M$. Note that local sections are not required to be compatible with the face and degeneracy maps.

If $\pi : Y_{\bullet} \to G^{\bullet} \times M$ is such, then the sequence $Y_{\bullet}^{[q]} = \{ Y_p^{[q]} \}_{p \ge 0}$ becomes a simplicial manifold in a natural way. For a principal $\T$-bundle $P \to Y_p^{[q]}$ we define a principal $\T$-bundle $\d P \to Y_{p+1}^{[q]}$ by $\d P = \d_0^*P \otimes \d_1^*P^{\otimes -1} \otimes \d_2^*P \otimes \cdots \otimes \d_{p+1}^*P^{\otimes (-1)^{p+1}}$, where $\d_i : Y_{p+1}^{[q]} \to Y_p^{[q]}$ is the face map. It is easy to see that $\d \d P \to Y_{p+2}^{[q]}$ is canonically isomorphic to the trivial bundle. Similar notations will be used for connections, sections, etc.

\begin{rem}
The family $\{ Y_p^{[q+1]} \}_{p, q \ge 0}$ forms a bi-simplicial manifold. In \cite{Me}, an equivariant bundle gerbe is defined in terms of the bi-simplicial manifold. To see the similarity between the definition of equivariant bundle gerbes and that of equivariant gerbes \cite{Bry2}, we use this fact rather implicitly.
\end{rem}

As we see in Section \ref{sec:EDC}, we have the subgroup $F^1\!A^k(G^p \times M)$ of $A^k(G^p \times M)$ for each $p$. We define the group $A^k(G^p \times M)_{rel}$ of relative $k$-forms with respect to the fibration $G^p \times M \to G^p \times pt$, or $k$-forms in the direction of $M$, by $A^k(G^p \times M)_{rel} = A^k(G^p \times M) / F^1\!A^k(G^p \times M)$. We write $[\omega]_{rel} \in A^k(G^p \times M)_{rel}$ for the element represented by $\omega \in A^k(G^p \times M)$.

It is worth while to note that $\omega \in A^k(M)$ is $G$-invariant if and only if $[\d \omega]_{rel} = 0$ in $A^k(G \times M)_{rel}$. On the other hand, $\omega \in A^k(M)$ is $G$-invariant and vanishes in the direction of $G$, if and only if $\d \omega = 0$ in $A^k(G \times M)$.

\smallskip

Let $\pi : Y_p \to G^p \times M$ be a surjective submersion admitting local sections. Since $\pi^* : A^k(G^p \times M) \to A^k(Y_p^{[q]})$ is injective, we can regard $F^1\!A^k(G^p \times M)$ as a subgroup of $A^k(Y_p^{[q]})$. We define the group $A^k(Y_p^{[q]})_{rel}$ by setting $A^k(Y_p^{[q]})_{rel} = A^k(Y_p^{[q]}) / F^1\!A^k(G^p \times M)$. 

Let $P \to Y_p^{[q]}$ be a principal $\T$-bundle, and $\A(P)$ the space of connections on $P$. Because $\A(P)$ is an affine space under $A^1(Y_p^{[q]})$, we define the space $\A(P)_{rel}$ of \textit{relative connections} on $P$ by $\A(P)_{rel} = \A(P) / F^1\!A^1(G^p \times M)$. It is clear that $\A(P)_{rel}$ is an affine space under $A^1(Y_p^{[q]})_{rel}$. For $\nabla \in \A(P)$, we write $[\nabla]_{rel} \in \A(P)_{rel}$ and call it a \textit{relative connection}. The curvature of $[\nabla]_{rel}$ is defined to be $F([\nabla]_{rel}) = [F(\nabla)]_{rel} \in \im A^2(Y_p^{[q]})_{rel}$.


\subsection{Equivariant bundle gerbes}

We introduce here equivariant bundle gerbes formulated by Meinrenken \cite{Me}.  However, bi-simplicial manifolds are not used, and the definition of connections is slightly different.

\begin{dfn}[\cite{Me}] \label{dfn:EBG}
Let $G$ be a Lie group acting on a smooth manifold $M$.

(a) A \textit{$G$-equivariant bundle gerbe} $\G_G$ over a manifold $M$ consists of the following data: a simplicial manifold $Y_{\bullet} = \{ Y_p \}_{p \ge 0}$ with a simplicial surjective submersion $\pi : Y_{\bullet} \to G^{\bullet} \times M$ which admits local sections; a bundle gerbe $(Y_0, P, s)$; a pseudo $\T$-bundle $(Q, t)$ for $(Y_1, \d P, \d s)$; a section $u : Y_2 \to \d Q$ such that $\delta u = \d t^{\otimes -1}$ and $\d u = 1$. We write $\G_G = (Y_\bullet, (Y_0, P, s), (Q, t), u)$.

(b) A \textit{($G$-invariant) connection} $\nabla_G$ on a $G$-equivariant bundle gerbe $\G_G$ consists of the following data: a connection $\nabla$ on $(Y_0, P, s)$; a relative connection $D_{rel}$ on the principal $\T$-bundle $Q \to Y_1$ such that $t^*(\delta D_{rel}^{\otimes -1} \otimes [\d \nabla]_{rel}) = 0$ and $u^*(\d D_{rel}) = 0$. We write $\nabla_G = (\nabla, D_{rel})$.

(c) A \textit{($G$-invariant) curving} $f$ for a connection $\nabla_G = (\nabla, D_{rel})$ is defined to be a curving $f$ for $\nabla$ such that $F(D_{rel}) = [\d f]_{rel}$.
\end{dfn}

We often mean an ``equivariant bundle gerbe with invariant connection and invariant curving'' by an ``equivariant bundle gerbe with connection and curving.'' (We drop the word ``invariant'' for brevity.)

If $G$ is compact, then each equivariant bundle gerbe admits invariant connections and invariant curvings \cite{Me}. The choice of invariant connections and invariant curvings are not unique. 

When $(\G_G, \nabla_G, f)$ is a $G$-equivariant bundle gerbe with connection and curving over $M$, the 3-curvature $\Omega$ of the curving $f$ is a closed $G$-invariant 3-form.

As in the case of ordinary bundle gerbes, the inverse and the product are defined in a similar fashion.

\begin{dfn}[\cite{Me}]
Let $\G_G = (Y_\bullet, (Y_0, P, s), (Q, t), u)$ be a $G$-equivariant bundle gerbe over $M$ equipped with a $G$-invariant connection $\nabla_G = (\nabla, D_{rel})$.

(a) A \textit{($G$-equivariant) pseudo $\T$-bundle} for $\G_G$ is defined to be a pseudo $\T$-bundle $(R, v)$ for $(Y_0, P, s)$ together with a section $r : Y_1 \to \d R \otimes Q^{\otimes -1}$ such that $\delta r = \d v^{\otimes -1} \otimes t$ and $\d r = u^{\otimes -1}$. 

(b) A \textit{($G$-invariant) connection} on $((R, v), r)$ is defined to be a connection $\eta$ on $(R, v)$ such that $r^*([\d \eta]_{rel} \otimes D_{rel}^{\otimes -1}) = 0$. 
\end{dfn}

A $G$-equivariant bundle gerbe with connection and curving $(\G_G, \nabla_G, f)$ is said to be \textit{trivial} if it admits a $G$-equivariant pseudo $\T$-bundle together with a $G$-invariant connection $\eta$ such that $F(\eta) = f$.

\begin{dfn}
We define a \textit{stable isomorphism} between $G$-equivariant bundle gerbes with connection and curving $(\G_G, \nabla_G, f)$ and $(\G'_G, \nabla'_G, f')$ over $M$ to be a pseudo $\T$-bundle with connection $((R, v), r)$ for $(\G_G, \nabla_G, f) \otimes (\G'_G, \nabla'_G, f')^{\otimes -1}$ equipped with a $G$-invariant connection $\eta$ such that $F(\eta) = f - f'$.
\end{dfn}

Let $G$ and $G'$ be Lie groups acting on smooth manifolds $M$ and $M'$, respectively. If we have a homomorphism $\phi : G' \to G$ and a map $\Phi : M' \to M$ compatible with the actions, then we naturally obtain a simplicial map ${G'}^{\bullet} \times M' \to G^{\bullet} \times M$. We also denote this simplicial map by $\Phi$. For a $G$-equivariant bundle gerbe $\G_G = (Y_\bullet, (Y_0, P, s), (Q, t), u)$ over $M$, the \textit{pull-back} $\Phi^*\G_G$ is defined to be $(\Phi^*Y_\bullet, \Phi^*(Y_0, P, s), \Phi^*(Q, t), \Phi^*u)$. This is a $G'$-equivariant bundle gerbe over $M'$. In addition, a $G$-invariant connection $\nabla_G$ on $\G_G$ induces a $G'$-invariant connection $\Phi^*\nabla_G$ on $\Phi^*\G_G$, and a $G$-invariant curving $f$ for $\nabla_G$ a $G'$-invariant curving $\Phi^*f$ for $\Phi^*\nabla_G$.


\subsection{Strongly equivariant bundle gerbes}

Let $M$ be a smooth manifold with a left action of a Lie group $G$, and $\pi : Y \to M$ a surjective submersion admitting a local section. Suppose that the action of $G$ on $M$ lifts to that on $Y$. In this case, we have a natural action of $G$ on the fiber product $Y^{[p]}$ compatible with the projections $\pi : Y^{[p]} \to M$ and $\pi_i : Y^{[p]} \to Y^{[p-1]}$.

\begin{dfn}[\cite{Ma-S}] \label{dfn:sEBG}
Let $G$ be a Lie group acting on a smooth manifold $M$. 

(a) A bundle gerbe $(Y, P, s)$ is said to be \textit{strongly $G$-equivariant} if the following holds: the action of $G$ lifts to that on $Y$; the action of $G$ on $Y^{[2]}$ lifts to that on the principal $\T$-bundle $P \to Y^{[2]}$ by bundle isomorphisms; and the section $s : Y^{[3]} \to \delta P$ is $G$-invariant.

(b) A \textit{($G$-invariant) connection} on a strongly $G$-equivariant bundle gerbe $(Y, P, s)$ is defined to be a connection $\nabla$ on $(Y, P, s)$ which is $G$-invariant as a 1-form on $P$.

(c) A \textit{($G$-invariant) curving} for a $G$-invariant connection $\nabla$ on a strongly $G$-equivariant bundle gerbe $(Y, P, s)$ is defined to be a curving $f$ for $\nabla$ which is $G$-invariant as a 2-form on $Y$.
\end{dfn}

\begin{lem}[\cite{Me}] \label{lem:sEBG_to_EBG}
(a) A strongly $G$-equivariant bundle gerbe $\G$ over $M$ induces a $G$-equivariant bundle gerbe $\G_G$ over $M$.

(b) A $G$-invariant connection $\nabla$ on $\G$ induces a $G$-invariant connection $\nabla_G$ on the induced $G$-equivariant bundle gerbe $\G_G$. 

(c) There exists a one to one correspondence between $G$-invariant curvings for $\nabla$ and those for $\nabla_G$.
\end{lem}

\begin{proof}
For (a), let $\G = (Y, P, s)$ be a strongly $G$-equivariant bundle gerbe over $M$. Then we have a simplicial manifold $G^{\bullet} \times Y$ and a simplicial surjective submersion $\pi : G^{\bullet} \times Y \to G^{\bullet} \times M$ which admits local sections. We define a section $t : G \times Y^{[2]} \to \d P$ by $t(g, y_1, y_2) = p_{12} \otimes (gp_{12})^{\otimes -1}$, where $p_{12}$ is a point on the fiber of $P$ at $(y_1, y_2) \in Y^{[2]}$. As is known \cite{Bry2}, the section obeys $\d t = 1$, since $P$ is $G$-equivariant. Because $s$ is $G$-invariant, we have $\delta t = \d s$. Thus, if we put $Y_{\bullet} = G^{\bullet} \times Y$, $Q = Y_1 \times \T$ and $u = 1$, then $\G_G = (Y_\bullet, (Y, P, s), (Q, t), u)$ is a $G$-equivariant bundle gerbe in the sense of Definition \ref{dfn:EBG}.

For (b), observe that a $G$-invariant connection $\nabla$ on the $G$-equivariant $\T$-bundle $P \to Y^{[2]}$ is the same thing as a connection $\nabla$ on $P$ such that $[t^*(\d \nabla)]_{rel} = 0$. If we take a relative connection $D_{rel}$ on the trivial bundle $Q$ to be the trivial connection, then $\nabla_G = (\nabla, D_{rel})$ is a connection on $\G_G$ in the sense of Definition \ref{dfn:EBG}.

For (c), recall that a 2-form $f$ on $Y$ is $G$-invariant if and only if $[\d f]_{rel} = 0$. Since we have $F(D_{rel}) = 0$, a $G$-invariant curving for $\nabla$ is the same thing as a $G$-invariant curving for $\nabla_G$.
\end{proof}

We remark that invariant connections on a strongly equivariant bundle gerbe $\G$ are not in one to one correspondence with those on $\G_G$ generally. 

In the same way as the ordinary bundle gerbes, we can define the product and the inverse for strongly equivariant bundle gerbes. Under these operations, strongly equivariant bundle gerbes are closed. We can readily see that these operations are compatible with the constructions in Lemma \ref{lem:sEBG_to_EBG}. The same thing holds, when we take connections and curvings into account.

\begin{dfn}
Let $\G$ be a strongly $G$-equivariant bundle gerbe over $M$, and $\nabla$ a $G$-invariant connection on $\G$.

(a)  We define a \textit{($G$-equivariant) pseudo $\T$-bundle} for $\G$ to be a pseudo $\T$-bundle $(R, v)$ for $\G$ (regarded as an ordinary bundle gerbe) such that $R$ is a $G$-equivariant $\T$-bundle and $v$ is a $G$-invariant section.

(b) We define a \textit{($G$-invariant) connection} on a $G$-equivariant pseudo $\T$-bundle $(R, v)$ for $\G$ to be a connection $\eta$ on $(R, v)$ which is $G$-invariant as a 1-form on $R$.
\end{dfn}

\begin{lem} \label{lem:pTsEBG_to_pTEBG}
Let $(\G, \nabla)$ be a strongly $G$-equivariant bundle gerbe with connection over $M$, and $(\G_G, \nabla_G)$ the $G$-equivariant bundle gerbe with connection induced from $(\G, \nabla)$. There is a one to one correspondence between $G$-equivariant pseudo $\T$-bundles with connection for $(\G, \nabla)$ and those for $(\G_G, \nabla_G)$.
\end{lem}

\begin{proof}
Let $(R, v)$ be a $G$-equivariant pseudo $\T$-bundle for $\G$ equipped with a $G$-invariant connection $\eta$. We define a section $r : G \times Y \to \d R$ by $r(g, y) = p \otimes (gp)^{\otimes -1}$, where $p \in R$ is a point on the fiber at $y$. Then we have $\d r = 1$ and $[r^*\d \eta]_{rel} = 0$. By means of the $G$-invariance of the section $v$, we can see $\delta r = \d v^{\otimes -1}$. Hence $((R, v), r)$ gives rise to a $G$-equivariant pseudo $\T$-bundle for $\G_G$, and $\eta$ is a $G$-invariant connection on $((R, v), r)$. Conversely, let $((R, v), r)$ is a $G$-equivariant pseudo $\T$-bundle for $\G_G$ equipped with a $G$-invariant connection $\eta$. The conditions $\d r = 1$ and $[r^*\d \eta]_{rel} = 0$ mean that $(R, \eta)$ is a $G$-equivariant $\T$-bundle with $G$-invariant connection. The condition $\delta r = \d v^{\otimes -1}$ implies the $G$-invariance of the section $v$.
\end{proof}

A strongly $G$-equivariant bundle gerbe with connection and curving $(\G, \nabla, f)$ over $M$ is said to be \textit{trivial}, when it admits a $G$-equivariant pseudo $\T$-bundle $(R, v)$ equipped with a $G$-invariant connection $\eta$ such that $F(\eta) = f$. As a direct consequence of Lemma \ref{lem:pTsEBG_to_pTEBG}, $(\G, \nabla, f)$ is trivial if and only if the induced $G$-equivariant bundle gerbe with connection and curving $(\G_G, \nabla_G, f)$ is trivial.

\begin{dfn}
Let $(\G, \nabla, f)$ and $(\G', \nabla', f')$ be strongly $G$-equivariant bundle gerbes with connection and curving. A \textit{stable isomorphism} from $(\G, \nabla, f)$ to $(\G', \nabla', f')$ is defined to be a $G$-equivariant pseudo $\T$-bundle with connection $((R, v), \eta)$ for $(\G, \nabla, f)^{\otimes -1} \otimes (\G', \nabla', f')$ such that $F(\eta) = f' - f$. 
\end{dfn}

The following is also a consequence of Lemma \ref{lem:pTsEBG_to_pTEBG}.

\begin{lem}
$(\G, \nabla, f)$ and $(\G', \nabla', f')$ are stably isomorphic if and only if the induced ones are.
\end{lem}

Let $G$ and $G'$ be Lie groups acting on smooth manifolds $M$ and $M'$, respectively. We also let $\phi : G' \to G$ be a homomorphism, and $\Phi : M' \to M$ a map compatible with the group actions. When $\G = (Y, P, s)$ is a strongly $G$-equivariant bundle gerbe over $M$, the pull-back $\Phi^*\G =(\Phi^*Y, (\til{\Phi}^{[2]})^*P, (\til{\Phi}^{[3]})^*s)$ becomes a strongly $G'$-equivariant bundle gerbe over $M'$ through the homomorphism $\phi$.

We denote by $(\Phi^*\G)_{G'}$ the $G'$-equivariant bundle gerbe induced from $\Phi^*\G$. Because we can identify $\Phi^*(G^\bullet \times Y)$ with ${G'}^\bullet \times \Phi^*Y$ naturally , we can also identify $(\Phi^*\G)_{G'}$ with $\Phi^*(\G_G)$. In other words, the construction in Lemma \ref{lem:sEBG_to_EBG} is compatible with the pull-back operations. The constructions of connections and curvings are compatible with the pull-back operations as well.


\section{Characteristic class}
\label{sec:characteristic_class}

The main purpose of this section is to construct a ``characteristic class'' for an equivariant bundle gerbe with connection and curving. We define it as an element belonging to the equivariant smooth Deligne cohomology group $H^2(G^\bullet \times M, \bF(2))$. The cohomology class is a key to the study of the obstruction to the reduction, and will be used in the next section. 

In view of the aim of this paper, it suffices to consider \textit{strongly} equivariant bundle gerbes only. However, we will deal with general equivariant bundle gerbes, because there is no essential difference in constructing the class. 


\subsection{Local data}

To define a cohomology class in $H^2(G^{\bullet} \times M, \bF(2))$ for each $G$-equivariant bundle gerbe with connection and curving over $M$, we use certain \textit{local data}. Although there does not generally exist a pseudo $\T$-bundle for an equivariant bundle gerbe, it exists locally. Such local pseudo $\T$-bundles constitute local data. The characteristic class will be defined as the obstruction to gluing local data together to construct a pseudo $\T$-bundle globally. (We can find this idea in \cite{Bry2}.)

\medskip

Let $\U^{\bullet} = \{ \U^{(p)} = \{ U^{(p)}_\alpha \}_{\alpha \in \AA^{(p)}} \}_{p \ge 0}$ be an open cover of $G^{\bullet} \times M$. For an equivariant bundle gerbe $\G_G = (Y_\bullet, (Y_0, P, s), (Q, t), u)$ equipped with a connection $\nabla_G = (\nabla, D_{rel})$ and a curving $f$, we consider the following data (i)--(iii).

\begin{list}{}{\parsep=-2pt\topsep=4pt}
\item[(i).]
For $\alpha \in \AA^{(0)}$, a pseudo $\T$-bundle $(R_\alpha, v_\alpha)$ for $(Y_0, P, s)|_{U^{(0)}_\alpha}$ and a connection $\eta_\alpha$ on $(R_\alpha, v_\alpha)$.

\item[(ii).]
For $\alpha, \beta \in \AA^{(0)}$ such that $U^{(0)}_{\alpha \beta} = \emptyset$, a section $w_{\alpha \beta} : Y_0|_{U^{(0)}_{\alpha \beta}} \to R_\alpha \otimes R_\beta^{\otimes -1}$ satisfying $\delta w_{\alpha \beta} = v_\alpha^{\otimes -1} \otimes v_\beta$ on $Y_0^{[2]}|_{U^{(0)}_{\alpha \beta}}$.

\item[(iii).] 
For $\alpha \in \AA^{(1)}$, a section $r_\alpha : Y_1|_{U^{(1)}_{\alpha}} \to \d R_\alpha \otimes Q^{\otimes -1}$ such that $\delta r_\alpha = \d v_\alpha^{\otimes -1} \otimes t$ on $Y_1^{[2]}|_{U^{(1)}_{\alpha}}$.
\end{list}

We make a remark on the notation $\d R_\alpha$ above. By the hypothesis, we have the following commutative diagram for $i = 0, 1$:
$$
\xymatrix{
Y_1|_{U^{(1)}_\alpha} \ar[r]^-{\d_i} \ar[d] & 
Y_0|_{U^{(0)}_{\d_i(\alpha)}} \ar[d] \\
U^{(1)}_\alpha \ar[r]^-{\d_i} & 
U^{(0)}_{\d_i(\alpha)}.
}
$$
We mean by $\d R_\alpha$ the principal $\T$-bundle $\d_0^*R_{\d_0(\alpha)} \otimes \d_1^*R_{\d_1(\alpha)}^{\otimes -1}$ over $Y_1|_{U^{(1)}_\alpha}$. The meaning of $\d v_\alpha$ is the same.

\begin{lem} \label{lem:construct_local_data}
Let $\U^{\bullet} = \{ \U^{(p)} \}_{p \ge 0}$ be an open cover of $G^{\bullet} \times M$. For $p = 0, 1$ we suppose that $\U^{(p)}$ is a sufficiently fine good cover of $G^p \times M$ so that we can take local sections $\psi^{(p)}_\alpha : U^{(p)}_\alpha \to Y_p|_{U^{(p)}_\alpha}$. Then there exist local data $\{ (R_\alpha, v_\alpha), \eta_\alpha, w_{\alpha \beta}, r_\alpha \}$ for $(\G_G, \nabla_G, f)$.
\end{lem}

We remark that $\psi^{(0)}_\alpha$ and $\psi^{(1)}_\alpha$ are not necessarily compatible with the face maps. (If $G$ is compact, then we can take such compatible sections \cite{Me}.) We also remark that, for an open cover $\U^{\bullet}$ given, we can construct an refinement of $\U^{\bullet}$ which satisfies the assumption in the lemma above.

\begin{proof}
Using a local section $\psi^{(0)}_\alpha : U^{(0)}_\alpha \to Y_0|_{U^{(0)}_\alpha}$, we define, for $p \ge 1$, a map
$$
(i_{12 \cdots p}, \ \psi^{(0)}_\alpha \circ \pi) : \
Y_0^{[p]}|_{U^{(0)}_\alpha} \longrightarrow 
Y_0^{[p+1]}|_{U^{(0)}_\alpha}
$$
by $(y_1, \ldots, y_p) \mapsto (y_1, \ldots, y_p, \psi^{(0)}_\alpha \circ \pi(y_1))$. For $p > 1$, these maps satisfy
$$
\pi_i \circ (i_{12 \cdots p}, \ \psi^{(0)}_\alpha \circ \pi) =
\left\{
\begin{array}{cl}
(i_{12 \cdots p-1}, \ \psi^{(0)}_\alpha \circ \pi) \circ \pi_i, & 
(i = 1, \ldots, p), \\
\id, & (i = p+1).
\end{array}
\right.
$$
Now we define a principal $\T$-bundle by $R_\alpha = (i_1, \ \psi^{(0)}_\alpha \circ \pi)^*P^{\otimes -1}$. Then a section $v_\alpha : Y_0^{[2]}|_{U^{(0)}_\alpha} \to \delta R_\alpha^{\otimes -1} \otimes P$ is obtained by setting $v_\alpha = (i_{12}, \psi^{(0)}_\alpha \circ \pi)^*s$. A computation shows that $(i_{123}, \psi^{(0)}_\alpha \circ \pi)^*\delta s = \delta v_\alpha \otimes s^{\otimes -1}$. Thus, $\delta s = 1$ leads to $\delta v_\alpha = s$, and we have constructed a pseudo $\T$-bundle $(R_\alpha, v_\alpha)$ for each $\alpha$. We also define a connection $\eta_\alpha$ on $R_\alpha$ by $\eta_\alpha = (i_1, \psi^{(0)}_\alpha \circ \pi)^*\nabla$. Because $v_\alpha^*(\delta \eta_\alpha^{\otimes -1} \otimes \nabla) = (i_{12}, \psi^{(0)}_\alpha \circ \pi)^*(s^*(\delta \nabla))$, we have a connection $\eta_\alpha$ on $(R_\alpha, v_\alpha)$.

\medskip

Next, we suppose that $U^{(0)}_{\alpha \beta} = U^{(0)}_\alpha \cap U^{(0)}_\beta$ is not empty. We take sections $\psi^{(0)}_\alpha : U^{(0)}_\alpha \to Y_0|_{U^{(0)}_\alpha}$ and $\psi^{(0)}_\beta : U^{(0)}_\beta \to Y_0|_{U^{(0)}_\beta}$. For $p \ge 1$ we define a map
$$
(i_{12 \cdots p}, \
\psi^{(0)}_\alpha \circ \pi, \ 
\psi^{(0)}_\beta \circ \pi) : \
Y_0^{[p]}|_{U^{(0)}_{\alpha \beta}} \longrightarrow 
Y_0^{[p+2]}|_{U^{(0)}_{\alpha \beta}}
$$
by $(y_1, \ldots, y_p) \mapsto (y_1, \ldots, y_p, \psi^{(0)}_\alpha \circ \pi(y_1), \psi^{(0)}_\beta \circ \pi(y_1))$. This map satisfies
$$
\pi_i \circ 
(i_{12 \cdots p}, \psi^{(0)}_\alpha \circ \pi, \psi^{(0)}_\beta \circ \pi) 
=
\left\{
\begin{array}{cl}
(i_{12 \cdots p-1}, \psi^{(0)}_\alpha \circ \pi, \psi^{(0)}_\beta \circ \pi) 
\circ \pi_i, &
(i = 1,\! .., p), \\
(i_{12 \cdots p}, \ \psi^{(0)}_\beta \circ \pi), & (i = p+1), \\
(i_{12 \cdots p}, \ \psi^{(0)}_\alpha \circ \pi), & (i = p+2), 
\end{array}
\right.
$$
In the case of $p = 1$, the relation above gives a canonical isomorphism
$$
(i_1, \ \psi^{(0)}_\alpha \circ \pi, \ \psi^{(0)}_\beta \circ \pi)^* \delta P
=
\pi^*(\psi^{(0)}_\alpha, \psi^{(0)}_\beta)^*P \otimes 
R_\beta \otimes 
R_\alpha^{\otimes -1}.
$$
By the assumption, $U^{(0)}_{\alpha \beta}$ is a contractible open set. Hence we take a section $\sigma_{\alpha \beta} : U^{(0)}_{\alpha \beta} \to (\psi^{(0)}_\alpha, \psi^{(0)}_\beta)^*P$, and define a section $w_{\alpha \beta} : Y_0|_{U^{(0)}_{\alpha \beta}} \to R_\alpha \otimes R_\beta^{\otimes -1}$ by setting $w_{\alpha \beta} = (\pi^*\sigma_{\alpha \beta}) \otimes (i_1, \psi^{(0)}_\alpha \circ \pi, \psi^{(0)}_\beta \circ \pi)^*s^{\otimes -1}$. By a calculation, we obtain $(i_{12}, \psi^{(0)}_\alpha \circ \pi, \psi^{(0)}_\beta \circ \pi)^*\delta s^{\otimes -1} = \delta w_{\alpha \beta} \otimes v_\alpha \otimes v_\beta^{\otimes -1}$. Because $\delta s = 1$, we have constructed a section $w_{\alpha \beta}$ such that $\delta w_{\alpha \beta} = v_\alpha^{\otimes -1} \otimes v_\beta$.

\medskip

Finally, we construct a section $r_\alpha : Y_1|_{U^{(1)}_{\alpha}} \to \d R_\alpha \otimes Q^{\otimes -1}$ such that $\delta r_\alpha = \d v_\alpha^{\otimes -1} \otimes t$ on $Y_1^{[2]}|_{U^{(1)}_{\alpha}}$ for $\alpha \in \AA^{(1)}$. Let $\psi^{(1)}_\alpha : U^{(1)}_\alpha \to Y_1|_{U^{(1)}_\alpha}$ be a local section. We define a map
$$
(i_1, \ \psi^{(1)}_\alpha \circ \pi) : \
Y_1^{[p]}|_{U^{(1)}_\alpha} \longrightarrow 
Y_1^{[p+1]}|_{U^{(1)}_\alpha}
$$
by $y_1 \mapsto (y_1, \psi^{(1)}_\alpha \circ \pi(y_1))$. For $i = 0, 1$ we also define a map
$$
(i_1 \circ \d_i, \
\psi^{(0)}_{\d_i(\alpha)} \circ \d_i \circ \pi, \
\d_i \circ \psi^{(1)}_\alpha \circ \pi) : \
Y_1|_{U^{(1)}_\alpha} \longrightarrow 
Y_0^{[3]}|_{U^{(0)}_{\d_i(\alpha)}}
$$
by $y \mapsto (\d_i(y), \ \psi^{(0)}_{\d_i(\alpha)} \circ \d_i \circ \pi(y), \ \d_i \circ \psi^{(1)}_\alpha \circ \pi(y))$. By computations, we obtain the following canonical isomorphism:
\begin{equation*}
\begin{split}
Q^{\otimes -1} \otimes \d R_\alpha
& \cong
(i_1, \ \psi^{(1)}_\alpha \circ \pi)^*(\delta Q \otimes \d P^{\otimes -1}) \\
& \otimes
(i_1 \circ \d_0, \
\psi^{(0)}_{\d_0(\alpha)} \circ \d_0 \circ \pi, \
\d_0 \circ \psi^{(1)}_\alpha \circ \pi)^* \delta P^{\otimes -1} \\
& \otimes
(i_1 \circ \d_1, \
\psi^{(0)}_{\d_1(\alpha)} \circ \d_1 \circ \pi, \
\d_1 \circ \psi^{(1)}_\alpha \circ \pi)^* \delta P \\
& \otimes
\pi^*(\psi^{(1)}_\alpha)^* Q^{\otimes -1} \\
& \otimes
\pi^*(\psi^{(0)}_{\d_0(\alpha)} \circ \d_0, \ 
\d_0 \circ \psi^{(1)}_\alpha)^* P \\
& \otimes
\pi^*(\psi^{(0)}_{\d_1(\alpha)} \circ \d_1, \ 
\d_1 \circ \psi^{(1)}_\alpha)^* P^{\otimes -1}.
\end{split}
\end{equation*}
By the assumption, we can take a section $\tau_\alpha$ of the principal $\T$-bundle 
$$
(\psi^{(1)}_\alpha)^* Q^{\otimes -1}
\otimes
(\psi^{(0)}_{\d_0(\alpha)} \circ \d_0, \ 
\d_0 \circ \psi^{(1)}_\alpha)^* P
\otimes
(\psi^{(0)}_{\d_1(\alpha)} \circ \d_1, \ 
\d_1 \circ \psi^{(1)}_\alpha)^* P^{\otimes -1}
$$
over $U^{(1)}_\alpha$. Now we have a section $r_\alpha : Y_1|_{U^{(1)}_{\alpha}} \to \d R_\alpha \otimes Q^{\otimes -1}$ by setting
\begin{equation*}
\begin{split}
r_\alpha 
& = 
(i_1, \ \psi^{(1)}_\alpha \circ \pi)^*t^{\otimes -1} 
\otimes
(i_1 \circ \d_0, \
\psi^{(0)}_{\d_0(\alpha)} \circ \d_0 \circ \pi, \
\d_0 \circ \psi^{(1)}_\alpha \circ \pi)^* s^{\otimes -1} \\
& \otimes
(i_1 \circ \d_1, \
\psi^{(0)}_{\d_1(\alpha)} \circ \d_1 \circ \pi, \
\d_1 \circ \psi^{(1)}_\alpha \circ \pi)^* s 
\otimes \pi^*\tau_\alpha.
\end{split}
\end{equation*}
To show that $r_\alpha$ satisfies $\delta r_\alpha = \d v_\alpha^{\otimes -1} \otimes t$, we define for $i = 0, 1$ a map
$$
(i_{12} \circ \d_i, \
\psi^{(0)}_{\d_i(\alpha)} \circ \d_i \circ \pi, \
\d_i \circ \psi^{(1)}_\alpha \circ \pi) : \
Y_1^{[2]}|_{U^{(1)}_\alpha} \longrightarrow
Y_0^{[4]}|_{U^{(1)}_{\d_i(\alpha)}}
$$
by $(y_1, y_2) \mapsto (\d_i(y_1), \d_i(y_2), \psi^{(0)}_{\d_i(\alpha)} \circ \d_i \circ \pi(y_1), \d_i \circ \psi^{(1)}_\alpha \circ \pi(y_1))$. Then we have
\begin{multline*}
\pi_j \circ 
(i_{12} \circ \d_i, \
\psi^{(0)}_{\d_i(\alpha)} \circ \d_i \circ \pi, \
\d_i \circ \psi^{(1)}_\alpha \circ \pi) \\
=
\left\{
\begin{array}{lc}
(i_1 \circ \d_i, \
\psi^{(0)}_{\d_i(\alpha)} \circ \d_i \circ \pi, \
\d_i \circ \psi^{(1)}_\alpha \circ \pi) \circ \pi_j, &
(j = 1, 2), \\
(i_{12} \circ \d_i, \ \d_i \circ \psi^{(1)}_\alpha \circ \pi), &
(j = 3), \\
(i_{12} \circ \d_i, \ \psi^{(0)}_{\d_i(\alpha)} \circ \d_i \circ \pi), &
(j = 4).
\end{array}
\right.
\end{multline*}
Because $\delta t = \d s$ holds, computations show
\begin{equation*}
\begin{split}
\delta r_\alpha \otimes \d v_\alpha \otimes t^{\otimes -1} 
& = 
(i_{12} \circ \d_0, \
\psi^{(0)}_{\d_0(\alpha)} \circ \d_0 \circ \pi, \
\d_0 \circ \psi^{(1)}_\alpha \circ \pi)^* \delta s^{\otimes -1} \\
& \otimes 
(i_{12} \circ \d_1, \
\psi^{(0)}_{\d_1(\alpha)} \circ \d_1 \circ \pi, \
\d_1 \circ \psi^{(1)}_\alpha \circ \pi)^* \delta s.
\end{split}
\end{equation*}
Since $\delta s = 1$, we have $\delta r_\alpha = \d v_\alpha^{\otimes -1} \otimes t$.
\end{proof}


\subsection{Construction}

We continue to use the notations in the previous subsection. When local data $\{ (R_\alpha, v_\alpha), \eta_\alpha, w_{\alpha \beta}, r_\alpha \}$ are given, we define a cochain
\begin{equation}
\left(
\begin{array}{c}
f_{\alpha \beta \gamma} \\
\theta^1_{\alpha \beta}, \
g_{\alpha \beta} \\
\theta^2_{\alpha}, \quad \
\omega^1_{\alpha}, \quad \
h_{\alpha} 
\end{array}
\right)
\in C^2(\U^{\bullet}, \bF(2)) = 
\begin{array}{c}
K^{0, 2, 0} \oplus \\
K^{0, 1, 1} \oplus 
K^{1, 1, 0} \oplus \\
K^{0, 0, 2} \oplus 
K^{1, 0, 1} \oplus
K^{2, 0, 0}
\end{array}
\label{cochain:degree2}
\end{equation}
by the formulae
\begin{align}
\pi^*f_{\alpha \beta \gamma}
& = 
w_{\beta \gamma} \otimes
w_{\alpha \gamma}^{\otimes -1} \otimes
w_{\alpha \beta}, 
\label{cocycle:f} \\
\pi^*g_{\alpha \beta}
& = 
\d w_{\alpha \beta} \otimes r_\alpha^{\otimes -1} \otimes r_\beta, 
\label{cocycle:g} \\
\pi^*h_{\alpha} 
& = 
u^{\otimes -1} \otimes \d r_\alpha^{\otimes -1}, 
\label{cocycle:h} \\
\pi^*\theta^1_{\alpha \beta}
& = 
\frac{-1}{2\pi\im}
w_{\alpha \beta}^* (\eta_\alpha \otimes \eta_\beta^{\otimes -1}), 
\label{cocycle:theta1} \\
\pi^*\omega^1_{\alpha}
& = 
\frac{1}{2\pi\im}
r_\alpha^*([\d \eta_\alpha]_{rel} \otimes D_{rel}^{\otimes -1}), \\
\pi^*\theta^2_{\alpha} 
& =  
\frac{1}{2\pi\im} (F(\eta_\alpha) - f).
\label{cocycle:theta2}
\end{align}
It is easy to verify that the definition above is well-defined and (\ref{cochain:degree2}) is a cocycle. 

\begin{lem} \label{lem:indep_local_data}
Let us fix an open cover $\U^\bullet = \{ \U^{(p)} \}_{p \ge 0}$ of $G^\bullet \times M$. Suppose that $\U^{(0)}$ is a good cover of $M$ and that there exists a local section $\psi_\alpha : U^{(0)}_\alpha \to Y_0|_{U^{(0)}_\alpha}$ for each $\alpha \in \AA^{(0)}$. For $(\G_G, \nabla_G, f)$, the class in $H^2(\U^{\bullet}, \bF(2))$ represented by (\ref{cochain:degree2}) is independent of the choices of local data $\{ (R_\alpha, v_\alpha), \eta_\alpha, w_{\alpha \beta},r_\alpha \}$.
\end{lem}

\begin{proof}
Let $\{ (R'_\alpha, v'_\alpha), \eta'_\alpha, w'_{\alpha \beta}, r'_\alpha \}$ be the other data. First of all, we construct a section $\rho_\alpha : Y_0|_{U^{(0)}_\alpha} \to R_\alpha \otimes {R'}_\alpha^{\otimes -1}$ such that $\delta \rho_\alpha = v_\alpha^{\otimes -1} \otimes v'_\alpha$ as follows. Recall the map used in the proof of Lemma \ref{lem:construct_local_data}:
$$
(i_{12 \cdots p}, \ \psi_\alpha \circ \pi) : 
Y_0|_{U^{(p)}_\alpha} \longrightarrow Y_0^{[p+1]}|_{U^{(0)}_\alpha}.
$$
There is a canonical isomorphism
$$
(i_1, \psi_\alpha \circ \pi)^* 
\delta(R_\alpha \otimes {R'}_\alpha^{\otimes -1}) 
=
\pi^*(\psi_\alpha^*(R_\alpha^{\otimes -1} \otimes R'_\alpha)) \otimes
(R_\alpha \otimes {R'}_\alpha^{\otimes -1}).
$$
Since $U^{(0)}_\alpha$ is contractible by the assumption, there is a section $\varrho_\alpha : U^{(0)}_\alpha \to \psi_\alpha^*(R_\alpha^{\otimes -1} \otimes R'_\alpha)$. If we define $\rho_\alpha$ by $\rho_\alpha = \pi^*\varrho_\alpha \otimes (i_1, \psi_\alpha \circ \pi)^*(v_\alpha^{\otimes -1} \otimes v_\alpha)^{\otimes -1}$, then we have $\delta \rho_\alpha = v_\alpha^{\otimes -1} \otimes v_\alpha$.

By means of the section $\rho_\alpha$, we can define a cochain
\begin{equation}
\left(
\begin{array}{c}
k_{\alpha \beta} \\
\epsilon_{\alpha}, \quad \ \l_{\alpha}
\end{array}
\right)
\in C^1(\U^{\bullet}, \bF(2)) = 
\begin{array}{c}
K^{0, 1, 0} \oplus \\
K^{0, 0, 1} \oplus 
K^{1, 0, 0}
\end{array}
\end{equation}
by setting
\begin{align*}
\pi^*k_{\alpha \beta} 
& = 
\rho_\alpha \otimes \rho_\beta^{\otimes -1} \otimes
w_{\alpha \beta}^{\otimes -1} \otimes w'_{\alpha \beta}, \\
\pi^* \l_\alpha
& = 
r_\alpha \otimes {r'}_{\alpha}^{\otimes -1} \otimes 
\d \rho_\alpha^{\otimes -1}, \\
\pi^*\epsilon_\alpha 
& = 
\frac{-1}{2\pi\im} 
\rho_\alpha^*(\eta_\alpha \otimes {\eta'}_\alpha^{\otimes -1}).
\end{align*}
Then we can prove that
$$
\left(
\begin{array}{c}
f'_{\alpha \beta \gamma} \\
{\theta'}^1_{\alpha \beta}, \
g'_{\alpha \beta} \\
{\theta'}^2_{\alpha}, \quad \
{\omega'}^1_{\alpha}, \quad \
h'_{\alpha} 
\end{array}
\right)
=
\left(
\begin{array}{c}
f_{\alpha \beta \gamma} \\
\theta^1_{\alpha \beta}, \
g_{\alpha \beta} \\
\theta^2_{\alpha}, \quad \
\omega^1_{\alpha}, \quad \
h_{\alpha} 
\end{array}
\right)
+ D
\left(
\begin{array}{c}
k_{\alpha \beta} \\
\epsilon_{\alpha}, \quad \ \l_{\alpha}
\end{array}
\right),
$$
where the left hand side is the cocycle given by $\{ (R'_\alpha, v'_\alpha), \eta'_\alpha, w'_{\alpha \beta}, r'_\alpha\}$. Thus, the resulting cohomology class in $H^2(\U^{\bullet}, \bF(2))$ is independent of the choice of local data.
\end{proof}

\begin{prop} \label{prop:char_class_ebgcc}
To a $G$-equivariant bundle gerbe with connection and curving $(\G_G, \nabla_G, f)$ over $M$ we can assign a well-defined cohomology class 
$$
\delta_G(\G_G, \nabla_G, f) \in H^2(G^{\bullet} \times M, \bF(2)).
$$
\end{prop}

\begin{proof}
Note that, for a given open cover of $G^{\bullet} \times M$, we can take a refinement $\U^{\bullet}$ such that $\U^{(p)}$ is good for $p = 0, \ldots, 3$. For such an open cover $\U^{\bullet}$, we have $H^2(G^{\bullet} \times M, \bF(2)) \cong H^2(\U^{\bullet}, \bF(2))$. By Lemma \ref{lem:construct_local_data} we can take local data, and obtain the cocycle (\ref{cochain:degree2}). The cocycle defines a cohomology class $\delta(\U^{\bullet}) \in H^2(\U^{\bullet}, \bF(2))$, which is independent of the choice of local data by Lemma \ref{lem:indep_local_data}.

Let ${\U'}^{\bullet}$ be the other open cover such that ${\U'}^{(p)}$ is good for $p = 0, \ldots, 3$. Let $\V^{\bullet}$ be a common refinement of $\U^{\bullet}$ and ${\U'}^{\bullet}$ such that $\V^{(0)}$ is good. We denote the induced homomorphisms by $\phi : H^2(\U^{\bullet}, \bF(2)) \to H^2(\V^{\bullet}, \bF(2))$ and $\phi' : H^2({\U'}^{\bullet}, \bF(2)) \to H^2(\V^{\bullet}, \bF(2))$. The local data $\{ (R_\alpha, v_\alpha), \eta_\alpha, w_{\alpha \beta}, r_\alpha \}$ on $\U^{\bullet}$ induce those on $\V^{\bullet}$ which give $\phi(\delta(\U^{\bullet}))$. Similarly, we obtain local data on $\V^{\bullet}$ which give $\phi'(\delta({\U'}^{\bullet}))$. Because $\V^{(0)}$ is good, we have $\phi(\delta(\U^{\bullet})) = \phi'(\delta({\U'}^{\bullet}))$ by Lemma \ref{lem:indep_local_data}. This implies that, in $H^2(G^{\bullet} \times M, \bF(2))$, the class $\delta(\U^{\bullet})$ is independent of the choice of the open cover $\U^{\bullet}$. By setting $\delta_G(\G_G, \nabla_G, f) = \delta(\U^{\bullet})$, we complete the proof of this proposition.
\end{proof}

\begin{prop} \label{prop:formulae_equiv_Deligne_class}
Let $(\G_G, \nabla_G, f)$ and $(\G'_G, \nabla'_G, f')$ be $G$-equivariant bundle gerbes with connection and curving over $M$. We have
\begin{align}
\delta_G( (\G_G, \nabla_G, f)^{\otimes -1} )
& = 
- \delta_G(\G_G, \nabla_G, f), \\
\delta_G( (\G_G, \nabla_G, f) \otimes (\G'_G, \nabla'_G, f') ) 
& = 
\delta_G(\G_G, \nabla_G, f) + 
\delta_G(\G'_G, \nabla'_G, f').
\end{align}
\end{prop}

\begin{proof}
This proposition can be shown by fixing a suitable open cover $\U^{\bullet}$ of $G^{\bullet} \times M$. Let $\{ (R_\alpha, v_\alpha), \eta_\alpha, w_{\alpha \beta}, r_\alpha \}$ be local data on $\U^{\bullet}$ which give the class $\delta_G(\G_G, \nabla_G, f)$. Then $\{ (R_\alpha^{\otimes -1}, v_\alpha^{\otimes -1}), \eta_\alpha^{\otimes -1}, w_{\alpha \beta}^{\otimes -1}, r_\alpha^{\otimes -1} \}$ give rise to local data for $(\G_G, \nabla_G, f)^{\otimes -1}$. It is direct to see that the data for $(\G_G, \nabla_G, f)^{\otimes -1}$ give $ - \delta_G(\G_G, \nabla_G, f)$. Let $\{ (R'_\alpha, v'_\alpha), \eta'_\alpha, w'_{\alpha \beta}, r'_\alpha \}$ be local data on $\U^{\bullet}$ which give $\delta_G(\G'_G, \nabla'_G, f')$. Then $\{ (R_\alpha \otimes R'_\alpha, v_\alpha \otimes v'_\alpha), \eta_\alpha \otimes \eta'_\alpha, w_{\alpha \beta} \otimes w'_{\alpha \beta}, r_\alpha \otimes r'_\alpha \}$ are local data for $(\G_G, \nabla_G, f) \otimes (\G'_G, \nabla'_G, f')$. The data give the class $\delta_G(\G_G, \nabla_G, f) + \delta_G(\G'_G, \nabla'_G, f')$.
\end{proof}

\begin{prop} \label{prop:trivial_ebgcc}
A $G$-equivariant bundle gerbe with connection and curving $(\G_G, \nabla_G, f)$ over $M$ is trivial if and only if $\delta_G(\G_G, \nabla_G, f) = 0$.
\end{prop}

\begin{proof}
If $(\G_G, \nabla_G, f)$ is trivial, then we have a pseudo $\T$-bundle $((R, v), r)$ with a connection $\eta$ such that $F(\eta) = f$. By restriction, we obtain local data $(R_\alpha, v_\alpha)$, $\eta_\alpha$ and $r_\alpha$. We put $w_{\alpha \beta} = 1$. Now it is clear that the associated cocycle (\ref{cochain:degree2}) is zero. Conversely, if $\delta_G(\G_G, \nabla_G, f)$ vanishes, then (\ref{cochain:degree2}) is written as a coboundary. By using the coboundary, we can glue the local data, which define (\ref{cochain:degree2}), together to construct globally a pseudo $\T$-bundle with a connection whose curvature coincides with the curving $f$.
\end{proof}

Proposition \ref{prop:formulae_equiv_Deligne_class} and Proposition \ref{prop:trivial_ebgcc} guarantee that the notion of stable isomorphisms is an equivalence relation on equivariant bundle gerbes with connection and curving. The propositions also imply:

\begin{cor} 
Suppose that $(\G_G, \nabla_G, f)$ and $(\G'_G, \nabla'_G, f')$ are $G$-equivariant bundle gerbes with connection and curving over $M$. They are stably isomorphic if and only if $\delta_G(\G_G, \nabla_G, f) = \delta_G(\G'_G, \nabla'_G, f')$.
\end{cor}

\medskip

Now we return to strongly equivariant bundle gerbes. As we see in Lemma \ref{lem:sEBG_to_EBG}, a strongly $G$-equivariant bundle gerbe with connection and curving $(\G, \nabla, f)$ induces a $G$-equivariant bundle gerbe with connection and curving $(\G_G, \nabla_G, f)$. We define $\delta_G(\G, \nabla, f) \in H^2(G^{\bullet} \times M, \bF(2))$ by $\delta_G(\G, \nabla, f) = \delta_G(\G_G, \nabla_G, f)$. As a corollary to Proposition \ref{prop:formulae_equiv_Deligne_class}--\ref{prop:trivial_ebgcc}, we summarize the properties of the class.

\begin{cor} \label{cor:properties_char_class_sEBGcc}
Let $(\G, \nabla, f)$ and $(\G', \nabla', f')$ be strongly $G$-equivariant bundle gerbes with connection and curving over $M$.

(a) We have
\begin{align*}
\delta_G( (\G, \nabla, f)^{\otimes -1} )
& = 
- \delta_G(\G, \nabla, f), \\
\delta_G( (\G, \nabla, f) \otimes (\G', \nabla', f') ) 
& = 
\delta_G(\G, \nabla, f) + 
\delta_G(\G', \nabla', f').
\end{align*}

(b) $(\G, \nabla, f)$ is trivial if and only if $\delta_G(\G, \nabla, f) = 0$.

(c) $(\G, \nabla, f)$ and $(\G', \nabla', f')$ are stably isomorphic if and only if 
$$\delta_G(\G, \nabla, f) = \delta_G(\G', \nabla', f').$$
\end{cor}

By using the last property, we prove some lemmas for later convenience. An \textit{isomorphism} of strongly $G$-equivariant bundle gerbes with connection and curving is defined to be an isomorphism of the underlying bundle gerbes with connection and curving $(\varphi, \til{\varphi}) : (\G, \nabla, f) \to (\G', \nabla', f')$ such that $\varphi$ and $\til{\varphi}$ are $G$-equivariant.

\begin{lem} \label{lem:sEBGcc_iso}
If $(\G, \nabla, f)$ and $(\G', \nabla', f')$ are isomorphic, then they are stably isomorphic.
\end{lem}

\begin{proof}
Let $(\G'_G, \nabla_G', f')$ be induced from $(\G', \nabla', f')$, $\U^{\bullet}$ a sufficiently fine open cover of $G^{\bullet} \times M$, and $\{ (R'_\alpha, v'_\alpha), \eta'_\alpha, w'_{\alpha \beta}, r'_\alpha \}$ local data for $(\G'_G, \nabla_G', f')$ giving the cohomology class $\delta_G(\G', \nabla', f')$. Because $\varphi^{[p]}$ and $\til{\varphi}$ are $G$-equivariant, $\{ \varphi^*(R'_\alpha, v'_\alpha), \varphi^*\eta'_\alpha, \varphi^*w'_{\alpha \beta}, \varphi_1^*r'_\alpha \}$ give rise to local data for $(\G_G, \nabla_G, f)$, where $\varphi_1 : G \times Y \to G \times Y'$ is given by $\varphi_1 = \id \times \varphi$. If we use the local data to give $\delta_G(\G, \nabla, f)$, then we have $\delta_G(\G, \nabla, f) = \delta_G(\G', \nabla', f')$.
\end{proof}

Let $\kappa \in \im A^1(Y)$ be a $G$-invariant 1-form on $Y$. Then the connection $\nabla - \delta \kappa$ on $P$ gives rise to a $G$-invariant connection on $\G$. Similarly, the 2-form $f - d\kappa$ gives rise to a $G$-invariant curving for $\nabla - \delta \kappa$.

\begin{lem} \label{lem:shift_connection_and_curving}
The strongly $G$-equivariant bundle gerbes with connection and curving $(\G, \nabla, f)$ and $(\G, \nabla - \delta \kappa, f - d\kappa)$ are stably isomorphic.
\end{lem}

\begin{proof}
Let $\U^{\bullet}$ be a sufficiently fine open cover of $G^{\bullet} \times M$. We take local data $\{ (R_\alpha, v_\alpha), \eta_\alpha, w_{\alpha \beta}, r_\alpha \}$ for $(\G_G, \nabla_G, f)$ induced from $(\G, \nabla, f)$. We also take a local section $\psi_\alpha : U^{(0)}_\alpha \to Y|_{U^{(0)}_\alpha}$ and put $\kappa_\alpha = \psi_\alpha^*\kappa$. Since $\kappa$ is $G$-invariant, $\{ (R_\alpha, v_\alpha), \eta_\alpha - \pi^*\kappa_\alpha, w_{\alpha \beta}, r_\alpha \}$ give rise to local data for the equivariant bundle gerbe with connection and curving induced from $(\G, \nabla - \delta \kappa, f - d\kappa)$. We can verify that the \Cech cocycles defined by these local data are the same. 
\end{proof}

\begin{rem}
In Section \ref{sec:bundle_gerbes}, we mentioned the Dixmier-Douady class associated to a bundle gerbe $\G$ over $M$. Suppose that $\U$ is a good cover of $M$. Then the cocycle $(f_{\alpha \beta \gamma}) \in Z^2(\U, \u{\T})$ defined by (\ref{cocycle:f}) represents the Dixmier-Douady class $\delta(\G) \in H^2(M, \u{\T}) \cong H^3(M, \Z)$. Similarly, the cocycle $(f_{\alpha \beta \gamma}, \theta^1_{\alpha \beta}, \theta^2_\alpha) \in Z^2(\U, \F(2))$ defined by (\ref{cocycle:f}), (\ref{cocycle:theta1}) and (\ref{cocycle:theta2}) represents the class $\delta(\G, \nabla, f)$ associated to a bundle gerbe with connection and curving $(\G, \nabla, f)$ over $M$.
\end{rem}


\subsection{Classification}

Here we classify $G$-equivariant bundle gerbes with connection and curving by using the cohomology class $\delta_G(\G_G, \nabla_G, f)$. However, the classification is not necessary in the sequel, so the reader can skip this subsection.

\begin{thm} 
The group of stable isomorphism classes of $G$-equivariant bundle gerbes with connection and curving over $M$ is isomorphic to the cohomology group $H^2(G^{\bullet} \times M, \bF(2))$.
\end{thm}

\begin{proof}
By Proposition \ref{prop:formulae_equiv_Deligne_class} and Proposition \ref{prop:trivial_ebgcc}, the assignment 
$$
(\G_G, \nabla_G, f) \mapsto \delta_G(\G_G, \nabla_G, f)
$$ 
induces a monomorphism from the group of stable isomorphism classes of $G$-equivariant bundle gerbes with connection and curving over $M$ to the cohomology group $H^2(G^{\bullet} \times M, \bF(2))$. To prove that the monomorphism is surjective, we take and fix a cocycle
\begin{equation}
\left(
\begin{array}{c}
f_{\alpha \beta \gamma} \\
\theta^1_{\alpha \beta}, \
g_{\alpha \beta} \\
\theta^2_{\alpha}, \quad \
\omega^1_{\alpha}, \quad \
h_{\alpha} 
\end{array}
\right)
\in Z^2(\U^{\bullet}, \bF(2)), \label{cocycle:given}
\end{equation}
where $\U^{\bullet} = \{ \U^{(p)} \}_{p \ge 0}$ is an open cover of $G^{\bullet} \times M$. 

\smallskip

First, we construct a bundle gerbe over $M$ as follows. Using the open cover $\U^{(p)} = \{ U^{(p)} \}_{\alpha \in \AA^{(p)}}$, we put $Y_p = \bigsqcup_{\alpha \in \AA^{(p)}} U^{(p)}_{\alpha}$. By virtue of the face and degeneracy maps of $G^p \times M$, the sequence $Y_{\bullet} = \{ Y_p \}_{p \ge 0}$ gives rise to a simplicial manifold. Apparently, we have a simplicial surjective submersion $\pi : Y_{\bullet} \to G^{\bullet} \times M$ admitting local sections. We put $P = Y_0^{[2]} \times \T$. Note that $\delta P = Y_0^{[3]} \times \T$ and $Y_0^{[3]} = \bigsqcup_{\alpha, \beta, \gamma \in \AA^{(0)}} U^{(0)}_{\alpha \beta \gamma}$. We denote by $s = \bigsqcup_{\alpha, \beta, \gamma \in \AA^{(0)}} (\id \times f_{\alpha \beta \gamma}^{-1})$ the section $s : Y_0^{[3]} \to \delta P$ that carries $x \in U^{(0)}_{\alpha \beta \gamma}$ to $(x, f^{-1}_{\alpha \beta \gamma}(x)) \in U^{(0)}_{\alpha \beta \gamma} \times \T$. Then $\G = (Y_0, P, s)$ is a bundle gerbe over $M$.

\smallskip

Second, we make this bundle gerbe into $G$-equivariant. We define a principal $\T$-bundle $Q \to Y_1$ by $Q = Y_1 \times \T$. Since $Y_1^{[2]} = \bigsqcup_{\alpha, \beta \in \AA^{(1)}} U^{(1)}_{\alpha \beta}$ and $\delta Q^{\otimes -1} \otimes \d P = Y_1^{[2]} \times \T$, we define a section $t : Y_1^{[2]} \to \delta Q^{\otimes -1} \otimes \d P$ by $t = \bigsqcup_{\alpha, \beta \in \AA^{(1)}} (\id \times g_{\alpha \beta}^{-1})$. Similarly, because $\d Q = Y_2 \times \T$, we define a section $u : Y_2 \to \d Q$ by $u = \bigsqcup_{\alpha \in \AA^{(2)}} (\id \times h_\alpha^{-1})$. Then $\G_G = (Y_{\bullet}, (Y_0, P, s), (Q, t), u)$ is a $G$-equivariant bundle gerbe over $M$.

\smallskip

Third, we give $\G_G$ a $G$-invariant connection and a $G$-invariant curving. A connection $\nabla$ on $(Y_0, P, s)$ is given by $\nabla = \bigsqcup_{\alpha, \beta \in \AA^{(0)}} (-2\pi\im\theta^1_{\alpha \beta} + u^{-1}du)$, where $u^{-1}du$ stands for the Maurer-Cartan form on $\T$. If we define a relative connection $D_{rel}$ on $Q$ by $D_{rel} = \bigsqcup_{\alpha \in \AA^{(1)}} (-2\pi\im\omega^1_\alpha + [u^{-1}du]_{rel})$, then $\nabla_G = (\nabla, D_{rel})$ is a $G$-invariant connection on $\G_G$. A $G$-invariant curving $f$ for $\nabla_G$ is defined by $f = \bigsqcup_{\alpha \in \AA^{(0)}} (-2\pi\im\theta^2_\alpha)$.

\smallskip

Finally, we take local data $\{ (R_\alpha, v_\alpha), \eta_\alpha, w_{\alpha \beta}, r_\alpha \}$ for $(\G_G, \nabla_G, f)$ with respect to the open cover $\U^{\bullet}$ as follows. 
\begin{align*}
R_\alpha 
& =  
Y_0|_{U_\alpha} \times \T = 
\bigsqcup_{i \in \AA^{(0)}} U_{i \alpha}^{(0)} \times \T, \\
v_\alpha 
& =  
\bigsqcup_{i, j \in \AA^{(0)}} (\id \times f_{i j \alpha}^{-1}), \\
\eta_\alpha
& = 
\bigsqcup_{i \in \AA^{(0)}} (2\pi\im \theta^1_{i \alpha} + u^{-1}du), \\
w_{\alpha \beta} 
& = 
\bigsqcup_{i \in \AA^{(0)}} (\id \times f_{i \alpha \beta}), \\
r_{\alpha}
& = 
\bigsqcup_{i \in \AA^{(1)}} (\id \times g_{i \alpha}).
\end{align*}
By direct computations, we can see that the cocycle defined by the local data above coincides with the given cocycle (\ref{cocycle:given}).
\end{proof}

We note that the cohomology $H^2(G^{\bullet} \times M, \bF(2))$ also classifies the isomorphism classes of $G$-equivariant gerbes with connection and curvings \cite{Go1}

Forgetting connections and curvings, we can obtain the following classification of $G$-equivariant bundle gerbes.

\begin{cor} 
The group of stable isomorphism classes of $G$-equivariant bundle gerbes over $M$ is isomorphic to $H^2(G^{\bullet} \times M, \u{\T}_{G^{\bullet} \times M})$.
\end{cor}

The isomorphism in the corollary above is induced by the assignment to $\G_G$ of the class $\delta_G (\G_G) \in H^2(G^{\bullet} \times M, \u{\T})$ represented by the cocycle $(f_{\alpha \beta \gamma}, g_{\alpha \beta}, h_\alpha) \in Z^2(\U^{\bullet}, \u{\T})$ determined by (\ref{cocycle:f}), (\ref{cocycle:g}) and (\ref{cocycle:h}).

If $G$ is compact, then $H^m(G^{\bullet} \times M, \u{\T}_{G^{\bullet} \times M})$ is isomorphic to the \textit{equivariant cohomology group} (\cite{A-B}) $H^{m+1}_G(M, \Z)$ for $m > 1$. Thus, in this case, the group of stable isomorphism classes of $G$-equivariant bundle gerbes over $M$ is isomorphic to $H^3_G(M, \Z)$.


\section{Obstruction}
\label{sec:obstruction}

In this section, we prove Theorem \ref{ithm:obstruction}: the element $\beta(\G, \nabla, f) \in \mathcal{Z}/\mathcal{B}$ given in Section \ref{sec:introduction} is the obstruction to the reduction of a strongly equivariant bundle gerbe with connection and curving. For the purpose, we relate $\beta_G(\G, \nabla, f)$ with the characteristic class $\delta_G(\G, \nabla, f)$. Then the results on equivariant smooth Deligne cohomology groups lead to the proof.

Throughout this section, $G$ denotes a Lie group acting on a smooth manifold $M$ by left. 


\subsection{Relation to the characteristic class}

To begin with, we recall the definition of $\beta_G(\G, \nabla, f) \in \mathcal{Z} / \mathcal{B}$.

Apart from the dual space $\g^* = \Hom(\g, \R)$ of $\g$, we introduce a similar vector space $\g^\dagger = \Hom(\g, \im\R)$. This will help to suppress notations in the below. We usually identify these vector spaces thorough the map $\g^\dagger \to \g^*$ given by $z \mapsto \frac{-1}{2\pi\im}z$. We also denote by $\langle \ | \ \rangle : \g \otimes \g^\dagger \to \im\R$ the natural contraction between $\g$ and $\g^\dagger$. 

Recall the \textit{moment} \cite{B-V} associated to an equivariant principal $\T$-bundle with connection. Let $\mathcal{P} \to M$ be a $G$-equivariant principal $\T$-bundle equipped with a $G$-invariant connection $\Theta$. The moment associated to $(\mathcal{P}, \Theta)$ is a function $\mu : M \to \g^\dagger$ defined by $\langle X | \mu(x) \rangle = \Theta(p; X^*)$, where $p \in \mathcal{P}$ is a point lying on the fiber of $x \in M$ and $X^* \in T_p\mathcal{P}$ is the tangent vector generated by the infinitesimal action of $X \in \g$. It is easy to verify that $\langle X | d \mu \rangle + \iota_{X^*}F(\Theta) = 0$ and $g^*\mu = \Ad_g \mu$ for $X \in \g$ and $g \in G$.

\begin{lem}[\cite{Ma-S}] \label{lem:lambda}
Let $(Y, P, s)$ be a strongly $G$-equivariant bundle gerbe over $M$, $\nabla$ a $G$-invariant connection on $(Y, P, s)$, and $\til{\lambda} : Y^{[2]} \to \g^\dagger$ the moment associated to $(P, \nabla)$. Then the space $\{ \lambda : Y \to \g^\dagger |\ \delta \lambda = \til{\lambda} \}$ is non-empty and forms an affine space under $A^0(M, \g^\dagger)$.
\end{lem}

\begin{proof}
Since $\g^\dagger = \Hom(\g, \im\R)$ is a vector space, Lemma \ref{lem:Murray} provides the following exact sequence:
$$
0 \to 
A^0(M, \g^\dagger) \stackrel{\pi^*}{\to}
A^0(Y, \g^\dagger) \stackrel{\delta}{\to}
A^0(Y^{[2]}, \g^\dagger) \stackrel{\delta}{\to}
A^0(Y^{[3]}, \g^\dagger) \to \cdots.
$$
A computation shows that $\til{\lambda} \in A^0(Y^{[2]}, \g^\dagger)$ satisfies $\delta \til{\lambda} = 0$. Hence the lemma follows from the above exact sequence. 
\end{proof}

\begin{dfn} 
Let $\G = (Y, P, s)$ be a strongly $G$-equivariant bundle gerbe over $M$ equipped with a $G$-invariant connection $\nabla$ and a $G$-invariant curving $f$, and $\til{\lambda} : Y^{[2]} \to \g^\dagger$ the moment associated to $(P, \nabla)$. Taking a map $\lambda : Y \to \g^\dagger$ such that $\delta \lambda = \til{\lambda}$, we define $(E, \zeta) \in \mathcal{C} = A^1(M, \g^*) \oplus C^{\infty}(G, A^0(M, \g^*))$ by
\begin{align}
\langle X | \pi^*E \rangle 
& = 
\frac{1}{2\pi\im}
\left( \langle X | d \lambda \rangle + \iota_{X^*} f \right), 
\label{formula:dfn_E_by_lambda} \\
\langle X | \pi^*\zeta(g) \rangle
& = 
\frac{1}{2\pi\im}
\langle X | g^*\lambda - \Ad_g \lambda \rangle,
\label{formula:dfn_zeta_by_lambda}
\end{align}
where $\iota_{X^*}$ is the contraction by the tangent vector field on $Y$ generated by the infinitesimal action of $X \in \g$. 
\end{dfn}

Once $\lambda : Y \to \g^\dagger$ is chosen, we can verify the following formulae for $(E, \zeta)$:
\begin{align*}
- \frac{1}{2\pi\im} \iota_{X^*}\Omega &= \langle X | d E \rangle, \\
g^*E - \Ad_gE &= d \zeta(g), \\
\Ad_g\zeta(h) - \zeta(gh) + h^*\zeta(g) &= 0,
\end{align*}
where $\Omega$ is the 3-curvature of $(\G, \nabla, f)$. Hence $(E, \zeta)$ belongs to $\mathcal{Z}$. We can also verify that $[E, \zeta] \in \mathcal{Z} / \mathcal{B}$ is independent of the choice of $\lambda : Y \to \g^\dagger$.

\begin{dfn}
We define $\beta_G(\G, \nabla, f) \in \mathcal{Z}/\mathcal{B}$ by $\beta_G(\G, \nabla, f) = [E, \zeta]$.
\end{dfn}

\medskip

Now we recall the connecting homomorphism induced by (\ref{exact_seq:simplicial_complex}):
$$
\beta : \
H^2(G^{\bullet} \times M, \bF(2)) \longrightarrow
H^3(G^{\bullet} \times M, F^1\!\F(2)). 
$$
Since the characteristic class $\delta_G(\G, \nabla, f)$ belongs to $H^2(G^{\bullet} \times M, \bF(2))$, we have $\beta(\delta_G(\G, \nabla, f)) \in H^3(G^{\bullet} \times M, F^1\!\F(2))$.

\begin{prop} \label{prop:expression_beta}
Let $(\G, \nabla, f)$ be a strongly $G$-equivariant bundle gerbe with connection and curving over $M$. Under the isomorphism in Lemma \ref{lem:iso_H3F1F2}, we have $\beta(\delta_G(\G, \nabla, f)) = \beta_G(\G, \nabla, f)$.
\end{prop}

To prove Proposition \ref{prop:expression_beta}, we need some preparations. Recall that a strongly $G$-equivariant bundle gerbe $(Y, P, s)$ over $M$ defines a $G$-equivariant bundle gerbe $(G^{\bullet} \times Y, (Y, P, s), (Q, t), u)$. In particular, we have a section $t : G \times Y^{[2]} \to \d P$ such that $\d t = 1$. If $\nabla$ is a $G$-invariant connection on $(Y, P, s)$, then we have $t^*[\d \nabla]_{rel} = 0$ in $\im A^1(G \times Y^{[2]})_{rel}$, that is, $t^*(\d \nabla) \in \im F^1\!A^1(G \times Y^{[2]})$. 

\begin{lem}[\cite{Ma-S}] \label{lem_exist_tau}
There is $\tau \in \im F^1\!A^1(G \times Y)$ such that $\delta \tau = t^*(\d \nabla)$.
\end{lem}

\begin{proof}
Note that $\delta (t^* \d \nabla) = 0$. Thus, by means of Lemma \ref{lem:Murray}, there is a 1-form $\tau \in \im A^1(G \times Y)$ such that $\delta \tau = t^*(\d \nabla)$. Because $t^*(\d \nabla)$ belongs to $\im F^1\!A^1(G \times Y^{[2]})$, the 1-form $\tau$ belongs to $\im F^1\!A^1(G \times Y)$.
\end{proof}

\begin{lem} \label{lem:surjection_tau_to_lambda}
There exists a surjection
$$
\varpi : 
\{ \tau \in \im F^1\!A^1(G \times Y)|\ \delta \tau = t^*(\d \nabla) \} 
\longrightarrow
\{ \lambda \in A^0(Y, \g^\dagger)|\ \delta \lambda = \til{\lambda} \}.
$$
\end{lem}

\begin{proof}
We define a map $\varpi : \im F^1\!A^1(G \times Y) \to A^0(Y, \g^\dagger)$ by $\langle X | \varpi(\tau)(y) \rangle = - \tau((e, y); X \oplus 0)$. Notice the following formula:
$$
\left(t^*(\d \nabla) \right) \left((g, (y_1, y_2)); gX \oplus W \right) = 
- \langle X | g^*\til{\lambda}(y_1, y_2) \rangle,
$$
where $W \in T_{(y_1, y_2)}Y^{[2]}$ is a tangent vector. Thus, if $\tau \in \im F^1\!A^1(G \times Y)$ satisfies $\delta \tau = t^*(\d \nabla)$, then $\delta \varpi(\tau) = \til{\lambda}$. Suppose that $\lambda \in A^0(Y, \g^\dagger)$ satisfies $\delta \lambda = \til{\lambda}$. If we define $\tau \in \im F^1\!A^1(G \times Y)$ by $\tau((g, y); gX \oplus V) = - \langle X | \lambda(gy) \rangle$, then $\delta \tau = t^*(\d \nabla)$ and $\varpi(\tau) = \lambda$.
\end{proof}

\begin{proof}[The proof of Proposition \ref{prop:expression_beta}]
First, let $\U^{\bullet}$ be a sufficiently fine open cover of $G^{\bullet} \times M$. For the $G$-equivariant bundle gerbe with connection and curving induced from the strongly $G$-equivariant bundle gerbe $(\G, \nabla, f)$, we take local data $\{ (R_\alpha, v_\alpha), \eta_\alpha, w_{\alpha \beta}, r_\alpha \}$ and define 
\begin{equation}
\left(
\begin{array}{c}
f_{\alpha \beta \gamma} \\
\theta^1_{\alpha \beta}, \
g_{\alpha \beta} \\
\theta^2_{\alpha}, \quad \
\omega^1_{\alpha}, \quad \
h_{\alpha} 
\end{array}
\right)
\in Z^2(\U^{\bullet}, \bF(2))
\label{cocycle:thm_expression_beta}
\end{equation}
by (\ref{cocycle:f}) -- (\ref{cocycle:theta2}). 

Second, we lift the cocycle (\ref{cocycle:thm_expression_beta}) in $Z^2(\U^{\bullet}, \bF(2))$ to a cochain in $C^2(\U^{\bullet}, \F(2))$. For this aim, we note that $Q$ and $D_{rel}$ are trivial. Thus $\omega^1_\alpha \in A^1(U^{(1)}_\alpha)_{rel}$ is
$$
\pi^*\omega^1_\alpha  = \frac{1}{2\pi\im} r_\alpha^*[\d \eta_\alpha]_{rel}.
$$
By the help of Lemma \ref{lem_exist_tau}, we can take a 1-form $\tau \in \im F^1\!A^1(G \times Y)$ such that $\delta \tau = t^*(\d \nabla)$. Using this 1-form, we define $\til{\omega}^1_\alpha \in A^1(U^{(1)}_\alpha)$ by
$$
\pi^*\til{\omega}^1_\alpha  = \frac{1}{2\pi\im}
(r_\alpha^*(\d \eta_\alpha) - \tau).
$$
Clearly, $[\til{\omega}^1_\alpha]_{rel} = \omega^1_\alpha$. Thus we have a lift of the cocycle (\ref{cocycle:thm_expression_beta}):
\begin{equation}
\left(
\begin{array}{c}
f_{\alpha \beta \gamma} \\
\theta^1_{\alpha \beta}, \
g_{\alpha \beta} \\
\theta^2_{\alpha}, \quad \
\til{\omega}^1_{\alpha}, \quad \
h_{\alpha} 
\end{array}
\right)
\in C^2(\U^{\bullet}, \F(2)).
\label{cochain:thm_expression_beta}
\end{equation}

Third, we compute the coboundary of the cochain (\ref{cochain:thm_expression_beta}). By direct computations, we can prove the following formulae:
\begin{align*}
\d \theta^1_{\alpha \beta} 
+ \frac{1}{2\pi\im} d\log g_{\alpha \beta}
- \delta \til{\omega}^1 
& = 
0, \\
\d \theta^2_\alpha 
- d \til{\omega}^1_\alpha
& =  
\frac{1}{2\pi\im} (d \tau - \d f)|_{U^{(1)}_\alpha}, \\
\d \til{\omega}^1_\alpha
+ \frac{1}{2\pi\im} d\log h_\alpha
& =  
- \frac{1}{2\pi\im} \d \tau |_{U^{(2)}_\alpha}.
\end{align*}
By using these formulae, the coboundary of (\ref{cochain:thm_expression_beta}) yields the following cocycle in the complex (\ref{double_complex:F1FN}): 
$$
\left(
\frac{1}{2\pi\im} (d \tau - \d f), 
- \frac{1}{2\pi\im} \d \tau
\right) \in F^1\!A^2(G \times M) \oplus F^1\!A^1(G^2 \times M).
$$

Finally, we compute the corresponding cocycle $(E, \zeta) \in \mathcal{Z}$ along Lemma \ref{lem:iso_H3F1F2}. Then we find that the cocycle $(E, \zeta)$ is expressed as (\ref{formula:dfn_E_by_lambda}) and (\ref{formula:dfn_zeta_by_lambda}) by using $\lambda = \varpi(\tau)$, where $\varpi$ is the map given in Lemma \ref{lem:surjection_tau_to_lambda}.
\end{proof}


\subsection{Proof of Theorem \ref{ithm:obstruction}}

By the help of Proposition \ref{prop:quotient_Deligne}, the following exact sequence is induced from (\ref{exact_seq:simplicial_complex}):
\begin{equation}
\xymatrix@C=17pt@R=2pt{
&
H^1(G^{\bullet} \! \times \! M, \bF(2)) \ar[r]^-{\beta} &
H^2(G^{\bullet} \! \times \! M, F^1\!\F(2)) \ar[r] &
H^2(M/G, \F(2)) \\
\ar[r]^-{\varphi \circ q^*} &
H^2(G^{\bullet} \! \times \! M, \bF(2)) \ar[r]^-{\beta} &
H^3(G^{\bullet} \! \times \! M, F^1\!\F(2)). &
} \label{exact_seq:N=2}
\end{equation}

\begin{thm} \label{thm:cohomological_reduction_differentiable}
Let $G$ and $M$ be as in Proposition \ref{prop:quotient_Deligne}, and $(\G, \nabla, f)$ a strongly $G$-equivariant bundle gerbe with connection and curving over $M$.

(a) There exists a bundle gerbe with connection and curving $(\bar{\G}, \bar{\nabla}, \bar{f})$ over $M/G$ whose pull-back under the projection $q: M \to M/G$ is stably isomorphic to $(\G, \nabla, f)$, if and only if $\beta_G(\G, \nabla, f) = 0$.

(b) The stable isomorphism classes of such $(\bar{\G}, \bar{\nabla}, \bar{f})$ as above are in one to one correspondence with $Coker\{ \beta : H^1(G^{\bullet} \times M, \bF(2)) \to H^2(G^{\bullet} \times M, F^1\!\F(2)) \}$. In particular, if the Lie algebra $\g$ of $G$ is such that $[\g, \g] = \g$, then the stable isomorphism class of $(\bar{\G}, \bar{\nabla}, \bar{f})$ is unique.
\end{thm}

This is the strongly equivariant bundle gerbe version of Theorem \ref{ithm:coh_reduction_gerbe}, and includes Theorem \ref{ithm:obstruction} as (a).

\begin{proof}
First of all, we notice that the isomorphism $q^*$ in Proposition \ref{prop:quotient_Deligne} is induced by the simplicial map $q : G^\bullet \times M \to \{ e \}^\bullet \times (M/G)$. Recall that, for a bundle gerbe with connection and curving $(\bar{\G}, \bar{\nabla}, \bar{f})$ over $M/G$, the constructions in Lemma \ref{lem:sEBG_to_EBG} and the pull-back operations are compatible. Hence we have $\varphi \circ q^*(\delta(\bar{\G}, \bar{\nabla}, \bar{f})) = \delta_G(q^*(\bar{\G}, \bar{\nabla}, \bar{f}))$ in $H^2(G^{\bullet} \! \times \! M, \bF(2))$.

Using Proposition \ref{prop:expression_beta}, we identify $\beta_G(\G, \nabla, f)$ with $\beta(\delta_G(\G, \nabla, f))$. By the exactness of (\ref{exact_seq:N=2}), there is a class $\bar{c} \in H^2(M/G, \F(2))$ such that $\varphi \circ q^*(\bar{c}) = \delta_G(\G, \nabla, f)$ if and only if $\beta(\delta_G(\G, \nabla, f)) = 0$. By Proposition \ref{prop:classification_bgcc}, the class $\bar{c}$ is realized by a bundle gerbe with connection and curving $(\bar{\G}, \bar{\nabla}, \bar{f})$ over $M/G$. Because $\delta_G(q^*(\bar{\G}, \bar{\nabla}, \bar{f})) = \delta_G(\G, \nabla, f)$, the pull-back of $(\bar{\G}, \bar{\nabla}, \bar{f})$ is stably isomorphic to $(\G, \nabla, f)$ by Corollary \ref{cor:properties_char_class_sEBGcc} (c), which establishes (a).

For (b), it is clear that the set of the stable isomorphism classes of such $(\bar{\G}, \bar{\nabla}, \bar{f})$ as in (a) is identified with the kernel of $\varphi \circ q^*$, which is also identified with the cokernel of $\beta$. The last part of (b) follows from Corollary \ref{cor:vanishing_H2F1F2}. 
\end{proof}

We notice that the cohomology $H^1(G^{\bullet} \times M, \bF(2))$ is isomorphic to the isomorphism classes of $G$-equivariant flat $\T$-bundles \cite{Go1}. The homomorphism $\beta : H^1(G^{\bullet} \times M, \bF(2)) \to H^2(G^{\bullet} \times M, F^1\!\F(2))$ is interpreted as the assignment of $\frac{-1}{2\pi\im}$ times the moment. 

\medskip

Before we proceed to the next section, we consider the case where connections and curvings are absent.

\begin{prop} \label{prop:cohomological_reduction_topological}
Let $G$ and $M$ be as in Proposition \ref{prop:quotient_Deligne}, and $\G$ a strongly $G$-equivariant bundle gerbe over $M$.

(a) There exists a bundle gerbe $\bar{\G}$ over $M/G$ whose pull-back under the projection $q: M \to M/G$ is stably isomorphic to $\G$.

(b) The stable isomorphism class of such $\bar{\G}$ as above is unique.
\end{prop}

\begin{proof}
In the same way as $\delta_G(\G, \nabla, f)$, we have a class $\delta_G(\G) \in H^2(G^\bullet \times M, \u{\T})$ associated to $\G$. This is the complete invariant of stable isomorphism classes of $G$-equivariant bundle gerbes over $M$. It is shown \cite{Go2} that the projection $q : M \to M/G$ induces an isomorphism $q^* : H^2(M/G, \u{\T}) \to H^2(G^{\bullet} \times M, \u{\T})$ under the same assumption as in Proposition \ref{prop:quotient_Deligne}. Thus, by Proposition \ref{prop:classification_bg}, there always exists such a bundle gerbe $\bar{\G}$ as in (a), and the stable isomorphism class of such $\bar{\G}$ is unique.
\end{proof}


\section{Reduction}
\label{sec:reduction}

As is seen, Theorem \ref{thm:cohomological_reduction_differentiable} is proved in a rather abstract way by using the result on equivariant smooth Deligne cohomology groups. To understand it more geometrically, we describe in this section a reduction of a strongly equivariant bundle gerbe with connection. Theorem \ref{ithm:reduction_sebgcc} and Theorem \ref{ithm:difference} stated in Section \ref{sec:introduction} will be proved as Theorem \ref{thm:red_sEBGcc} and Theorem \ref{thm:different_lambda}, respectively. After that, a reduction of equivariant pseudo $\T$-bundle is also described. 

The reduction presented here is based on that of strongly equivariant bundle gerbes given by Mathai and Stevenson \cite{Ma-S}. When a finite group acts, the other method of reduction is known by Gawedzki and Reis \cite{Gaw-R}. Their reduction is applicable to equivariant bundle gerbes which are not necessarily strongly equivariant.


\subsection{Reduction of strongly equivariant bundle gerbes}
\label{subsec:reduction_sebg}

Let $G$ be a Lie group acting on a smooth manifold $M$, and $(\G, \nabla, f)$ a strongly $G$-equivariant bundle gerbe with connection and curving over $M$. Throughout this section, we put the following assumptions:

\begin{list}{}{\parsep=-2pt\topsep=4pt}
\item[(A1)] 
the action of $G$ on $M$ is free and locally trivial;

\item[(A2)] 
the quotient space $M/G$ is a smooth manifold in such a way that the natural projection $q : M \to M/G$ is smooth;

\item[(A3)] 
the action of $G$ on $Y$ satisfies the same assumptions as (A1) and (A2).
\end{list}

The assumptions (A1) and (A2) are those put in Proposition \ref{prop:quotient_Deligne}.

\medskip

First of all, we describe a reduction of the strongly $G$-equivariant bundle gerbe $\G = (Y, P, s)$, forgetting $\nabla$ and $f$. By (A3), $\bar{Y} = Y/G$ is a smooth manifold. The surjection $\pi : Y \to M$ induces $\bar{\pi} : \bar{Y} \to M/G$, which is also a surjective submersion admitting local sections. For each positive integer $p$, we can naturally identify the fiber product $\bar{Y}^{[p]}$ with the quotient space $(Y^{[p]})/G$. Hence we have a principal $\T$-bundle $\bar{P} = P/G$ over $\bar{Y}^{[2]}$. We can also identify $\delta \bar{P}$ with $(\delta P)/G$. Thus, the invariant section $s : Y^{[3]} \to \delta P$ uniquely corresponds to a section $\bar{s} : \bar{Y}^{[3]} \to \delta \bar{P}$. Because $\delta s = 1$, we have $\delta \bar{s} = 1$. Hence $\bar{\G} = (\bar{Y}, \bar{P}, \bar{s})$ is a bundle gerbe over $M/G$, which is a reduction of $\G$:

\begin{lem} 
The pull-back of the bundle gerbe $\bar{\G} = (\bar{Y}, \bar{P}, \bar{s})$ under the projection $q : M \to M/G$ is stably isomorphic to $\G$.
\end{lem}

We omit the proof of this lemma, since it follows from Theorem \ref{thm:red_sEBGcc}.

Notice that the stable isomorphism class of $\bar{\G}$ is uniquely determined by that of $\G$, as a result of Proposition \ref{prop:cohomological_reduction_topological}.

\medskip

Next we construct a connection and a curving on $\bar{\G}$ from $\nabla$ and $f$. By (A1) and (A2), we have a principal $G$-bundle $q : M \to M/G$. We take and fix a connection $\Xi$ on this $G$-bundle.

Now we suppose the vanishing of the obstruction: $\beta_G(\G, \nabla, f) = 0$. This is equivalent to the existence of a map $\lambda : Y \to \g^{\dagger}$ such that $\delta \lambda = \til{\lambda}$ and $(E, \zeta) = 0$, where $(E, \zeta) \in \mathcal{Z}$ is defined by (\ref{formula:dfn_E_by_lambda}) and (\ref{formula:dfn_zeta_by_lambda}). Fix such a map for a while. We define a 1-form $\kappa \in \im A^1(Y)$ by $\kappa = \langle \pi^*\Xi | \lambda \rangle$. Since $\zeta = 0$, the 1-form $\kappa$ is $G$-invariant, so that we have a strongly $G$-equivariant bundle gerbe with connection and curving $(\G, \nabla - \delta \kappa, f - d\kappa)$ over $M$. We can easily verify
\begin{eqnarray*}
\iota_{X^*} (\nabla - \delta \kappa) 
& = & 
\langle X | \til{\lambda} - \delta \lambda \rangle = 0, \\
\iota_{X^*} (f - d \kappa)
& = &
2\pi\im \langle X | E \rangle = 0.
\end{eqnarray*}
Hence $\nabla - \delta \kappa$ descends to give a connection $\bar{\nabla}$ on $\bar{\G}$, and $f - d \kappa$ descends to give a curving $\bar{f}$ for $\bar{\nabla}$. This construction of $(\bar{\G}, \bar{\nabla}, \bar{f})$ is a reduction of $(\G, \nabla, f)$:

\begin{thm} \label{thm:red_sEBGcc}
If we fix a map $\lambda : Y \to \g^{\dagger}$ such that $\delta \lambda = \til{\lambda}$ and $(E, \zeta) = 0$, then $q^*(\bar{\G}, \bar{\nabla}, \bar{f})$ is stably isomorphic to $(\G, \nabla, f)$ as a strongly $G$-equivariant bundle gerbe with connection and curving.
\end{thm}

\begin{proof}
We write the pull-back $q^*\bar{\G}$ as $(q^*\bar{Y}, (\til{q}^{[2]})^*\bar{P}, (\til{q}^{[3]})^*\bar{s})$. We have the natural $G$-equivariant maps $\varphi : Y \to q^*\bar{Y}$ and $\til{\varphi} : P \to (\til{q}^{[2]})^*\bar{P}$. These maps form an isomorphism of strongly equivariant bundle gerbes with connection and curving $(\varphi, \til{\varphi}) : (\G, \nabla - \delta \kappa, f - d \kappa) \to q^*(\bar{\G}, \bar{\nabla}, \bar{f})$. Thus Lemma \ref{lem:sEBGcc_iso} implies that $q^*(\bar{G}, \bar{\nabla}, \bar{f})$ and $(\G, \nabla - \delta \kappa, f - d\kappa)$ are stably isomorphic. Because $\kappa$ is $G$-invariant, Lemma \ref{lem:shift_connection_and_curving} implies that $(\G, \nabla - \delta \kappa, f - d \kappa)$ and $(\G, \nabla, f)$ are stably isomorphic.
\end{proof}

Our reduction of $(\G, \nabla, f)$ is performed under a choice of $\lambda$. In general, the choice is not unique. We can readily see the next lemma:

\begin{lem} \label{lem:expression_beta}
If $\beta_G(\G, \nabla, f) = 0$, then the set of maps $\lambda : Y \to \g^{\dagger}$ such that $\delta \lambda = \til{\lambda}$ and $(E, \zeta) = 0$ is an affine space under the vector space
$$
H^2(G^{\bullet} \times M, F^1\!\F(2)) \cong
\left\{ 
\mu \in A^0(M, \g^{\dagger}) |\ g^*\mu = \Ad_g\mu, \ d\mu = 0 
\right\}.
$$
Thus, If $[\g, \g] = \g$, then the choice of such $\lambda$ is unique.
\end{lem}

Let $\lambda' : Y \to \g^\dagger$ be the other choice, and $(\bar{\G}, \bar{\nabla}', \bar{f}')$ the bundle gerbe with connection and curving over $M/G$ obtained as the reduction with respect to $\lambda'$. By the unique map $\mu : M \to \g^{\dagger}$ such that $\pi^*\mu = \lambda' - \lambda$, we obtain $\kappa' - \kappa = \pi^*\langle \Xi | \mu \rangle$, where $\kappa' = \langle \pi^*\Xi | \lambda' \rangle$. So we have 
$$
(\G, \nabla - \delta \kappa', f - d \kappa') = 
(\G, \nabla - \delta \kappa, f - d \kappa - d \pi^*\langle \Xi | \mu \rangle).
$$
Thus the choice of $\lambda$ does not affect connections on $\bar{\G}$: we always have $\bar{\nabla}' = \bar{\nabla}$. On the other hand, the choice of $\lambda$ affects curvings on $\bar{\G}$. Though $\langle \Xi | \mu \rangle$ does not descend to $M/G$ in general, $d \langle \Xi | \mu \rangle$ descends to give the 2-form $\bar{\sigma} \in \im A^2(M/G)$ such that $q^*\bar{\sigma} = d \langle \Xi | \mu \rangle$. So we can write $\bar{f}' = \bar{f} - \bar{\pi}^*\bar{\sigma}$.

\begin{thm} \label{thm:different_lambda}
Suppose that $(\bar{\G}, \bar{\nabla}, \bar{f})$ and  $(\bar{\G}, \bar{\nabla}, \bar{f}')$ are obtained as the reductions of $(\G, \nabla, f)$ with respect to $\lambda$ and $\lambda'$, respectively. They are stably isomorphic if and only if the difference $\lambda' - \lambda$ is the moment of a $G$-equivariant flat principal $\T$-bundle over $M$. 
\end{thm}

\begin{proof}
Suppose that $(\bar{\G}, \bar{\nabla}, \bar{f})$ and $(\bar{\G}, \bar{\nabla}, \bar{f} - \bar{\pi}^*\bar{\sigma})$ are stably isomorphic. By Corollary \ref{cor:different_curving}, there is a principal $\T$-bundle $(\bar{R}, \bar{\eta})$ over $M/G$ such that $F(\bar{\eta}) = \bar{f'} - \bar{f} = - \bar{\sigma}$. The pull-back $q^*\bar{R}$ is naturally a $G$-equivariant principal $\T$-bundle over $M$. Here we consider a $G$-invariant connection $q^*\bar{\eta} + \langle \Xi | \mu \rangle$ on $q^*\bar{R}$. A direct calculation shows that the connection is flat and its moment is $\mu$.

Conversely, when $\mu$ is the moment of a $G$-equivariant flat principal $\T$-bundle $(R, \eta)$ over $M$, we consider the $G$-invariant connection $\eta - \langle \Xi | \mu \rangle$ on $R$. Because the moment of $\eta - \langle \Xi | \mu \rangle$ vanishes, $(R, \eta - \langle \Xi | \mu \rangle)$ reduces to give a principal $\T$-bundle with connection $(\bar{R}, \bar{\eta})$ over $M/G$. Clearly, the curvature is $F(\bar{\eta}) = - \bar{\sigma}$. Hence $(\bar{\G}, \bar{\nabla}, \bar{f})$ and $(\bar{\G}, \bar{\nabla}, \bar{f} - \bar{\pi}^*\bar{\sigma})$ are stably isomorphic by Corollary \ref{cor:different_curving}.
\end{proof}

Recalling Theorem \ref{thm:cohomological_reduction_differentiable} (b), we find that our reduction yields, by varying the choice of $\lambda$, all the stable isomorphism classes of bundle gerbes with connection and curving over $M$ whose pull-backs are stably isomorphic to $(\G, \nabla, f)$.

\smallskip

So far, we fixed a choice of a connection $\Xi$ on the $G$-bundle $q : M \to M/G$. We here vary the choice of $\Xi$, fixing a choice of $\lambda$.

\begin{prop}
As far as $\lambda$ is fixed, the stable isomorphism class of the bundle gerbe with connection and curving over $M/G$ obtained by the reduction is independent of the choice of a connection on the $G$-bundle $q : M \to M/G$.
\end{prop}

\begin{proof}
Let $\Xi$ and $\Xi'$ be connections on the $G$-bundle $q : M \to M/G$. We put $\kappa = \langle \pi^*\Xi | \lambda \rangle$ and $\kappa' = \langle \pi^*\Xi' | \lambda \rangle$. The reduction yields $(\bar{\G}, \bar{\nabla}, \bar{f})$ and $(\bar{\G}, \bar{\nabla}', \bar{f}')$ by taking the quotient of $(\G, \nabla - \delta \kappa, f - d \kappa)$ and $(\G, \nabla - \delta \kappa', f - d \kappa')$, respectively. Denote by $\xi \in A^0(M, \g)$ the difference $\xi = \Xi' - \Xi$. Using the $G$-invariant 1-form $\alpha \in \im A^1(M)$ defined by $\alpha = \langle \xi | \lambda \rangle$, we can write $\kappa' - \kappa = \pi^*\alpha$. Hence we have $\nabla - \delta \kappa' = \nabla - \delta \kappa$, so that $\bar{\nabla} = \bar{\nabla'}$. Notice that $\alpha$ vanishes in the direction of $G$. Thus we have $\bar{f}' = \bar{f} - \bar{\pi}^*d\bar{\alpha}$, where $\bar{\alpha} \in \im A^1(M/G)$ is the 1-form such that $q^*\bar{\alpha} = \alpha$. Since $- d\bar{\alpha}$ is the curvature of the connection $- \bar{\alpha}$ on the trivial $\T$-bundle over $M/G$, Corollary \ref{cor:different_curving} completes the proof.
\end{proof}


\subsection{Reduction of pseudo $\T$-bundles with connection}

Let $\G = (Y, P, s)$ be a strongly $G$-equivariant bundle gerbe over $M$ equipped with a $G$-invariant connection $\nabla$ and a $G$-invariant curving $f$. 

\begin{lem} \label{lem:data_pseudo_T_bundle_conn}
Let $((R, v), \eta)$ be a $G$-equivariant pseudo $\T$-bundle with connection for $(\G, \nabla, f)$, and $\rho : Y \to g^\dagger$ the moment associated with $(R, \eta)$. We fix a map $\lambda : Y \to \g^\dagger$ such that $\delta \lambda = \til{\lambda}$, and define a $G$-invariant 2-form $\omega \in \im A^2(M)$ and a map $\mu : M \to \g^\dagger$ by
\begin{align}
\pi^*\omega & =  F(\eta) - f, \label{formula:omega} \\
\pi^*\mu & =  \rho - \lambda. \label{formula:mu}
\end{align}
For $X \in \g$ and $g \in G$ we have
\begin{align}
d \omega & =  - \Omega,  \\
\langle X | d \mu \rangle + \iota_{X^*}\omega & = 
- \langle X | 2\pi\im E \rangle, \\
\langle X | g^*\mu - \Ad_g \mu \rangle  & = 
- \langle X | 2\pi\im \zeta(g) \rangle,
\end{align}
where $\Omega$ is the 3-curvature of $f$, and $(E, \zeta) \in \mathcal{Z}$ is defined by using $\lambda$.
\end{lem}

\begin{proof}
One can easily verify this lemma by straightforward computations.
\end{proof}

Recall the reduction of $(\G, \nabla, f)$. We take and fix a connection $\Xi$ on the $G$-bundle $q : M \to M/G$. Choosing a map $\lambda$ such that $\delta \lambda = \til{\lambda}$ and $(E, \zeta) = 0$, we put $\kappa = \langle \pi^*\Xi | \lambda \rangle$. Then the quotient of $(\G, \nabla - \delta \kappa, f - d \kappa)$ gives $(\bar{\G}, \bar{\nabla}, \bar{f})$, which is the reduction of $(\G, \nabla, f)$ with respect to $\lambda$. 

We can see that $((R, v), \eta - \kappa)$ is a $G$-equivariant pseudo $\T$-bundle with connection for $(\G, \nabla - \delta \kappa, f - d \kappa)$. So one may expect that $((R, v), \eta - \kappa)$ induces a pseudo $\T$-bundle with connection for $(\bar{\G}, \bar{\nabla}, \bar{f})$. However, it does not hold in general.

\begin{prop} \label{prop:reduction_pseudo_T_bundle}
If $\mu$ vanishes, then we can construct a pseudo $\T$-bundle with connection $((\bar{R}, \bar{v}), \bar{\eta})$ for $(\bar{\G}, \bar{\nabla}, \bar{f})$ from $((R, v), \eta)$.
\end{prop}

\begin{proof}
By taking the quotient of $(R, v)$, we obtain a pseudo $\T$-bundle $(\bar{R}, \bar{v})$ for $\bar{\G}$. It is straightforward to see that the moment of $\eta - \kappa$ is $\mu$. Thus, under the assumption of this proposition, the $G$-invariant connection $\eta - \kappa$ descends to give a connection $\bar{\eta}$ on the $\T$-bundle $\bar{R} = R/G$. We can easily verify that $\bar{\eta}$ gives rise to a connection on the pseudo $\T$-bundle $(\bar{R}, \bar{v})$.
\end{proof}

If $\mu = 0$ (and $(E, \zeta) = 0$), then the 2-form $\omega = (f - d \kappa) - F(\eta - d \kappa)$ satisfies $\iota_{X^*}\omega = 0$ for all $X \in \g$. Hence the 2-form induces $\bar{\omega} \in \im A^2(M/G)$ such that $q^*\bar{\omega} =  \omega$. If $\bar{\Omega}$ is the 3-curvature of $\bar{f}$, then we have $d\bar{\omega} = - \bar{\Omega}$.


\section{Example}
\label{sec:example}

In this section, we give some examples of reductions of strongly equivariant bundle gerbes. The first example deals with an action of $S^1$ on $S^3$. In this case, we can obtain distinct stable isomorphism classes of bundle gerbes with connection and curving over the quotient space by the reductions. In the second example, we will define an equivariant bundle gerbe over an infinite dimensional space as a \textit{lifting bundle gerbe} \cite{Mu}, and will obtain the ``basic bundle gerbe'' on $SU(2)$ by the reduction. The third example is a reduction of a pseudo $\T$-bundle with connection, which is tied up with the second example.


\subsection{The trivial bundle gerbe over $S^3$}

Let us consider the diagonal embedding of $S^1 = \T$ into $S^3 = SU(2)$. This gives a free action of $G = S^1$ on $M = S^3$. The quotient space is $M/G \cong S^2$, and $q : M \to M/G$ is the Hopf fibration.

\medskip

We define a trivial bundle gerbe $\G = (Y, P, s)$ over $M$ as follows. We put $Y = M$ and define $\pi : Y \to M$ to be the identity. For all $p$ we can identify $Y^{[p]}$ with $M$ by the diagonal embedding of $M$ into $Y^p$. We put $P = Y^{[2]} \times \T$ and $s = 1$. By means of the trivial action of $G$ on $P$, the bundle gerbe $\G$ is strongly $G$-equivariant. 

Let $\nabla$ be the trivial connection on the trivial $\T$-bundle $P$. It is obvious that $\nabla$ is a $G$-invariant connection on $\G$. If we define $f \in \im A^2(Y)$ by $f = 0$, then $f$ is a $G$-invariant curving for $\nabla$. 

\medskip

As is clear, the strongly $G$-equivariant bundle gerbe with connection and curving $(\G, \nabla, f)$ gives $\delta_G(\G, \nabla, f) = 0$, so that $\beta_G(\G, \nabla, f) = 0$. To perform the reduction presented in Section \ref{subsec:reduction_sebg}, we choose a map $\lambda : Y \to \g^{\dagger}$ such that $\delta \lambda = \til{\lambda}$ and $(E, \zeta) = 0$, where $\til{\lambda} : Y^{[2]} \to \g^{\dagger}$ is the moment associated to $(P, \nabla)$, and $(E, \zeta) \in \mathcal{Z}$ is given by (\ref{formula:dfn_E_by_lambda}) and (\ref{formula:dfn_zeta_by_lambda}). It is clear that $\tilde{\lambda} = 0$. If we take $\lambda_0 : Y \to \g^{\dagger}$ to be $\lambda_0 = 0$, then we have $\delta \lambda_0 = \til{\lambda}$ and $(E, \zeta) = 0$. So we can take $\lambda_0$ as a choice. 

We also have the other choice in the present case. Note that $G = S^1$ is abelian and that the adjoint action on $\g = \im\R$ is trivial. For $r \in \R$ we define $\lambda_r : \g \to \im\R$ by $\lambda_r(z) = rz$. If we regard $\lambda_r \in \g^{\dagger}$ as an element $\lambda_r \in A^0(Y, \g^{\dagger})$ in the natural fashion, then $\lambda_r$ is the other choice of $\lambda$ such that $\delta \lambda = \til{\lambda}$ and $(E, \zeta) = 0$. Because $M = S^3$ is connected, we have $H^2(G^{\bullet} \times M, F^1\bF(2)) \cong \g^{\dagger} \cong \R$. Hence, by Lemma \ref{lem:expression_beta}, any such choice of $\lambda$ is of the form $\lambda_r$.

\medskip

Now we perform the reduction of $(\G, \nabla, f)$ with respect to $\lambda_r$. Let $\Xi$ be a connection on the $S^1$-bundle $S^3 \to S^2$. We define $\kappa_r \in \im A^1(Y)$ by $\kappa_r = \langle \Xi | \lambda_r \rangle = r \Xi$. We take the quotient of $(\G, \nabla - \delta \kappa_r, f - d \kappa_r)$ to obtain $(\bar{\G}, \bar{\nabla}, \bar{f}_r)$ over $M/G \cong S^2$. By the construction of $\bar{\G} = (\bar{Y}, \bar{P}, \bar{s})$, we have $\bar{Y} = M/G$. Hence the curving $\bar{f}_r$ is a 2-form on $M/G \cong S^2$. If $F(\Xi) \in \im A^2(S^2)$ denotes the curvature of the connection $\Xi$, then we have the expression $f_r = r F(\Xi)$.

\begin{prop}
The bundle gerbe with connection and curving $(\bar{\G}, \bar{\nabla}, \bar{f}_r)$ over $M/G \cong S^2$ is trivial if and only if $r$ is an integer.
\end{prop}

Though we can prove this proposition appealing to Theorem \ref{thm:different_lambda}, we give an easier proof by using Corollary \ref{cor:different_curving}.

\begin{proof}
The Euler class of the $\T$-bundle $S^3 \to S^2$ is a generator of $H^2(S^2, \Z) \cong \Z$. Since $F(\Xi)$ is the curvature of a connection on this bundle, $f_r = r F(\Xi)$ is the curvature of a connection on a $\T$-bundle over $S^2$ if and only if $r \in \Z$. 
\end{proof}

By Proposition \ref{prop:Deligne_coh:manifold} (b) and Proposition \ref{prop:classification_bgcc}, the the group of stable isomorphism classes of bundle gerbes with connection and curving over $S^2$ is isomorphic to $H^2(S^2, \T) \cong H^2(S^2, \R) / H^2(S^2, \Z)$. Since $\frac{-r}{2\pi\im}F(\Theta)$ represents the de Rham cohomology class corresponding to $r \in \R \cong H^2(S^2, \R)$, all the stable isomorphism classes of bundle gerbes with connection and curving over $S^2$ are obtained by the reduction of $(\G, \nabla, f)$.


\subsection{Chern-Simons bundle gerbes}
\label{subsec:CS_bundle_gerbe}

As is pointed out by Freed \cite{F2}, higher gerbes appear as basic ingredients when we consider higher codimensions in field theories. In the context of $SU(2)$ Chern-Simons theory \cite{F1,Wi}, such a gerbe is formulated over the space of connections on an $SU(2)$-bundle over an oriented closed 1-manifold \cite{Go3}. The specialization in the case of the trivial $SU(2)$-bundle over $S^1$ provides the example below. We note that it is related to the construction of the basic gerbes over compact simple Lie groups due to Behrend, Xu and Zhang \cite{B-X-Z}.

\medskip

Let $\su(2)$ be the Lie algebra of $SU(2)$. We denote the space of connections on the trivial $SU(2)$-bundle over $S^1$ by $M = A^1(S^1, \su(2))$. The group $G = LSU(2) = C^\infty(S^1, SU(2))$ of free loops in $SU(2)$ acts on the space of connections by the gauge transformation. According to the convention of the present paper, we consider the left action of $g \in LSU(2)$ on $A \in A^1(S^1, \su(2))$ given by $A \mapsto gAg^{-1} + g dg^{-1}$.

Let us choose the standard base point on $S^1$, and denote the based loop group by $G_0 = \Omega SU(2)$. As is well-known, the action of $G_0 \subset G$ on $M$ is free. The quotient space $M/G_0$ is identified with $SU(2)$, and the projection map $q : M \to M/G_0$ with the holonomy of a connection around $S^1$. Note that $G/G_0 \cong SU(2)$ acts on $M/G_0 \cong SU(2)$ by the adjoint action.

\medskip

Let $\pi : Y \to M$ be the trivial $LSU(2)$-bundle $Y = M \times LSU(2)$. We write $\Gamma = LSU(2)$ for the structure group of $Y$ in order to distinguish it from the group $G = LSU(2)$ acting on the base space $M$. It is known \cite{P-S} that the loop group $\Gamma = LSU(2)$ has a central extension $\hGamma^k$ for each integer $k \in \Z$:
$$
1 \longrightarrow
\T \longrightarrow
\hGamma^k \stackrel{\varpi}{\longrightarrow}
\Gamma \longrightarrow
1.
$$

\begin{dfn}
We define a bundle gerbe $\G^k$ over $M$ by the \textit{lifting bundle gerbe} \cite{Mu} associated with the central extension $\hGamma^k$ and the trivial $\Gamma$-bundle $\pi : Y \to M$. We call $\G^k$ the \textit{Chern-Simons bundle gerbe}.
\end{dfn}

The concrete description of $\G^k = (Y, P, s)$ is as follows. We introduce a map $\tau : Y^{[2]} \to \Gamma$ by $y_ 1 \tau(y_1, y_2) = y_2$ for $(y_1, y_2) \in Y^{[2]}$. When we write $y_i = (A, \gamma_i)$, we have $\tau(y_1, y_2) = \gamma_1^{-1}\gamma_2$. The principal $\T$-bundle $P \to Y^{[2]}$ is defined by $P = \tau^*\hGamma^k$. The section $s : Y^{[3]} \to \delta P$ is defined by $s(y_1, y_2, y_3) = \h{\gamma}_{23} \otimes (\h{\gamma}_{12} \h{\gamma}_{23})^{\otimes -1} \otimes \h{\gamma}_{12}$, where $\h{\gamma}_{ij} \in \hGamma^k$ is an element such that $\varpi(\h{\gamma}_{ij}) = \tau(y_i, y_j)$.

\begin{prop}
The Chern-Simons bundle gerbe $\G^k = (Y, P, s)$ can be made into a strongly $G$-equivariant bundle gerbe over $M$.
\end{prop}

\begin{proof}
In general, the lifting bundle gerbe associated with an equivariant principal bundle can be made into strongly equivariant. We define a left action of $G = LSU(2)$ on $Y = M \times \Gamma = A^1(S^1, \su(2)) \times LSU(2)$ by $(A, \gamma) \mapsto (g A g^{-1} + g dg^{-1}, g \gamma)$. This is a lift of the action of $G$ on $M$ to that on $Y$. In order to make the $\T$-bundle $P \to Y^{[2]}$ into $G$-equivariant, we note that the left action of $G$ on the $\Gamma$-bundle $Y$ commutes with the right action of $\Gamma$. We define an action of $G$ on $Y^{[2]} \times \hGamma^k$ by the induced $G$-action on $Y^{[2]}$ and by the trivial action on $\hGamma^k$. Because $P = \tau^*\hGamma^k$ is a $G$-invariant subspace of $Y^{[2]} \times \hGamma^k$, we obtain a $G$-action on $P$ which makes $P \to Y^{[2]}$ into a $G$-equivariant $\T$-bundle. It is easy to see that the section $s : Y^{[2]} \to \delta P$ is $G$-invariant.
\end{proof}

Since $G_0$ is a subgroup of $G$, the Chern-Simons bundle gerbe $\G^k$ is a strongly $G_0$-equivariant bundle gerbe as well. Since the action of $G_0$ on $Y$ is free, we can perform the reduction of $\G^k$. Then we obtain a bundle gerbe $\bar{\G}^k$ over $M/G_0 \cong SU(2)$. We will see later that the Dixmier-Douady class of the bundle gerbe is $k \in \Z \cong H^3(SU(2), \Z)$. 

\bigskip

To construct a connection and a curving on $\G^k$, we use results in \cite{Go4,Mu-S2}. We denote by $\Lie\Gamma = L\su(2) = C^\infty(S^1, \su(2))$ the Lie algebra of $\Gamma = LSU(2)$. The Lie algebra $\Lie \hGamma^k$ of $\hGamma^k$ is naturally identified with the vector space $L\su(2) \oplus \im\R$ on which the Lie bracket is defined by
$$
[X_1 \oplus z_1, X_2 \oplus z_2] = [X_1, X_2] \oplus c(X_1, X_2),
$$
where the Lie algebra 2-cocycle $c$ is
$$
c(X_1, X_2) = - \frac{k\im}{2\pi} \int_{S^1} \Tr(X_1 dX_2).
$$
Because $\Gamma$ acts on $\Lie\hGamma^k$ through the adjoint action of $\hGamma^k$, it is possible to define a map $Z : \Gamma \to \Hom(\Lie\Gamma, \im\R)$ by $(Z(\gamma) | X ) = \Ad_\gamma(X \oplus 0) - (\Ad_\gamma X) \oplus 0$, where $( \ | \ ) : \Hom(\Lie\Gamma, \im\R) \otimes \Lie\Gamma\to \im\R$ is the natural contraction. An explicit formula for $Z$ is 
$$
( Z(\gamma) | X )= 
\frac{k\im}{2\pi}\int_{S^1}\Tr(\gamma^{-1}d\gamma \ X).
$$
This follows from the simplicity of $\su(2)$ and the following formula:
$$
c(\Ad_\gamma X_1, \Ad_\gamma X_2)
= c(X_1, X_2) + ( Z(\gamma) | [X_1, X_2] ).
$$
Let $\theta_\Gamma$ and $\h{\theta}_\Gamma$ be the left invariant Maurer-Cartan forms on $\Gamma$ and $\h{\Gamma}^k$, respectively. Then the 1-form $\nu = \h{\theta}_\Gamma - \varpi^*\theta_\Gamma \oplus 0$ is a connection on the principal $\T$-bundle $\varpi : \h{\Gamma}^k \to \Gamma$, and its curvature is $F(\nu) = - \frac{1}{2}c(\theta_\Gamma, \theta_\Gamma)$.

Since $\pi : Y \to M$ is trivial as a $\Gamma$-bundle, the Maurer-Cartan form $\theta_\Gamma$ defines the trivial connection $\Theta$ on $Y$. We define $L : Y \to \Hom(\Lie\hGamma^k, \im\R)$ by
$$
( L(A, \gamma) | X \oplus z ) =
\frac{k\im}{2\pi} \int_{S^1} 
\Tr((\gamma^{-1} A \gamma + \gamma^{-1}d\gamma) X)
+ z.
$$
This gives a section $L \in A^0(M, Y \times_{Ad} \Hom(\Lie\hGamma^k, \im\R))$ called a \textit{splitting} \cite{Bry1}. A bundle gerbe connection $\nabla \in \im A^1(P)$ and a bundle gerbe curving $f \in \im A^2(Y)$ are given by
\begin{align*}
\nabla &= \tau^*\nu + ( Z(\tau^{-1}) | \pi_2^*\Theta ), \\
f &= - ( L | F(\Theta \oplus 0) ),
\end{align*}
where $F(\Theta \oplus 0) = d (\Theta \oplus 0) +\frac{1}{2}[\Theta \oplus 0, \Theta \oplus 0 ] = F(\Theta) \oplus \frac{1}{2}c(\Theta, \Theta)$. In the present case, we have $F(\Theta) = 0$, so that $f = - \frac{1}{2}c(\Theta, \Theta)$ and $(\G, \nabla, f)$ is flat.

\begin{prop}
The connection $\nabla$ and the curving $f$ are $G$-invariant.
\end{prop}

\begin{proof}
Note that the connection $\Theta$ on $Y$ and its curvature $F(\Theta)$ are $G$-invariant. Because $\tau : Y^{[2]} \to \Gamma$ and $\pi_2 : Y^{[2]} \to Y$ are $G$-equivariant, the connection $\nabla$ is $G$-invariant. By a direct computation, we can prove that $L$ is $G$-invariant: $g^*L = L$ for $g \in G$. Hence the curving is also $G$-invariant.
\end{proof}

The moment $\til{\lambda} \in A^0(Y^{[2]}, \g^{\dagger})$ associated with the $G$-equivariant $\T$-bundle $(P, \nabla)$ is expressed as
\begin{equation}
\begin{split}
\langle X | \til{\lambda}((A, \gamma_1), (A, \gamma_2)) \rangle 
& = 
\left( Z(  \tau( (A, \gamma_1), (A, \gamma_2) )^{-1}) | 
\Theta((A, \gamma_1); X^*) \right) \\
& = 
\frac{k\im}{2\pi} \int_{S^1}
\Tr ( \gamma_2 d\gamma_2^{-1} X - \gamma_1 d\gamma_1^{-1} X),
\end{split} \label{formula:CS_tilde_lambda}
\end{equation}
where $X^*$ is the tangent vector on $Y$ generated by the infinitesimal action of $X \in \g = L\su(2)$. One choice of $\lambda \in A^0(Y, \g^{\dagger})$ such that $\delta \lambda = \til{\lambda}$ is
\begin{equation}
\begin{split}
\langle X | \lambda (A, \gamma) \rangle 
& = 
- (L(A, \gamma) | \Theta((A, \gamma); X^*) \oplus 0) \\
& = 
- \frac{k\im}{2\pi} \int_{S^1}
\Tr \left( (A - \gamma d\gamma^{-1}) X \right).
\end{split} \label{formula:CS_lambda}
\end{equation}
Using this map $\lambda$, we compute $(E, \zeta) \in A^1(M, \g^*) \oplus C^\infty(G, A^0(M, \g^*))$ according to (\ref{formula:dfn_E_by_lambda}) and (\ref{formula:dfn_zeta_by_lambda}). Then we obtain the following formulae:
\begin{align}
\langle X | E(A; \alpha) \rangle 
& = 
- \frac{k}{4\pi^2} \int_{S^1} \Tr (\alpha X), 
\label{formula:CS_E} \\
\langle X | \zeta(g) \rangle 
& = 
0.
\label{formula:CS_zeta}
\end{align}
We write $\g_0 = \Omega \su(2)$ for the Lie algebra of $G_0 = \Omega SU(2)$. The moment $\til{\lambda}_0 \in A^0(Y^{[2]}, \g_0^{\dagger})$ associated with the $G_0$-equivariant $\T$-bundle $(P, \nabla)$ is given by the same formula as (\ref{formula:CS_tilde_lambda}). Hence the map $\lambda_0 \in A^0(Y, \g_0^{\dagger})$ given by the same formula as (\ref{formula:CS_lambda}) satisfies $\delta \lambda_0 = \til{\lambda}_0$. If we use $\lambda_0$ to define $(E_0, \zeta_0) \in A^1(M, \g_0^*) \oplus C^\infty(G_0, A^0(M, \g_0^*))$ by (\ref{formula:dfn_E_by_lambda}) and (\ref{formula:dfn_zeta_by_lambda}), then it has the same expression as (\ref{formula:CS_E}) and (\ref{formula:CS_zeta}). 

\begin{lem}
If $k \neq 0$, then $\beta_G(\G^k, \nabla, f) \neq 0$ and $\beta_{G_0}(\G^k, \nabla, f) \neq 0$.
\end{lem}

\begin{proof}
For any $\mu \in A^0(M, \g^*)$, cocycles $(E, \zeta(g))$ and $(E + d\mu, \zeta(g) + g^*\mu - \Ad_g\mu)$ define the same class in $H^3(G^{\bullet} \times M, F^1\!\F(2)) \cong \mathcal{Z}/\mathcal{B}$. If we define $\mu$ by
$$
\langle X | \mu(A) \rangle =
\frac{k}{4\pi^2} \int_{S^1} \Tr(A X), 
$$
then we have $E + d\mu = 0$ and $\langle X | \zeta(g) + g^*\mu - \Ad_g\mu \rangle = \frac{1}{2\pi\im} (Z(g^{-1}) | X)$. 

Now we assume $[E, \zeta]$ to be trivial. Then there exists $\mu' \in A^0(M, \g^*)$ such that $d\mu' = 0$ and $\zeta(g) + g^*\mu - \Ad_g\mu = g^*\mu' - \Ad_g\mu'$. In other words, we can express $Z : \Gamma \to \Hom(\Lie\Gamma, \im\R)$ as $Z(g) = 2\pi\im (\mu' - \Ad_g\mu')$ using $\mu' \in \Hom(\Lie\Gamma, \im\R)$. Notice that the Lie algebra 2-cocycle $c$ of $\Lie\hGamma^k$ is obtained by $c(X, Y) = \frac{d}{dt}\big|_{t = 0}(Z(e^{-tX}) | Y)$. Thus, if such $\mu'$ as above exists, then $c$ is expressed as a coboundary, so that the central extension $\Lie\hGamma^k$ is trivial. This is a contradiction, and $\beta_G(\G^k, \nabla, f) = [E, \zeta]$ is a non-trivial class. Because $\hGamma^k$ induces a non-trivial central extension of $G_0 = \Omega SU(2)$ by restriction, the same argument shows that $\beta_{G_0}(\G^k, \nabla, f) = [E_0, \zeta_0]$ is also non-trivial.
\end{proof}

By this lemma, it turns out that we cannot perform the reduction of the strongly $G_0$-equivariant bundle gerbe with connection and curving $(\G^k, \nabla, f)$.

\medskip

In order to obtain a connection and curving on $\bar{\G}^k$ by a reduction, we replace the curving $f$ by the other one using a 2-form on $M = A^1(S^1, \su(2))$. As a preliminary, we denote by $\theta = g^{-1}dg$ the left invariant Maurer-Cartan form on $SU(2)$, and by $\bar{\theta} = dgg^{-1}$ the right invariant one. We define an adjoint invariant closed 3-form $\chi \in A^3(SU(2))$ and a 1-form $e \in A^1(SU(2), \su(2)^*)$ by
\begin{align*}
\chi &= - \frac{k}{24\pi^2} \Tr(\theta \wedge \theta \wedge \theta), \\
\langle X | e \rangle &= - \frac{k}{24\pi^2} \Tr((\bar{\theta} + \theta) X ).
\end{align*}
When we consider the adjoint action of $SU(2)$ on itself, we obtain $\langle X | de \rangle = \iota_{X^*}\chi$ and $g^*e = \Ad_g e$ for $X \in \su(2)$ and $g \in SU(2)$.

\begin{lem}[\cite{A-M-M}] \label{lem:A-M-M}
There is a $G$-invariant 2-form $\Upsilon \in \im A^2(M)$ such that
\begin{align*}
\frac{-1}{2\pi\im} d \Upsilon
& = q^*\chi, \\
\langle X | E \rangle + \frac{1}{2\pi\im} \iota_{X^*} \Upsilon
& = \langle X(0) | q^*e \rangle,
\end{align*}
where $q : M \to M/G_0$ is the projection, and we write $X(0) \in \su(2)$ for the evaluation of $X \in \g = L\su(2)$ by the base point on $S^1$.
\end{lem}

We refer the reader to \cite{A-M-M,G-S} for an explicit description of $\Upsilon$ and the proof of the lemma above.

Since $\Upsilon$ is $G$-invariant, we obtain a strongly $G$- and $G_0$-equivariant bundle gerbe with connection and curving $(\G^k, \nabla, f + \Upsilon)$ over $M$. We use the curving $f + \Upsilon$ and $\lambda$ given by (\ref{formula:CS_lambda}) in order to define $(E', \zeta') \in A^1(M, \g^*) \oplus C^\infty(G, A^0(M, \g^*))$. Similarly, $(E'_0, \zeta'_0) \in A^1(M, \g_0^*) \oplus C^\infty(G_0, A^0(M, \g_0^*))$ is defined. 

\begin{prop}
The cocycle $(E', \zeta')$ is expressed as
\begin{align}
\langle X | E' \rangle & =  \langle X(0) | q^*e \rangle, 
\label{formula:modified_E_CS} \\
\langle X | \zeta' \rangle & =  0.
\end{align}
Thus, we have $(E'_0, \zeta'_0) = 0$ and $\beta_{G_0}(\G^k, \nabla, f + \Upsilon) = 0$.
\end{prop}

\begin{proof}
By (\ref{formula:dfn_zeta_by_lambda}), we have $\langle X | E' \rangle = \langle X | E \rangle + \frac{1}{2\pi\im}\iota_{X^*}\Upsilon$. Hence Lemma \ref{lem:A-M-M} gives (\ref{formula:modified_E_CS}). Because (\ref{formula:dfn_E_by_lambda}) does not involve curvings, we have $\zeta' = \zeta = 0$.
\end{proof}

By the proposition above, we can perform the reduction of the strongly $G_0$-equivariant bundle gerbe with connection and curving $(\G^k, \nabla, f + \Upsilon)$ with respect to $\lambda_0$. Because $\su(2)$ is simple, the Lie algebra $\g_0 = \Omega\su(2)$ is such that $[\g_0, \g_0] = \g_0$. Thus, by Lemma \ref{lem:expression_beta}, the possible choice of a $\lambda$ required in our reduction is $\lambda = \lambda_0$ only.

Using $\lambda_0$ (and a connection on $q : M \to M/G_0$), we obtain a bundle gerbe with connection and curving $(\bar{\G}^k, \bar{\nabla}, \bar{f})$ over $M/G_0 \cong SU(2)$ by the reduction. Lemma \ref{lem:A-M-M} implies that the 3-curvature of $(\bar{\G}^k, \bar{\nabla}, \bar{f})$ is $-2\pi\im\chi$. The de Rham cohomology class of $\chi$ corresponds to $k \in \R \cong H^3(SU(2), \R)$. Hence we see that the Dixmier-Douady class $\delta(\bar{\G}^k)$ corresponds to $k \in \Z \cong H^3(SU(2), \Z)$.


\subsection{Chern-Simons pseudo $\T$-bundles}
\label{subsec:CS_pseudo_T_bundle}

We give here an example of the reduction of an equivariant pseudo $\T$-bundle. The pseudo $\T$-bundle we consider is related to the Chern-Simons bundle gerbe, and is interpreted as the the \textit{Chern-Simons line bundle for a 2-manifold with boundary} \cite{F1}. The reduction of the pseudo $\T$-bundle provides the data of the \textit{quasi-Hamiltonian space} given by Alekseev, Malkin and Meinrenken in the finite dimensional construction of the symplectic form on the moduli space of flat connections \cite{A-M-M}.

Notations in the previous subsection will be used without change.

\medskip

Let $\Sigma$ be a compact oriented 2-manifold whose boundary $\partial \Sigma$ is identified with $S^1$. We denote the space of connections on the trivial $SU(2)$-bundle over $\Sigma$ by $N = A^1(\Sigma, \su(2))$, and the gauge transformation group by $H = C^\infty(\Sigma, SU(2))$. As in the previous subsection, we consider the left action of $g \in H$ on $A \in N$ given by $A \mapsto g A g^{-1} + g dg^{-1}$. We denote by $H_0$ the subgroup of $H$ consisting of maps $g : \Sigma \to SU(2)$ which carry the base point on $\partial \Sigma = S^1$ to the identity element in $SU(2)$.

The restriction of a connection to the boundary gives a map $r : N \to M$. Similarly, we obtain homomorphisms $r : H \to G$ and $r : H_0 \to G_0$. (We use the same symbol $r$ for the maps induced by the restriction to the boundary.) 

\medskip

Recall the strongly $G$-equivariant bundle gerbe with connection and curving $(\G^k, \nabla, f + \Upsilon)$ over $M$. Through the homomorphism $r : H \to G$, we can regard the pull-back $r^*(\G^k, \nabla, f + \Upsilon)$ as a strongly $H$-equivariant bundle gerbe with connection and curving over $N$.

\begin{prop} 
There is an $H$-equivariant pseudo $\T$-bundle $(R, v)$ for $r^*\G^k$.
\end{prop}

\begin{proof}
First of all, we remark that the bundle gerbe $r^*\G^k = (r^*Y, r^*P, r^*s)$ is canonically identified with the lifting bundle gerbe associated with the $\Gamma$-bundle $r^*Y$ over $N$. Since $Y = M \times \Gamma$, we have $r^*Y = N \times \Gamma$. Let us consider the obvious $\T$-bundle $R = N \times \h\Gamma$ over $r^*Y$. We define a section $v : r^*Y^{[2]} \to \delta R^{\otimes -1} \otimes r^*P$ by $v((A, \gamma_1), (A, \gamma_2)) = (\h{\gamma}_2)^{\otimes -1} \otimes \h{\gamma}_1 \otimes (\h{\gamma}_1^{-1} \h{\gamma}_2 )$, where $\h{\gamma}_i \in \h{\Gamma}$ is an element such that $\varpi(\h{\gamma}_i) = \gamma_i$. This satisfies $\delta v = r^*s$, and we constructed a pseudo $\T$-bundle $(R, v)$ for $r^*\G^k$.

Next we introduce an action of $H$ on $R$. For $g \in H = C^\infty(\Sigma, SU(2))$, the \textit{Wess-Zumino term} \cite{F1} is an element $e^{W_\Sigma(g)} \in \h{\Gamma}^k$ such that $\varpi(e^{W_\Sigma(g)}) = r(g)$. The following formula is known:
$$
e^{W_\Sigma(gh)} =
\exp\frac{- k\im}{4\pi}\int_\Sigma \Tr (g^{-1}dg \wedge dh h^{-1}) \cdot
e^{W_\Sigma(g)} \cdot e^{W_\Sigma(h)},
$$
where the dot means the multiplication in $\h{\Gamma}^k$. Let $C : H \times N \to \h\Gamma^k$ be
$$
C(g, A) = 
\exp \frac{k\im}{4\pi} \int_{\Sigma} \Tr(A \wedge g^{-1}dg) \cdot 
e^{W_{\Sigma}(g)}.
$$
This obeys $C(gh, A) = C(g, hAh^{-1} + hdh^{-1}) \cdot C(h, A)$. So we can define an action of $H$ on $R = N \times \hGamma^k$ by $(A, \h{\gamma}) \mapsto (gAg^{-1} + gdg^{-1}, C(g, A) \cdot \h{\gamma})$. This makes $R$ into an $H$-equivariant principal $\T$-bundle over $r^*Y$. It is direct to see that the section $v$ is $H$-invariant.
\end{proof}

Let $\mathcal{L}$ be the Hermitian line bundle over $r^*Y = A^1(\Sigma, \su(2)) \times C^\infty(\Sigma, SU(2))$ associated to the $\T$-bundle $R \to r^*Y$. The line bundle $\mathcal{L}$ can be identified with the \textit{Chern-Simons line bundle for the 2-manifold $\Sigma$ with boundary} \cite{F1}. So we call $(R, v)$ the \textit{Chern-Simons pseudo $\T$-bundle}.

\begin{prop} 
There is an $H$-invariant connection $\eta$ on $(R, v)$.
\end{prop}

\begin{proof}
First, let $\vartheta \in \im A^1(N)$ be a 1-form given by
$$
\vartheta(A; \alpha) = - \frac{k\im}{4\pi} \int_\Sigma \Tr (A \wedge \alpha).
$$
Using the connection $\nu$ on the $\T$-bundle $\hGamma^k \to \Gamma$, we define a connection $\eta$ on the $\T$-bundle $R = N \times \h\Gamma^k$ by $\eta = \vartheta + \nu$. By computations, we can see $v^*(\delta \eta^{\otimes -1} \otimes \nabla) = 0$. Hence $\eta$ is a connection on the pseudo $\T$-bundle $(R, v)$.

Next we prove that $\eta$ is $H$-invariant. For $g \in H$, the difference $g^*\nu - \nu$ is a 1-form on $r^*Y = N \times \Gamma$, and is expressed as
$$
(g^*\nu - \nu)((A, \gamma); \alpha \oplus \gamma X) =
\frac{k\im}{4\pi} \int_\Sigma \Tr (\alpha \wedge g^{-1}dg)
$$
for a tangent vector $\alpha \oplus \gamma X \in T_{(A, \gamma)}(N \times \Gamma)$. The difference $g^*\vartheta - \vartheta$ is 
$$
(g^*\vartheta - \vartheta) (A; \alpha) =
- \frac{k\im}{4\pi} \int_\Sigma \Tr (\alpha \wedge g^{-1}dg).
$$
Therefore $\eta$ is an $H$-invariant connection on $R$.
\end{proof}

Using the $H$-equivariant pseudo $\T$-bundle with connection $((R, v), \eta)$ for the strongly $H$-equivariant bundle gerbe with connection $r^*(\G^k, \nabla, f+ \Upsilon)$, we define $\omega \in \im A^2(N)$ and $\mu \in A^0(N, \ha^\dagger)$ as in Lemma \ref{lem:data_pseudo_T_bundle_conn}, where we write $\ha = C^\infty(\Sigma, \su(2))$ for the Lie algebra of $H = C^\infty(\Sigma, SU(2))$. Then we have
\begin{align*}
\omega &= \sigma - r^*\Upsilon, \\
\langle X | \mu(A) \rangle  &= \frac{k\im}{2\pi} \int_{S^1} \Tr(F(A) X),
\end{align*}
where $\sigma = d\vartheta \in \im A^1(N)$ is given by
$$
\sigma(A; \alpha_1, \alpha_2) = 
- \frac{k\im}{2\pi} \int_\Sigma \Tr(\alpha_1 \wedge \alpha_2).
$$
As is well-known \cite{A-B}, $\frac{-1}{2\pi\im}\sigma$ is a (formal) symplectic form on the space of connections on the $SU(2)$-bundle over $\Sigma$. 

\medskip

Now we consider the reduction of the pseudo $\T$-bundle. Let us introduce the following $H$-invariant subspace of $N$:
$$
\mu^{-1}(0) = \{ A \in A^1(\Sigma, \su(2)) |\ F(A) = 0 \}.
$$
This is the space of flat connections. Because the action of $H_0$ on $N$ is free, so is the action on $\mu^{-1}(0)$. The quotient space $\mu^{-1}(0)/H_0$ is a finite dimensional smooth manifold \cite{A-M-M}. In fact, we can identify $\mu^{-1}(0)/H_0$ with $SU(2)^{2 g(\Sigma)}$, where $g(\Sigma)$ stands for the genus of $\Sigma$. As is clear, the restriction of $r^*(\G^k, \nabla, f + \Upsilon)$ to $\mu^{-1}(0)$ is a strongly $H_0$-equivariant bundle gerbe with connection and curving, and its obstruction to the reduction vanishes. We can naturally identify the bundle gerbe with connection and curving over $\mu^{-1}(0)/H_0$ obtained by the reduction of $r^*(\G^k, \nabla, f + \Upsilon)|_{\mu^{-1}(0)}$ with the pull-back of $(\bar{\G}^k, \bar{\nabla}, \bar{f})$ under the map $\bar{r} : \mu^{-1}(0)/H_0 \to M/G_0$ induced from $r : N \to M$. Thus, applying Proposition \ref{prop:reduction_pseudo_T_bundle}, we obtain a pseudo $\T$-bundle with connection $((\bar{R}, \bar{v}), \bar{\eta})$ for $\bar{r}^*(\bar{\G}^k, \bar{\nabla}, \bar{f})$ by the reduction of $((R, v), \eta)$.

\bigskip

Finally, we see that the Chern-Simons bundle gerbe and the Chern-Simons pseudo $\T$-bundle provide us data of a \textit{quasi-Hamiltonian space} \cite{A-M-M,G-S}.

Recall that the action of $G$ on $M$ induces the adjoint action of $G/G_0 \cong SU(2)$ on $M/G_0 \cong SU(2)$. The residual action makes $(\bar{\G}^k, \bar{\nabla}, \bar{f})$ into a strongly $G/G_0$-equivariant bundle gerbe with connection and curving. As associated data, we can recover the adjoint invariant closed 3-form $\chi \in A^0(SU(2))$ and the 1-form $e \in A^0(SU(2), \su(2)^*)$. They satisfy $\langle X | de \rangle = \iota_{X^*}\chi$ and $g^*e = \Ad_ge$ for $X \in \su(2)$ and $g \in SU(2)$.

Similarly, the action of $H$ on $N$ induces the action of $H/H_0$ on $\mu^{-1}(0)/H_0$. Note that $H/H_0 \cong G/G_0 \cong SU(2)$, and $\bar{r} : \mu^{-1}(0)/H_0 \to M/G_0$ is an equivariant map. The residual action makes $((\bar{R}, \bar{v}), \bar{\eta})$ into an $H/H_0$-equivariant pseudo $\T$-bundle with connection for $\bar{r}^*(\bar{\G}^k, \bar{\nabla}, \bar{f})$. Let $\bar{\omega} \in \im A^2(\mu^{-1}(0)/H_0)$ be the $H/H_0$-invariant 2-form defined by $\bar{\omega} = F(\bar{\eta}) - \bar{f}$. Note that $\bar{\omega}$ is the unique 2-form such that $q^*\bar{\omega} = \omega|_{\mu^{-1}(0)} = (\sigma - r^*\Upsilon)|_{\mu^{-1}(0)}$.

\begin{prop}[\cite{A-M-M}] 
We consider the space $\mu^{-1}(0)/H_0$ on which $SU(2)$ acts, the $SU(2)$-invariant 2-form $\bar{\omega}$ on $\mu^{-1}(0)/H_0$, and the $SU(2)$-equivariant map $\bar{r} : \mu^{-1}(0)/H_0 \to M/G_0 \cong SU(2)$. Then we have 
\begin{gather}
d \bar{\omega} = 2\pi\im \bar{r}^*\chi, \label{QH1_CS} \\
\iota_{X^*}\bar{\omega} = 
- 2\pi\im \langle X | \bar{r}^*e \rangle, \label{QH2_CS} \\
\begin{split}
& \{ V \in T_x(\mu^{-1}(0)/H_0) |\ 
\bar{\omega}(x; V, W) = 0 \ 
\textit{for} \ W \in T_x(\mu^{-1}(0)/H_0) \} \\
& =
\{ X^* \in T_x(\mu^{-1}(0)/H_0) |\ 
\langle X | \bar{r}^*e(x; W) \rangle = 0\ 
\textit{for} \ W \in T_x(\mu^{-1}(0)/H_0) \}.
\end{split} 
\end{gather}
Therefore $(\mu^{-1}(0)/H_0, \frac{-1}{2\pi\im}\bar{\omega})$ gives rise to a quasi-Hamiltonian $SU(2)$-space with its $SU(2)$-valued moment map $\bar{r}$.
\end{prop}

We remark that (\ref{QH1_CS}) and (\ref{QH2_CS}) follow from Lemma \ref{lem:data_pseudo_T_bundle_conn}.

The quasi-Hamiltonian space above is a one considered by Alekseev, Malkin and Meinrenken in the finite dimensional construction of the symplectic form on the moduli space of flat connections \cite{A-M-M}: let $\mathcal{C} \subset SU(2)$ be a conjugacy class of $SU(2)$, and $\mathcal{N}_\Sigma(\mathcal{C})$ the moduli space of flat connections on the trivial $SU(2)$-bundle over $\Sigma$ whose holonomies around $\partial \Sigma \cong S^1$ are in $\mathcal{C}$. We can describe the moduli space as $\mathcal{N}_\Sigma(\mathcal{C}) = \bar{r}^{-1}(\mathcal{C})/SU(2)$, and the symplectic structure is induced on $\bar{r}^{-1}(\mathcal{C})/SU(2)$ from the 2-form $\bar{\omega}$ via the quasi-Hamiltonian reduction.

\smallskip

As a framework including both the ordinary Hamiltonian spaces and the quasi-Hamiltonian spaces, the notion of \textit{Hamiltonian spaces for quasi-symplectic groupoids} is introduced by Xu \cite{X}. Its prequantization is also introduced by Laurent-Gengoux and Xu \cite{L-X}. It would be interesting to understand the reduction of strongly equivariant bundle gerbes and equivariant pseudo $\T$-bundles by using the framework.



\begin{flushleft}
Graduate school of Mathematical Sciences, University of Tokyo, \\
Komaba 3-8-1, Meguro-Ku, Tokyo, 153-8914 Japan. \\
e-mail: kgomi@ms.u-tokyo.ac.jp
\end{flushleft}


\begin{thebibliography}{99}

\bibitem{A-M-M}A. Alekseev, A. Malkin and E. Meinrenken, 
\textit{Lie group valued moment maps}.
J. Differential Geom. 48 (1998), no. 3, 445--495.
dg-ga/9707021.

\bibitem{A-B}M. F. Atiyah and R. Bott,
\textit{The Yang-Mills equations over Riemann surfaces}.
Philos. Trans. Roy. Soc. London Ser. A 308 (1983), no. 1505, 523--615.

\bibitem{B-X-Z}K. Behrend, P. Xu and B. Zhang,
\textit{Equivariant gerbes over compact simple Lie groups}.
C. R. Math. Acad. Sci. Paris 336 (2003), no. 3, 251--256. 
math.SG/0306183.

\bibitem{B-V}N. Berline and M. Vergne,
\textit{Classes caract\'eristiques \'equivariantes. Formule de localisation en cohomologie \'equivariante}.
C. R. Acad. Sci. Paris S\'erie. I Math. 295 (1982), no. 9, 539--541.

\bibitem{B-T}R. Bott and L. Tu,
\textit{Differential forms in algebraic topology}.
Graduate Texts in Mathematics, 82, Springer-Verlag, New York-Berlin, 1982.

\bibitem{Bry3}J-L. Brylinski,
\textit{Differentiable cohomology of gauge groups}.\\
math.DG/0011069.

\bibitem{Bry2}J-L. Brylinski,
\textit{Gerbes on complex reductive Lie groups}.
math.DG/0002158.

\bibitem{Bry1}J-L. Brylinski,
\textit{Loop spaces, Characteristic Classes and Geometric Quantization}.
Birkh$\ddot{\textrm{a}}$user Boston, Inc., Boston, MA, 1993.

\bibitem{Bry-M1}J-L. Brylinski and D. A. McLaughlin,
\textit{The geometry of degree-four characteristic classes and of line bundles on loop spaces. I}.
Duke Math. J. 75 (1994), no. 3, 603--638.

\bibitem{De}P. Deligne,
\textit{Th\'eorie de Hodge, III}.
Inst. Hautes \'Etudes Sci. Publ. Math. No. 44 (1974), 5--77.

\bibitem{De-F}P. Deligne and D. S. Freed,
\textit{Classical field theory}.
Quantum fields and strings: a course for mathematicians, Vol. 1 
(Princeton, NJ, 1996/1997), 137--225, Amer. Math. Soc., Providence, RI, 1999.


\bibitem{E-V}H. Esnault and E. Viehweg,
\textit{Deligne-Beilinson cohomology}.
Beilinson's conjectures on special values of $L$-functions, 43--91, 
Perspect. Math., 4, Academic Press, Boston, MA, 1988.

\bibitem{F1}D. S. Freed,
\textit{Classical Chern-Simons Theory. I},
Adv. Math. 113(1995), no.2, 237--303. 
hep-th/9206021.

\bibitem{F2}D. S. Freed,
\textit{Higher algebraic structures and quantization}.
Comm. Math. Phys. 159 (1994), no.2, 343--398. 
hep-th/9212115.

\bibitem{Gaw-R}K. Gawedzki and N. Reis,
\textit{WZW branes and gerbes}.
Rev. Math. Phys. 14 (2002), no. 12, 1281--1334.
hep-th/0205233.

\bibitem{Gi}J. Giraud,
\textit{Cohomologie  non-ab\'elienne}.
Grundl. 179, Springer Verlag (1971).

\bibitem{Go4}K. Gomi,
\textit{Connections and curvings on lifting bundle gerbes}. 
J. London Math. Soc. (2) 67 (2003), no. 2, 510--526. 
math.DG/0107175.

\bibitem{Go1}K. Gomi,
\textit{Equivariant smooth Deligne cohomology}.
Osaka J. Math. (to appear) math.DG/0307373.

\bibitem{Go3}K. Gomi,
\textit{Geometry of gerbes in Chern-Simons theory}.
PhD thesis, University of Tokyo.

\bibitem{Go2}K. Gomi,
\textit{Relationship between equivariant gerbes with connection and gerbes with connection over the quotient spaces}.
Commun. Contemp. Math. (to appear) 
math.DG/0308032.

\bibitem{G-S}V. W. Guillemin and S. Sternberg,
\textit{Supersymmetry and equivariant de Rham theory}.
Springer-Verlag, Berlin, 1999.

\bibitem{Ko}B. Kostant,
\textit{Quantization and unitary representations. I. Prequantization}.
Lectures in modern analysis and applications, III, pp. 87--208. 
Lecture Notes in Math., Vol. 170, Springer, Berlin, 1970.

\bibitem{L-X}C. Laurent-Gengoux and P. Xu,
\textit{Quantization of quasi-presymplectic groupoids and 
their Hamiltonian spaces}.
math.SG/0311154.

\bibitem{L-U}E. Lupercio and B. Uribe,
\textit{Deligne cohomology for orbifolds, discrete torsion and $B$-fields}.
Geometric and topological methods for quantum field theory 
(Villa de Leyva, 2001), 468--482, World Sci. Publishing, River Edge, NJ, 2003.
hep-th/0201184.

\bibitem{Ma-S}V. Mathai and D. Stevenson,
\textit{Chern character in twisted $K$-theory: 
equivariant and holomorphic cases}. 
Comm. Math. Phys. 236 (2003), no. 1, 161--186. 
hep-th/0201010.

\bibitem{Me}E. Meinrenken,
\textit{The basic gerbe over a compact simple Lie group}.
Enseign. Math. (2) 49 (2003), no. 3-4, 307--333.
math.DG/0209194.

\bibitem{Mu}M. K. Murray,
\textit{Bundle gerbes}.
J. London Math. Soc. (2) 54 (1996), no.2, 403-416.
dg-ga/9407015.

\bibitem{Mu-S}M. K. Murray and D. Stevenson,
\textit{Bundle gerbes: stable isomorphism and local theory}.
J. London Math. Soc. (2) 62 (2000), no.3, 925-937.
math.DG/9908135.

\bibitem{Mu-S2}M. K. Murray and D. Stevenson,
\textit{Higgs fields, bundle gerbes and string structures}.
Comm. Math. Phys. 243 (2003), 541--555.
math.DG/0106179.

\bibitem{P-S}A. Pressley and G. Segal,
\textit{Loop groups}. 
Oxford Mathematical Monographs. Oxford Science Publications. The Clarendon Press, Oxford University Press, New York, 1986. 

\bibitem{Se}G. Segal,
\textit{Classifying spaces and spectral sequences}.
Inst. Hautes \'Etudes Sci. Publ. Math. No. 34 1968 105--112.

\bibitem{We}A. Weil,
\textit{Introduction \`a l'\'etude des vari\'et\'es k\"ahl\'eriennes}.
Publications de l'Institut de Math\'ematique de l'Universit\'e de Nancago, VI.
Actualit\'es Sci. Ind. no. 1267. Hermann, Paris 1958.

\bibitem{Wi}E. Witten, 
\textit{Quantum field theory and the Jones polynomial}. 
Comm. Math. Phys. 121 (1989), no. 3, 351--399. 

\bibitem{X}P. Xu,
\textit{Momentum maps and Morita equivalence}.
math.SG/0307319.

\end{thebibliography}
\end{document}